\documentclass[11pt, eqno]{article}
\usepackage{bbm}
\usepackage{mathrsfs}
\usepackage{amsfonts}
\usepackage{amssymb}
\usepackage{graphicx}
\usepackage[all]{xy}

\usepackage{psfrag}
\usepackage{subfigure}
\usepackage{color}

\usepackage{amssymb,latexsym}
\usepackage{amsmath,latexsym}
\usepackage{amscd}

\newcommand{\K}{{\mathbb K}}
\newcommand{\F}{{\mathbb F}}
\newcommand{\N}{{\mathbb N}}
\newcommand{\C}{{\mathbb C}}
\newcommand{\R}{{\mathbb R}}
\newcommand{\Q}{{\mathbb Q}}
\newcommand{\Z}{{\mathbb Z}}
\newcommand{\CP}{{\mathbb CP}}

\newtheorem{theorem}{Theorem}[section]

\newtheorem{corollary}[theorem]{Corollary}

\newtheorem{remark}[theorem]{Remark}

\newtheorem{lemma}[theorem]{Lemma}

\newtheorem{proposition}[theorem]{Proposition}
\newtheorem{claim}[theorem]{Claim}

\begin{document}

\title{Splitting lemmas for the Finsler energy functional\\  on the space of $H^1$-curves}
\author{\\Guangcun Lu
\thanks{\endgraf
 Partially supported
by the NNSF   10971014 and 11271044 of China.}
\\
{\normalsize School of Mathematical Sciences, Beijing Normal University},\\
{\normalsize Laboratory of Mathematics
 and Complex Systems,  Ministry of
  Education},\\
  {\normalsize Beijing 100875, The People's Republic
 of China}\\
{\normalsize (gclu@bnu.edu.cn)}}
\date{
 Revised version,  April 28, 2016} \maketitle \vspace{-0.5in}

\abstract{We establish the splitting
lemmas (or generalized Morse lemmas) for the energy functionals of Finsler metrics on the
natural Hilbert manifolds of $H^1$-curves around a critical point or
a critical $\R^1$ orbit of a Finsler isometry invariant closed geodesic. They are the desired
generalization on Finsler manifolds of the corresponding Gromoll-Meyer's splitting lemmas
 on Riemannian manifolds (\cite{GM1, GM2}). As an application we extend to Finsler manifolds a result by
Grove and Tanaka \cite{GroTa78, Tan82} about the existence of infinitely many, geometrically distinct,
isometry invariant closed geodesics on a closed Riemannian manifold.} \vspace{-0.1in}
\medskip\vspace{12mm}

 \vspace{2mm}



\section{Introduction and results}
\setcounter{equation}{0}

Let $M$ be a smooth manifold of dimension $n$ endowed with a $C^k$
Finsler metric $F$, where $k\ge 5$ is an integer or $\infty$.  The
fundamental tensor $g^F$ of  $(M,F)$ is a Riemannian
metric on the pulled-back bundle
$\pi^\ast(T^\ast M)$ over $TM\setminus 0_{TM}$ defined by
$$
g^F(x,y): T_xM\times T_xM\to\R,\,(u,
v)\mapsto\frac{1}{2}\frac{\partial^2}{\partial s\partial
t}\left[F^2(x, y+ su+ tv)\right]\bigl|_{s=t=0}
$$
for any $(x,y)\in TM\setminus 0_{TM}$. It is of class $C^{k-2}$, and satisfies  $g^F(x,y)[y,y]=(F(x,y))^2$
and $g^F(x, \lambda y)=g^F(x,y)$ for $\lambda>0$.
  The length of a Lipschitz
continuous curve $\gamma:[a,b]\to M$ on $(M, F)$  is defined by $l_F(\gamma)=\int^b_aF(\gamma(t),\dot
\gamma(t))dt$.  A differentiable curve $\gamma=\gamma(t)$ is said to have
constant speed if $F(\gamma(t), \dot\gamma(t))$ is constant along
$\gamma$. Call a regular piecewise $C^k$ curve in $(M, F)$ a (Finslerian) geodesic
 if it minimizes the length between two sufficiently close points on
  the curve (\cite{BaChSh}).  From the
viewpoint of the calculus of variations the constant speed geodesics are the critical points
of the energy (or length) functional of the Finsler metric $F$ on a suitable Hilbert manifold.
Precisely, given a Riemannian metric $g$ on $M$ let
$W^{1,2}(I, M)$ denote  the space of absolutely continuous curves $\gamma$ from $I:=[0,1]$ to $M$ such that
$\int^1_0\langle\dot\gamma(t),\dot\gamma(t)\rangle dt<\infty$, where
$\langle u, v\rangle =g_x(u,v)$ for $u,v\in T_xM$. Then $W^{1,2}(I,
M)\subset C^0(I,M)$ and it
 is a Riemannian-Hilbert  manifold with the Riemannian structure on the tangent space
 $T_\gamma W^{1,2}(I, M)=W^{1,2}(\gamma^\ast TM)$ (consisting of
all $W^{1,2}$-sections of the pull-back bundle $\gamma^\ast TM\to I$)
  given by
\begin{equation}\label{e:1.1}
\langle\xi,\eta\rangle_1=\int^1_0\langle\xi(t),\eta(t)\rangle dt+
\int^1_0\langle\nabla^g_{\dot\gamma}\xi(t),\nabla^g_{\dot\gamma}\xi(t)\rangle dt
\end{equation}
 (using the $L^2$ covariant
derivative along $\gamma$ associated to the Levi-Civita connection $\nabla^g$ of the metric $g$).
 Let $\|\cdot\|_1=\sqrt{\langle\cdot,\cdot\rangle_1}$ be the induced norm.
 A smooth submanifold $N\subset M\times M$ determines
a Riemannian-Hilbert  submanifold of $W^{1,2}(I, M)$,
$$
\Lambda_N(M):=\{\gamma\in W^{1,2}(I, M)\,|\,(\gamma(0),\gamma(1))\in
N\}
$$
with tangent space $T_\gamma\Lambda_N(M)=W^{1,2}_N(\gamma^\ast TM)$ (consisting of
all $\xi\in W^{1,2}(\gamma^\ast TM)$  with $(\xi(0),\xi(1))\in T_{(\gamma(0),\gamma(1))}N$).
 The inclusion 
 $\Lambda_N(M)\hookrightarrow C^0_N(I,
M)=\{\gamma\in C^0(I,M)\,|\, (\gamma(0),\gamma(1))\\ \in N\}
$
 is a
homotopy equivalence (cf. \cite[Th.1.3]{Gro73I}).
Both $(W^{1,2}(I,M), \langle\cdot,\cdot\rangle_1)$ and
 $(\Lambda_N(M),\langle\cdot,\cdot\rangle_1)$ are complete Hilbert-Riemannian
manifolds if the metric $g$ is complete and the submanifold $N$ is
closed (as a subset).   See Appendix C of \cite{McSa} for  details.
The energy functional of the  metric $F$ given by
\begin{equation}\label{e:1.2}
{\cal L}:\Lambda_N(M)\to\R,\;\gamma\mapsto\int^1_0F^2(\gamma(t),
\dot\gamma(t))dt,
\end{equation}
is of class $C^{2-}$ by \cite{Me}, and cannot be
 twice differentiable at a regular curve
$\gamma\in\Lambda_N(M)$ with $N=\{p,q\}$ if
$F^2$ is not the square of the norm of a Riemannian metric along the curve,
where  a curve $\gamma\in\Lambda_N(M)$
is called \textsf{regular} if $\dot\gamma\ne 0$ a.e. in $[0,1]$ (\cite[page 271]{Ca}).
However, it  is of class $C^{k-2}$ on a dense open subset
 consisting of all regular curves in the Banach manifold
$$
{\cal
X}_N:=C^1_N(I, M)=\{\gamma\in C^1(I,M)\,|\, (\gamma(0), \gamma(1))\in
N\}.
$$
 A curve $\gamma\in\Lambda_N(M)$ is a constant (non-zero)
speed $F$-geodesic satisfying the boundary condition
\begin{equation}\label{e:1.3}
g^F(\gamma(0),\dot\gamma(0))[u,\dot\gamma(0)]=
g^F(\gamma(1),\dot\gamma(1))[v,\dot\gamma(1)]\quad\forall (u,v)\in
T_{(\gamma(0),\gamma(1))}N
\end{equation}
if and only if it is a (non constant) critical point of ${\cal L}$
(\cite{Me, KoKrVa, CaJaMa3}). In addition, suppose that $(M,F)$ is
forward (resp. backward) complete and that $N$ is  a closed
submanifold of $M\times M$ such that the first projection (resp. the
second projection) of $N$ to $M$ is compact; then ${\cal L}$
satisfies the Palais-Smale condition on $\Lambda_N(M)$ provided $(M,
g)$ is complete (cf. \cite{CaJaMa3}).

Because of the lack of $C^2$-smoothness for the functional ${\cal L}$ on $\Lambda_N(M)$,
the deep and very effective tool in the study of geodesics,  the famous splitting lemma of Gromoll-Meyer
\cite{GM1},  cannot be directly applied to  $\mathcal{L}$ as
in the Riemannian case \cite{GM2}. A workaround is to restrict ${\cal L}$ to
finite dimensional approximations of $\Lambda_N(M)$.
Our aim is to show that in the $H^1$ settings.
 When $N=\{p,q\}$
consists of two non-conjugate points (which means that all critical points
of $\mathcal{L}$ are nondegenerate),  Caponio, Javaloyes and Masiello \cite{CaJaMa1,CaJaMa2}
obtained the Morse inequalities without using
finite dimensional approximations recently. In the general case, the author \cite{Lu3}
developed an infinite dimensional method to obtain the shifting theorems of the critical groups of
critical points and critical orbits for the functional $\mathcal{L}$
 without involving  finite-dimensional approximations and any  Palais'
result in \cite{Pa}. In this paper we establish the generalized Morse lemmas (or splitting
lemmas) for the energy functional $\mathcal{L}$ in the $H^1$-topology around critical points or critical orbits
corresponding with  Gromoll-Meyer's splitting lemmas  in \cite{GM1, GM2}.

We only consider the following three cases:

\noindent{\bf Case 1.} $N=M_0\times M_1$, where $M_0$ and $M_1$ are
two  boundaryless and connected submanifolds of $M$. In this
case the boundary condition (\ref{e:1.3}) becomes
\begin{equation}\label{e:1.4}
\left\{\begin{array}{ll}
&g^F(\gamma(0),\dot\gamma(0))[u,\dot\gamma(0)]=0\quad\forall u\in
T_{\gamma(0)}M_0,\\
&g^F(\gamma(1),\dot\gamma(1))[v,\dot\gamma(1)]=0\quad\forall v\in
T_{\gamma(1)}M_1. \end{array}\right.
\end{equation}

\noindent{\bf Case 2.} $N=G(\mathbb{I}_g)$ the graph of an isometry
$\mathbb{I}_g$ on $(M, g)$,  in particular $G(id_M)=\triangle_M$.
The boundary condition (\ref{e:1.3}) becomes
\begin{equation}\label{e:1.5}
g^F(\gamma(0),\dot\gamma(0))[u,\dot\gamma(0)]=
g^F(\mathbb{I}_g(\gamma(0)),\dot\gamma(1))[\mathbb{I}_{g\ast}u,\dot\gamma(1)]
\quad\forall u\in T_{\gamma(0)}M.
\end{equation}
Furthermore, if  $F$ is $\mathbb{I}_{g}$-invariant, i.e.,
$F(\mathbb{I}_{g}(x),\mathbb{I}_{g\ast}(u))=F(x,u)\;\forall (x,u)\in
TM$, where $\mathbb{I}_{g\ast}:TM\to TM$ denotes the differential of
$\mathbb{I}_{g}$, we shall show above Theorem~\ref{th:1.5} that
(\ref{e:1.5}) becomes
$\mathbb{I}_{g\ast}(\dot\gamma(0))=\dot\gamma(1)$ and that $\gamma$
may be extended into an $\mathbb{I}_{g}$-invariant $F$-geodesic
$\gamma^\star:\R\to M$ via
\begin{equation}\label{e:1.6}
\gamma^\star(t)=\mathbb{I}^{[t]}_{g}(\gamma(t-[t]))\;\forall t\in\R,
\end{equation}
where $[t]$ denotes the greatest integer $\le t$,
 called the \textsf{corresponding (maximal)
$\mathbb{I}_{g}$-invariant $F$-geodesic} (determined by $\gamma$).
Here an $F$-geodesics $\alpha:\R\to M$ is said to be
\textsf{$\mathbb{I}_{g}$-invariant} if there exists a $s\ge 0$ such that
$\alpha(t+s)=\mathbb{I}_{g}(\alpha(t))\;\forall t\in\R$. (Even if
$F$ is the norm of a Riemannian metric different from $g$  this kind of case
seems not to be considered in Riemannian geometry.)

\noindent{\bf Case 3.} $N=G(\mathbb{I}_F)$ the graph of an isometry
$\mathbb{I}_F$ on $(M, F)$. 

 \vspace{2mm}

\noindent{\bf Notation.}  As in \cite{Lu3}, for a normed vector space $(E, \|\cdot\|)$
and $\delta>0$ we put ${\bf B}_\delta(E)=\{x\in
E\,|\,\|x\|<\delta\}$ and $\overline{{\bf B}_\delta(E)}=\{x\in
E\,|\,\|x\|\le\delta\}$ (in order to avoid confusion when
there are several spaces involved). Denote by $\mathscr{L}(E)$ the space of continuous linear operator
from $E$ to itself, and by $\mathscr{L}_s(E)$ the space of continuous linear self-adjoint operator
from $E$ to itself if $E$ is a Hilbert space. For a continuous symmetric bilinear form (or the
associated self-adjoint operator) $B$ on a Hilbert space we write
${\bf H}^-(B)$, ${\bf H}^0(B)$ and ${\bf H}^+(B)$ for the negative
definite, null and positive definite spaces of it.   $\K$ always
denotes an Abelian group without special statements.  \vspace{2mm}

\subsection{The case $N=M_0\times M_1$}\label{sec:1.1}

  Suppose that
$\gamma_0\in\Lambda_N(M)$ is a  critical point of ${\cal L}$ on
$\Lambda_N(M)$ with energy $c>0$.  Then it is a $C^k$-smooth
nonconstant $F$-geodesics with constant speed
$F(\gamma_0(t),\dot\gamma_0(t))\equiv \sqrt{c}>0$. Because of  technical reasons, we need to make
the following assumption on the Riemannian metric $g$ on $M$:
\begin{equation}\label{e:1.7}
\hbox{\textsf{ $M_0$ (resp. $M_1$) is totally geodesic near
$\gamma_0(0)$ (resp. $\gamma_0(1)$).}}
\end{equation}
 Let $\exp$ denote
the exponential map of $g$, and let $\|\cdot\|_1=\sqrt{\langle\cdot,\cdot\rangle_1}$
with $\langle\cdot,\cdot\rangle_1$ given by (\ref{e:1.1}).   Set
${\bf B}_{2\rho}(T_{\gamma_0}\Lambda_N(M)):= \{\xi\in
T_{\gamma_0}\Lambda_N(M)\,|\,\|\xi\|_{1}<2\rho\}$
for $\rho>0$. Since $\gamma_0^\ast TM$ is only $C^{k}$,
so is the map $\gamma_0^\ast TM\ni (t,v)\mapsto\exp_{\gamma_0(t)}v=\exp(\gamma_0(t),v)\in M$.
It follows from Lemma~\ref{lem:5.2} that for $\rho>0$ small
enough the map
\begin{equation}\label{e:1.8}
{\rm EXP}_{\gamma_0}: {\bf
B}_{2\rho}(T_{\gamma_0}\Lambda_N(M))\to\Lambda_N(M)
\end{equation}
given by ${\rm EXP}_{\gamma_0}(\xi)(t)=\exp_{\gamma_0(t)}(\xi(t))$,
is only  a $C^{k-3}$ coordinate chart around $\gamma_0$ on $\Lambda_N(M)$.
Therefore, $H^1(I, \gamma_0^\ast TM)$ is a $C^{k-3}$ Hilbert manifold, and
 $T_{\gamma_0}\Lambda_N(M)$ is a $C^{k-3}$ submanifold of $H^1(I, \gamma_0^\ast TM)$,
which implies that ${\cal L}\circ{\rm EXP}_{\gamma_0}$ is $C^{2-0}$ because $k\ge 5$.
Clearly, $0\in W^{1,2}_N(\gamma_0^\ast TM)$ is  a critical point of
 ${\cal L}\circ{\rm EXP}_{\gamma_0}$. Observe that the tangent space of  the Banach manifold
${\cal X}_N$ at $\gamma_0$ is
 $T_{\gamma_0}{\cal
X}_N=C^1_{TN}(\gamma_0^\ast TM)=\{\xi\in C^1(\gamma_0^\ast TM)\,|\,
(\xi(0), \xi(1))\in {TN} \}$ with usual $C^1$-norm, and that
${\bf
B}_{2\rho}(T_{\gamma_0}\Lambda_N(M))\cap T_{\gamma_0}{\cal X}_N$ is
an open neighborhood of $0$ in $T_{\gamma_0}{\cal X}_N$.
As above $T_{\gamma_0}{\cal
X}_N\ni\xi\mapsto{\rm EXP}_{\gamma_0}\xi\in {\cal
X}_N$ is also $C^{k-3}$.
 Let ${\cal A}$
be the restriction of the gradient of ${\cal L}\circ{\rm
EXP}_{\gamma_0}$ to ${\bf B}_{2\rho}(T_{\gamma_0}\Lambda_N(M))\cap
T_{\gamma_0}{\cal X}_N$. Since $\gamma_0$ is $C^k$ and $\dot{\gamma}_0(t)\ne 0$, $\forall t$,
${\cal A}$ is actually a $C^{k-3}$ map from a small neighborhood of $0\in
T_{\gamma_0}{\cal X}_N$ to $T_{\gamma_0}{\cal X}_N$, which can be derived from (\ref{e:2.29}) and
Lemma~\ref{lem:2.3}, moreover
${\cal A}(0)=\nabla{\cal L}(\gamma_0)|_{T_{\gamma_0}{\cal X}_N}$ and
$$
\langle d{\cal A}(0)[\xi],\eta\rangle_{1}=d^2{\cal L}|_{{\cal
X}_N}(\gamma_0)[\xi,\eta]\quad\forall\xi,\eta\in T_{\gamma_0}{\cal
X}_N.
$$
This extends results obtained in \cite{CaJaMa1} for the case $N=\{p,q\}$.
The key point here is that the continuous symmetric bilinear form $d^2{\cal L}|_{{\cal
X}_N}(\gamma_0)$ can be extended into a continuous symmetric bilinear form on
$T_{\gamma_0}\Lambda_N(M)$, also denoted by $d^2{\cal L}|_{{\cal
X}_N}(\gamma_0)$.  The self-adjoint operator associated to the
latter is Fredholm, and has finite dimensional negative definite and
null spaces ${\bf H}^-(d^2{\cal L}|_{{\cal X}_N}(\gamma_0))$ and
${\bf H}^0(d^2{\cal L}|_{{\cal X}_N}(\gamma_0))$, which are actually
contained in $T_{\gamma_0}{\cal X}_N$.  There exists an orthogonal
decomposition
\begin{equation}\label{e:1.9}
T_{\gamma_0}\Lambda_N(M)={\bf H}^-(d^2{\cal L}|_{{\cal
X}_N}(\gamma_0))\oplus {\bf H}^0(d^2{\cal L}|_{{\cal
X}_N}(\gamma_0))\oplus{\bf H}^+(d^2{\cal L}|_{{\cal
X}_N}(\gamma_0)),
\end{equation}
which  induces a (topological) direct sum decomposition of Banach
spaces
$$
T_{\gamma_0}{\cal X}_N={\bf H}^-(d^2{\cal L}|_{{\cal
X}_N}(\gamma_0))\dot{+} {\bf H}^0(d^2{\cal L}|_{{\cal
X}_N}(\gamma_0))\dot{+}\bigr({\bf H}^+(d^2{\cal L}|_{{\cal
X}_N}(\gamma_0))\cap T_{\gamma_0}{\cal X}_N\bigl).
$$
Using the implicit function theorem we get $\delta\in (0, 2\rho]$
and a unique $C^{k-3}$-map
\begin{equation}\label{e:1.10}
 h: {\bf B}_\delta\bigl({\bf H}^0(d^2{\cal L}|_{{\cal
X}_N}(\gamma_0))\bigr)\to {\bf H}^-(d^2{\cal L}|_{{\cal
X}_N}(\gamma_0))\dot{+}\bigl({\bf H}^+(d^2{\cal L}|_{{\cal
X}_N}(\gamma_0))\cap T_{\gamma_0}{\cal X}_N\bigr)
\end{equation}
such that $h(0)=0$, $dh(0)=0$ and
\begin{equation}\label{e:1.11}
(I-P^0){\cal A}(\xi+
h(\xi))=0\;\forall\xi\in {\bf B}_\delta\bigl({\bf H}^0(d^2{\cal
L}|_{{\cal X}_N}(\gamma_0))\bigr),
\end{equation}
where
$P^\star:T_{\gamma_0}\Lambda_N(M)\to {\bf H}^\star(d^2{\cal
L}|_{{\cal X}_N}(\gamma_0))$, $\star=-,0,+$, are the orthogonal
projections according to the decomposition (\ref{e:1.9}). Define
${\cal L}^\circ:{\bf B}_\delta\bigl({\bf H}^0(d^2{\cal L}|_{{\cal
X}_N}(\gamma_0))\bigr)\to\R$ by
\begin{equation}\label{e:1.12}
{\cal L}^\circ(\xi)= {\cal L}\circ{\rm EXP}_{\gamma_0}(\xi+  h(\xi)).
\end{equation}
It is  $C^{k-3}$ since
  $T_{\gamma_0}{\cal X}_N\ni\xi\mapsto{\rm EXP}_{\gamma_0}\xi\in {\cal
X}_N$ is also $C^{k-3}$;
moreover ${\cal L}^\circ$ has a
critical point $0$, and $d^2{\cal L}^\circ(0)=0$.
 Call
$$
m^-({\cal L},\gamma_0):=\dim {\bf H}^-(d^2{\cal L}|_{{\cal
X}_N}(\gamma_0))\quad\hbox{and}\quad m^0({\cal L},\gamma_0):=\dim {\bf
H}^0(d^2{\cal L}|_{{\cal X}_N}(\gamma_0))
$$
the {\bf Morse index} and the {\bf nullity} of $\gamma_0$, respectively.
(Actually, they  do not depend on the choice of the Riemannian
metric $g$).  Here is our first  \textsf{splitting lemma}.

\begin{theorem}\label{th:1.1}
Under the  above notation, if $\dim M>1$ for sufficiently small
  $\delta>0$  there exists  an origin-preserving
homeomorphism $\psi$ from ${\bf B}_\delta(T_{\gamma_0}\Lambda_N(M))$
to an open neighborhood of $0$ in $T_{\gamma_0}\Lambda_N(M)$
 such that for all $\xi\in{\bf B}_\delta(T_{\gamma_0}\Lambda_N(M))$,
$$
{\cal L}\circ{\rm EXP}_{\gamma_0}\circ\psi(\xi)=\|P^+\xi\|_1^2-
\|P^-\xi\|_1^2 + {\cal L}^\circ(P^0\xi).
$$
 Moreover, $\psi\left((P^-+P^0){\bf B}_\delta(T_{\gamma_0}\Lambda_N(M))\right)$
is contained in $T_{\gamma_0}{\cal X}_N$,
and $\psi$ is also a homeomorphism from $(P^-+P^0){\bf
B}_\delta(T_{\gamma_0}\Lambda_N(M))$ onto $\psi\left((P^-+P^0){\bf
B}_\delta(T_{\gamma_0}\Lambda_N(M))\right)$ even if the topology on
the latter is taken as the induced one by $T_{\gamma_0}{\cal X}_N$;
 $\psi$ restricts to a $C^{k-3}$ embedding from $P^0{\bf B}_\delta(T_{\gamma_0}\Lambda_N(M))$
to $T_{\gamma_0}{\cal X}_N$ (and so to $T_{\gamma_0}\Lambda_N(M)$) because $\psi(\xi)=\xi+h(\xi)$
 for $\xi\in P^0{\bf B}_\delta(T_{\gamma_0}\Lambda_N(M))$.
\end{theorem}

 If the critical point $\gamma_0$ is also isolated, then
 ${\cal L}\circ{\rm EXP}_{\gamma_0}$ and ${\cal L}^\circ$ have
 isolated critical points $0$. Let $C_\ast({\cal L},
\gamma_0;\K)$ (resp. $C_{\ast}({\cal L}^\circ, 0;\K)$) denote  the
critical group of the functional ${\cal L}$ (resp. ${\cal L}^\circ$)
at $\gamma_0$ (resp. $0$) with the coefficient group $\K$.
By Theorem~\ref{th:1.1} $\mathscr{N}:={\rm EXP}_{\gamma_0}\circ\psi\bigl(P^0{\bf B}_\delta(T_{\gamma_0}\Lambda_N(M))\bigr)$
is a $C^{k-3}$ (embedded) submanifold of ${\cal X}_N$ (and so $\Lambda_N(M)$)
containing $\gamma_0$ as an interior point. It is uniquely determined by $F$ and $g$
in the sense of germs (since our $h$ and $\psi$ may be explicitly determined by data only depending
on $g$ and $F$), and is called a \textsf{characteristic submanifold} of $\Lambda_N(M)$
for ${\cal L}$ at $\gamma_0$. Clearly, $C_{\ast}({\cal L}|_{\mathscr{N}},\gamma_0;\K)\cong
C_{\ast}({\cal L}^\circ, 0;\K)$. $\mathscr{H}^0_\ast({\cal L}, \gamma_0;\K):=C_{\ast}({\cal L}|_{\mathscr{N}},\gamma_0;\K)$ is often called the \textsf{characteristic invariant} of ${\cal L}$
at $\gamma_0$. We have the following \textsf{shifting theorem}.

\begin{theorem}\label{th:1.2}
$C_q({\cal L}, \gamma_0;\K)\cong C_{q-m^-({\cal L},\gamma_0)}({\cal L}^\circ,
0;\K)\cong \mathscr{H}^0_{q-m^-({\cal L},\gamma_0)}({\cal L},
\gamma_0;\K)\;\forall q=0,1,\cdots$.
\end{theorem}

\begin{remark}\label{rm:1.3}
{\rm
Any Riemannian metric $g'$ on $M$  can always be perturbed near $\gamma_0(0)$
and $\gamma_0(1)$ (with local charts around them) to yield a metric $g$ satisfying
(\ref{e:1.7}). Since $g$ and $g'$ induce the same Hilbert
manifold structures on $\Lambda_N(M)$ and equivalent Hilbert-Riemannian structures
on $T\Lambda_N(M)$ the assumption (\ref{e:1.7}) has almost no effect for application
ranges of our theory. It is also obvious that Theorem~\ref{th:1.1} is only related to the behavior
of $g$ in a  neighborhood of $\gamma_0([0,1])\subset M$.}
\end{remark}

\begin{theorem}\label{th:1.4}
Under the assumptions of Theorem~\ref{th:1.1}, let
$\tilde\Lambda$ be a Hilbert submanifold of $\Lambda_N(M)$ containing $\gamma_0$
and satisfying the following conditions:
\begin{description}
\item[(i)]  $\tilde{\cal X}_N:={\cal X}_N\cap\tilde\Lambda$ is a Banach submanifold of
${\cal X}_N$, and is also dense in $\tilde\Lambda$;

\item[(ii)]  the chart in (\ref{e:1.8}) restricts to a chart around $\gamma_0$
on $\tilde\Lambda$, $\widetilde{\rm EXP}_{\gamma_0}: {\bf
B}_{2\rho}(T_{\gamma_0}\tilde\Lambda)\to\tilde\Lambda$, and
a chart around $\gamma_0$
on $\tilde{\cal X}_N$, $\widetilde{\rm EXP}^X_{\gamma_0}: {\bf
B}_{2\rho}(T_{\gamma_0}\tilde\Lambda)\cap T_{\gamma_0}\tilde{\cal X}_N\to\tilde{\cal X}_N$;

\item[(iii)] $\nabla({\cal L}\circ\widetilde{\rm EXP}_{\gamma_0})(\xi)=\nabla({\cal L}\circ{\rm EXP}_{\gamma_0})(\xi)$ for any $\xi\in {\bf
B}_{2\rho}(T_{\gamma_0}\tilde\Lambda)$, where the inner product on $T_{\gamma_0}\tilde\Lambda$
is also given by (\ref{e:1.1});

\item[(iv)] ${\bf H}^0\bigl(d^2({\cal L}|_{{\cal X}_N})(\gamma_0)\bigr)\subset T_{\gamma_0}\tilde\Lambda$.
\end{description}
Then (since $d^2({\cal L}|_{\tilde{\cal X}_N})(\gamma_0)[u,v]
=d^2({\cal L}\circ\widetilde{\rm EXP}^X_{\gamma_0})(0)[u,v]=
d^2({\cal L}\circ{\rm EXP}^X_{\gamma_0})(0)[u,v]\linebreak=d^2({\cal L}|_{{\cal X}_N})(\gamma_0)[u,v]$
for $u,v\in T_{\gamma_0}\tilde{\cal X}_N$ implies
${\bf H}^\ast\bigl(d^2({\cal L}|_{\tilde{\cal X}_N})(\gamma_0)\bigr)
\subset{\bf H}^\ast\bigl(d^2({\cal L}|_{{\cal X}_N})(\gamma_0)\bigr)$, $\ast=0,+,-$), $\psi$ in Theorem~\ref{th:1.1} restricts to
 an origin-preserving
homeomorphism $\tilde\psi$ from ${\bf B}_\delta(T_{\gamma_0}\tilde\Lambda)$
to an open neighborhood of $0$ in $T_{\gamma_0}\tilde\Lambda$
 such that for all $\xi\in{\bf B}_\delta(T_{\gamma_0}\tilde\Lambda)$,
$$
{\cal L}\circ\widetilde{\rm EXP}_{\gamma_0}\circ\tilde\psi(\xi)=\|\tilde{P}^+\xi\|_1^2-
\|\tilde{P}^-\xi\|_1^2 + {\cal L}^\circ(\tilde{P}^0\xi),
$$
where $\tilde{P}^\ast$ are orthogonal projections from
$$
T_{\gamma_0}\tilde\Lambda={\bf H}^-(d^2{\cal L}|_{\tilde{\cal
X}_N}(\gamma_0))\oplus {\bf H}^0(d^2{\cal L}|_{\tilde{\cal
X}_N}(\gamma_0))\oplus{\bf H}^+(d^2{\cal L}|_{\tilde{\cal
X}_N}(\gamma_0))
$$
onto ${\bf H}^\ast(d^2{\cal L}|_{\tilde{\cal
X}_N}(\gamma_0))$, $\ast=0,-,+$, and we actually have
$\tilde{P}^\star\xi={P}^\star\xi$, $\star=-,0,+$.
 Moreover the inclusion $\tilde\Lambda\hookrightarrow\Lambda_N(M)$ induces isomorphisms
\begin{description}
\item[(a)] $\mathscr{H}^0_\ast({\cal L}|_{\tilde\Lambda},\gamma_0;\K)\cong \mathscr{H}^0_\ast({\cal L},
\gamma_0;\K)$ if $m^0({\cal L}|_{\tilde\Lambda},\gamma_0))=m^0({\cal L},\gamma_0))$;

\item[(b)] $C_\ast({\cal L}|_{\tilde\Lambda}, \gamma_0;\K)\cong C_\ast({\cal L}, \gamma_0;\K)$ if
$m^\star({\cal L}|_{\tilde\Lambda},\gamma_0))=m^\star({\cal L},\gamma_0))$, $\star=0,-$.
\end{description}
\end{theorem}

This result may be, in some sense, viewed as \textsf{heritability} of the splitting lemma.
In the Riemannian case,   ${\cal X}_N$ and $\tilde{\cal X}_N$
are chosen as $\Lambda_N(M)$ and $\tilde\Lambda$, respectively; the conditions (i)-(ii) are not needed;
(b) is  obvious and  the corresponding claim with (a) is a direct consequence of Lemma~7 in \cite{GM2}.
Such a result is necessary in many applications.

\subsection{The case $N=G(\mathbb{I}_g)$}\label{sec:1.2}

 There exists, as showed in
\cite[Lem.2.2]{Gro74}, a natural continuous $\R$-action by
Riemannian-isometries on $\Lambda_{G(\mathbb{I}_g)}(M)$ induced by
$g$, $\mu:\Lambda_{G(\mathbb{I}_g)}(M)\times\R\to
\Lambda_{G(\mathbb{I}_g)}(M)$ defined by
\begin{equation}\label{e:1.13}
\mu(\gamma, s)(t)=\left\{\begin{array}{cc} \gamma(t+s)\,&\,{\rm
for}\;
0\le t\le 1-s\\
\mathbb{I}_g(\gamma(t+s-1))\,&\,{\rm for}\; 1-s\le t\le 1
\end{array}\right.
\end{equation}
if $s\in [0, 1]$, and $\mu(\gamma,s)=\mu(\mathbb{I}_g^m\circ\gamma,
s-m))$ when $s\in [m, m+1]$. Correspondingly, we have also a
continuous $\R$-action by Hilbert-Riemannian bundle isomorphisms on
$T\Lambda_{G(\mathbb{I}_g)}(M)$,
$\hat\mu:T\Lambda_{G(\mathbb{I}_g)}(M)\times\R\to
T\Lambda_{G(\mathbb{I}_g)}(M)$ defined by
\begin{equation}\label{e:1.14}
\hat\mu(\xi, s)(t)=\left\{\begin{array}{cc} \xi(t+s)\,&\,{\rm for}\;
0\le t\le 1-s\\
\mathbb{I}_{g\ast}(\xi(t+s-1))\,&\,{\rm for}\; 1-s\le t\le 1
\end{array}\right.
\end{equation}
if $s\in [0, 1]$, and
$\hat\mu(\xi,s)=\mu(\mathbb{I}_{g\ast}^m\circ\xi, s-m))$ when $s\in
[m, m+1]$.

The isotropy group $\R_\gamma=\{s\in\R\,|\,\mu(\gamma,s)=\gamma\}$
at a point $\gamma\in \Lambda_{G(\mathbb{I}_g)}(M)$ is either
$\{0\}$, infinite cyclic or $\R$. The orbit space
$\Lambda_{G(\mathbb{I}_g)}(M)/\R$ is not Hausdorff. Moreover, the
orbit $\mu(\gamma,\R)$ of each $\gamma\in
C^r(I,M)\cap\Lambda_{G(\mathbb{I}_g)}(M)\setminus{\rm
Fix}(\mathbb{I}_g)$ with
$\dot\gamma(1)=\mathbb{I}_{g\ast}(\dot\gamma(0))$ is an immersed
$C^{r-1}$ submanifold of $\Lambda_{G(\mathbb{I}_g)}(M)$, where $r\ge 1$ is an integer or $\infty$. Note that
the functional ${\cal L}$ is not invariant under the above action
without further assumptions. When $N$ in
Theorem~\ref{th:1.1} is replaced by $G(\mathbb{I}_g)$ the
corresponding conclusions also hold, which can be seen from the proofs
of Theorem~\ref{th:1.1} and  Theorem~\ref{th:1.7}.

\begin{theorem}\label{th:1.5}
Let $F$ be a $C^k$ Finsler metric on a smooth manifold $M$ of dimension more than one. Suppose that
$\gamma_0$ is a nonconstant critical point  of ${\cal L}$ on
$\Lambda_{G(\mathbb{I}_g)}(M)$.  Then for sufficiently small
  $\delta>0$  there exists  an origin-preserving
homeomorphism $\psi$ from ${\bf
B}_\delta(T_{\gamma_0}\Lambda_{G(\mathbb{I}_g)}(M))$ to an open
neighborhood of $0$ in $T_{\gamma_0}\Lambda_{G(\mathbb{I}_g)}(M))$
 such that for all $\xi\in{\bf B}_\delta(T_{\gamma_0}\Lambda_{G(\mathbb{I}_g)}(M))$,
$$
{\cal L}\circ{\rm EXP}_{\gamma_0}\circ\psi(\xi)=\|P^+\xi\|_1^2-
\|P^-\xi\|_1^2 + {\cal L}^\circ(P^0\xi).
$$
 Moreover,  $\psi$ restricts to a $C^{k-3}$ embedding from $P^0{\bf
B}_\delta(T_{\gamma_0}\Lambda_{G(\mathbb{I}_g)}(M))$
to $T_{\gamma_0}{\cal X}_{G(\mathbb{I}_g)}$ (and so to $T_{\gamma_0}\Lambda_{G(\mathbb{I}_g)}(M)$);
$\psi\left((P^-+P^0){\bf
B}_\delta(T_{\gamma_0}\Lambda_{G(\mathbb{I}_g)}(M))\right)$ is
contained in $T_{\gamma_0}{\cal X}_{G(\mathbb{I}_g)}$, and $\psi$ is
also a homeomorphism from $(P^-+P^0){\bf
B}_\delta(T_{\gamma_0}\Lambda_{G(\mathbb{I}_g)}(M))$ onto
$\psi\left((P^-+P^0){\bf
B}_\delta(T_{\gamma_0}\Lambda_{G(\mathbb{I}_g)}(M))\right)$ even if
the topology on the latter is taken as the induced one by
$T_{\gamma_0}{\cal X}_{G(\mathbb{I}_g)}$.
\end{theorem}

If the critical point $\gamma_0$ is isolated, a shifting theorem
follows from Theorem~\ref{th:1.5} directly; and
a corresponding result with Theorem~\ref{th:1.4} is also easily given.

\textsf{From now on we suppose that $F$ is
$\mathbb{I}_g$-invariant}. In this case it is convenient to consider the functional ${\cal L}$
on the Hilbert-Riemannian manifold
$$
\Lambda(M,\mathbb{I}_g)=\{\gamma\in
W^{1,2}_{loc}(\R,M)\,|\,\gamma(t+1)=\mathbb{I}_g(\gamma(t))\;\forall
t\}
$$
with the natural Riemannian metric given by (\ref{e:1.1}) because
there exists a natural Hilbert-Riemannian isometry
$\Lambda_{G(\mathbb{I}_g)}(M)\ni\gamma\to\gamma^\star\in\Lambda(M,\mathbb{I}_g)$, where
$\gamma^\star$ is defined by (\ref{e:1.6}), whose
inverse is given by the restriction map
$\Lambda(M,\mathbb{I}_g)\ni\gamma\mapsto\gamma|_I\in
\Lambda_{G(\mathbb{I}_g)}(M)$. Under this correspondence the action in (\ref{e:1.13})
corresponds to  a continuous $\R$-action on $\Lambda(M,\mathbb{I}_g)$
by isometries, $\mathcal{T}: \Lambda(M, \mathbb{I}_g)\times\R\to \Lambda(M, \mathbb{I}_g)$ by
\begin{equation}\label{e:1.15}
\mathcal{T}(\gamma, s)(t)={\bf T}_s(\gamma)(t)=\gamma(t+ s)\quad\hbox{for any}\;t, s\in\R\;\hbox{and}\;\gamma\in\Lambda(M,\mathbb{I}_g).
\end{equation}
We also write ${\bf T}_s(\gamma)$ as $s\cdot\gamma$ sometimes.
(Clearly, $(s\cdot\gamma^\star)|_{[0,1]}=\mu(\gamma,s)$ for any $s\in\R$
if $\gamma\in \Lambda_{G(\mathbb{I}_g)}(M)$  corresponds to
 $\gamma^\star\in\Lambda(M,\mathbb{I}_g)$ as above.)  Note that if the isometry $\mathbb{I}_g$ is of finite order
$p\in\N$, i.e. $\mathbb{I}_g^p=id_M$ then the $\R$-action in
(\ref{e:1.15}) induces a continuous $\R/p\Z$-action by isometries on
$\Lambda(M,\mathbb{I}_g)$:
\begin{eqnarray}\label{e:1.16}
[s]_p\cdot \gamma(t)=\gamma(s+t), \quad \forall [s]_p\in \R/p\Z,
         \; \gamma\in \Lambda(M,\mathbb{I}_g).
  \end{eqnarray}
Clearly, the critical points of the functional $\mathcal{L}$ on
$\Lambda(M,\mathbb{I}_g)$  consist of  all
$\mathbb{I}_g$-invariant $F$-geodesics plus the constant maps in the
fixed point set ${\rm Fix}(\mathbb{I}_g)$ of $\mathbb{I}_g$. These
critical points are of class $C^k$, and in particular sit in the
Banach manifold
$$
{\cal X}={\cal X}(M, \mathbb{I}_g)=\{\gamma\in
C^1(\R,M)\,|\,\gamma(t+1)=\mathbb{I}_g(\gamma(t))\;\forall t\}.
$$
A nonconstant curve $\gamma\in \Lambda(M, \mathbb{I}_g)$ can be a non-isolated critical point
of  ${\cal L}$. For each critical point $\gamma$ of ${\cal L}$,
 corresponding to Proposition~2.3 in \cite{Gro74}  we have also:
\begin{description}
\item[(i)] $\R_\gamma=\{0\}$ if and only if the corresponding $\mathbb{I}_g$-invariant
 $F$-geodesic on $M$ is non-closed.
\item[(ii)] $\R_\gamma$ is infinite cyclic with generator $s$ if and only if the
corresponding $\mathbb{I}_g$-invariant $F$-geodesic is closed and
the prime period of the  geodesic is $s$.
\item[(iii)] $\R_\gamma=\R$  if and only $\gamma$ is a constant map to ${\rm Fix}(\mathbb{I}_g)$.
\end{description}
So the $\R$-orbit of a critical point $\gamma$ of ${\cal L}$ on $\Lambda(M, \mathbb{I}_g)$  is
either a point, an embedded $S^1$ or a $1-1$ immersed image of $\R$
(with constant speed)  according to the corresponding maximal
$\mathbb{I}_g$-invariant $F$-geodesic  being a fixed point, a closed
geodesic or a non-closed geodesic.
 The critical points in
the same $\R$-orbit correspond to the same $\mathbb{I}_g$-invariant
$F$-geodesic with only different initial point. A non-closed $\mathbb{I}_g$-invariant $F$-geodesic
 is in $1-1$ correspondence with a non-closed $\R$-orbit of critical points in
$\Lambda(M,\mathbb{I}_g)$; and a closed $\mathbb{I}_g$-invariant
$F$-geodesic is in $1-1$ correspondence with a ``tower" of closed
$\R$-orbits of critical points in $\Lambda(M,\mathbb{I}_g)$. (Note
that each orbit of the tower corresponds to a ``covering" of an
underlying prime closed $\mathbb{I}_g$-invariant geodesic.)

Carefully checking the proof of \cite[Th.2.4]{Gro74} and \cite[Th.4.1]{Gro85}
  it is not hard to see that the following result holds true.

\begin{proposition}\label{prop:1.5}
Let $\mathbb{I}_g$ be proper in the sense that the displacement
function for it, $\delta_{\mathbb{I}_g}:M\to\R,\;x\mapsto d_g(x,
\mathbb{I}_g(x))$ is proper. Suppose that $F$ is
$\mathbb{I}_g$-invariant and that there exists a non-closed (i.e.,
non-periodic) $\mathbb{I}_g$-invariant $F$-geodesic $\gamma:\R\to
M$.  Then in the closure $\overline{\gamma(\R)}$ there are
uncountably many $\mathbb{I}_g$-invariant $F$-geodesics
$\gamma_\tau:\R\to M$ with
$F(\gamma_\tau,\dot\gamma_\tau)=F(\gamma,\dot\gamma)$.
\end{proposition}



So if the isometry $\mathbb{I}_g$ is
  proper,  in the studies of multiplicity
of $\mathbb{I}_g$-invariant $F$-geodesics, we can assume that there exist
finitely many (geometrically different) $\mathbb{I}_g$-invariant geodesics
(which must all be closed).

Let $\gamma_0\in \Lambda(M,\mathbb{I}_g)$ be a closed
$\mathbb{I}_g$-invariant $F$-geodesic with ${\cal L}(\gamma_0)=c>0$
(which implies that $\R_{\gamma_0}$ is an infinite cyclic subgroup of $\R$ with generator $p>0$).
Then the orbit ${\cal
O}:=\R\cdot\gamma_0$ is an embedded circle $S^1(p):=\R/p\Z$. (The
$\R$-action reduces to $S^1(p)$-action on it).  We assume that ${\cal
O}$ is an isolated critical orbit below.
 It is an $\R$-invariant
compact connected $C^{k-1}$ submanifold of $\Lambda(M,\mathbb{I}_g)$ (cf. \cite[page 499]{GM2}).
Actually, ${\cal O}\subset C^k(\R, M)\cap \Lambda(M,\mathbb{I}_g)$ by Proposition~3.1 in
\cite{Lu3}, and it is a $C^{k-1}$ submanifold of the Banach manifold ${\cal
X}={\cal X}(\mathbb{I}_g,M)$ by Proposition~\ref{prop:5.1}.

 Let $\pi:N{\cal O}\to{\cal O}$ be the normal bundle
of ${\cal O}$ in $\Lambda(M,\mathbb{I}_g)$. This is a $C^{k-2}$
Hilbert vector bundle over ${\cal O}$ (because $T{\cal O}$ is a $C^{k-2}$ subbundle of
    $T_{\cal O}\Lambda(M,\mathbb{I}_g)$), and $XN{\cal O}:=T_{\cal
O}{\cal X}\cap N{\cal O}$ is a $C^{k-2}$ Banach vector subbundle of
$T_{\cal O}{\cal X}$ by Proposition~\ref{prop:5.1}. For
$\varepsilon>0$ we define
  $N{\cal O}(\varepsilon):=\{(\gamma,v)\in N{\cal
 O}\,|\,\|v\|_1<\varepsilon\}$,
$$
 N{\cal O}(\varepsilon)^X:=N{\cal O}(\varepsilon)\cap T_{\cal O}{\cal X}\quad\hbox{and}\quad
 XN{\cal O}(\varepsilon):= \{(\gamma,v)\in XN{\cal
O}\,|\,\|v\|_{C^1}<\varepsilon\}
$$
as open subsets of $T_{\cal O}{\cal X}$. Clearly, $XN{\cal
O}(\varepsilon)\subset N{\cal O}(\varepsilon)^X\subset N{\cal
O}(\varepsilon)$. Note that there exists a  natural induced
$\R$-actions on these bundles.
 For sufficiently small  $\varepsilon>0$, as above we derive from
  Lemma~\ref{lem:5.2} that the map
\begin{equation}\label{e:1.17}
{\rm EXP}:T\Lambda(M,\mathbb{I}_g)(\varepsilon)=\{(\gamma,v)\in
T\Lambda(M,\mathbb{I}_g)\,|\,\|v\|_1\}<\varepsilon\}\to
\Lambda(M,\mathbb{I}_g)
\end{equation}
defined by ${\rm EXP}(\gamma,v)(t)=\exp_{\gamma(t)}v(t)\;\forall
t\in\R$ via the exponential map $\exp$ of $g$,  restricts to a $C^{k-3}$
diffeomorphism
$\digamma: N{\cal O}(\varepsilon)\to {\cal
N}({\cal O},\varepsilon)$,
where ${\cal N}({\cal O},\varepsilon)$ is  an open
neighborhood of ${\cal O}$ in $\Lambda(M, \mathbb{I}_g)$.
Similarly, $\digamma$ restricts to a $C^{k-3}$
diffeomorphism from $XN{\cal O}(\varepsilon)$ to an open neighborhood
${\cal X}({\cal O},\varepsilon)$ of ${\cal O}$ in ${\cal X}(M,\mathbb{I}_g)$
(shrinking $\varepsilon>0$ if necessary), which implies that
${\cal L}\circ \digamma|_{XN{\cal O}(\varepsilon)}$ is $C^{k-3}$
as in the proof of Lemma~\ref{lem:2.3}.
For $\gamma\in{\cal O}$ let $N{\cal
O}(\varepsilon)_\gamma$, $N{\cal O}(\varepsilon)^X_\gamma$ and
$XN{\cal O}(\varepsilon)_\gamma$ be the fibers of fiber bundles
$N{\cal O}(\varepsilon)$, $N{\cal O}(\varepsilon)^X$ and $XN{\cal
O}(\varepsilon)$ at $\gamma\in{\cal O}$, respectively. Define
$$
{\cal
F}={\cal L}\circ \digamma,\qquad
{\cal F}_\gamma={\cal F}|_{N{\cal
O}(\varepsilon)_\gamma},\qquad {\cal F}^X={\cal F}|_{N{\cal O}(\varepsilon)^X}.
$$
 Denote by $A_\gamma$ the
restriction of the gradient $\nabla{\cal F}_\gamma$ to $N{\cal
O}(\varepsilon)^X_\gamma$. It takes values in $XN{\cal O}_\gamma$,
and for $\varepsilon>0$ small enough $A_\gamma$ is a $C^{k-3}$ map from
$XN{\cal O}(\varepsilon)_\gamma$ to $XN{\cal O}_\gamma$
by Lemma~\ref{lem:3.1}, (\ref{e:3.25}) and (\ref{e:3.35}). Moreover
the bundle map
\begin{equation}\label{e:1.18}
XN{\cal O}(\varepsilon)\ni(\gamma, v)\mapsto (\gamma,
A_{\gamma}(v))\in XN{\cal O}\quad\hbox{is}\quad C^{k-4}.
\end{equation}
(Indeed, since ${\cal F}|_{XN{\cal O}(\varepsilon)}$ is $C^{k-3}$, its fiberwise differential
$d_F({\cal F}|_{XN{\cal O}(\varepsilon)})$ is a $C^{k-4}$ bundle map from
$XN{\cal O}(\varepsilon)$ to $XN{\cal O}^\ast$. So any $C^{k-4}$ sections $\mathfrak{V}:{\cal O}\to
XN{\cal O}$ and $\mathfrak{U}:{\cal O}\to XN{\cal O}(\varepsilon)$ give rise to a $C^{k-4}$
 function
\begin{eqnarray*}
{\cal O}\ni \gamma\mapsto d_F({\cal F}|_{XN{\cal O}(\varepsilon)})(\mathfrak{V}(\gamma))[\mathfrak{U}(\gamma)]&=&
d({\cal F}|_{XN{\cal O}(\varepsilon)_\gamma})(\mathfrak{V}(x))[\mathfrak{U}(\gamma)]\\
&=&\langle A_\gamma(\mathfrak{V}(\gamma)), \mathfrak{U}(\gamma)\rangle_{1,\gamma}.
\end{eqnarray*}
By using a local trivialization argument, (\ref{e:1.18}) follows from this
because the Riemannian metric by given (\ref{e:1.1}) on $XN{\cal O}$ is of class $C^{k-2}$.)
It is clear that
\begin{equation}\label{e:1.19}
A_{s\cdot \gamma}(s\cdot v)=s\cdot A_\gamma(v)\quad\forall s\in
\R,\; v\in N{\cal O}(\varepsilon)^X_\gamma
\end{equation}
since ${\cal O}$ and ${\cal F}$ are $\R$-invariant.
  Denote by $B_\gamma$ the natural
extension of the symmetric bilinear form $d^2({\cal F}|_{XN{\cal
O}(\varepsilon)_\gamma})(0)$ on $N{\cal O}_\gamma$ and  the
associated self-adjoint operator with it. The latter is Fredholm, and
has finite dimensional negative definite and null spaces ${\bf
H}^-(B_\gamma)$ and ${\bf H}^0(B_\gamma)$. Moreover, ${\bf
H}^-(B_\gamma)+ {\bf H}^0(B_\gamma)$ is contained in $XN{\cal
O}_\gamma$, and $B_\gamma$ (as an element of $\mathscr{L}(N{\cal O}_\gamma)$) also restricts to a bounded linear
operator $dA_\gamma(0)\in\mathscr{L}(XN{\cal O}_\gamma)$. Because
$B_{s\cdot \gamma}(s\cdot u)=s\cdot B_\gamma(u)$ for any $s\in \R$
and $(\gamma,u)\in N{\cal O}$, we deduce that $\dim{\bf
H}^0(B_\gamma)$ and  $\dim{\bf H}^-(B_\gamma)$ are independent of
$\gamma\in{\cal O}$, and call
\begin{equation}\label{e:1.20}
m^0({\cal L}, {\cal O}):=\dim{\bf
H}^0(B_\gamma)\quad\hbox{and}\quad m^-({\cal L}, {\cal O}):=\dim{\bf H}^-(B_\gamma)
\end{equation}
the {\bf nullity} and {\bf Morse index} of ${\cal O}$, respectively.
In the case $m^0({\cal L}, {\cal O})=0$ the critical orbit
${\cal O}$ is said to be {\bf nondegenerate}. Note that we have
always $0\le m^0({\cal L}, {\cal O})\le 2n-1$. Observe that
there exists $\mathfrak{a}>0$ such that
 \begin{equation}\label{e:1.21}
\sigma(B_\gamma)\cap([-2\mathfrak{a},
2\mathfrak{a}]\setminus\{0\})=\emptyset,\quad\forall\gamma\in{\cal
O}
\end{equation}
 since $B_{s\cdot\gamma}(s\cdot u)=s\cdot B_\gamma(u)$ for any $s\in \R$ and
$(\gamma,u)\in N{\cal O}$.
  By this and Lemma~7.3 on the page 70 of \cite{Ch93}, for the critical manifold ${\cal O}$
   we have a natural $C^{k-2}$ Hilbert vector bundle orthogonal decomposition
\begin{equation}\label{e:1.22}
N{\cal O}={\bf H}^-(B)\oplus {\bf
H}^0(B)\oplus{\bf H}^+(B)
\end{equation}
with ${\bf H}^\star(B)_\gamma={\bf H}^\star(B_\gamma)$ for
$\gamma\in{\cal O}$ and $\star=+,0,-$, and hence (by
Proposition~\ref{prop:5.1}) a $C^{k-2}$ Banach vector bundle direct
sum decomposition
\begin{equation}\label{e:1.23}
 XN{\cal O}={\bf H}^-(B)\dot{+} {\bf
H}^0(B)\dot{+}({\bf H}^+(B)\cap XN{\cal O}).
\end{equation}
Let ${\bf P}^\star$ be the orthogonal bundle projections from
$N{\cal O}$ onto ${\bf H}^\star(B)$, $\star=+,0,-$, and let ${\bf
H}^0(B)(\delta)={\bf H}^0(B)\cap N{\cal O}(\delta)$ for $\delta>0$.
Note that ${\bf H}^0(B)(\delta)\subset XN{\cal O}$ and that
$\delta>0$ can be chosen so small that ${\bf H}^0(B)(\delta)\subset XN{\cal
O}(\varepsilon)$ since ${\bf H}^0(B)$ has finite rank and $B_{s\cdot
\gamma}(s\cdot u)=s\cdot B_\gamma(u),\;\forall (s,u,\gamma)$.  By
a local trivialization argument  and the compactness of ${\cal O}$ we
may use the implicit function theorem and (\ref{e:1.18}) to get a sufficiently small
$\delta>0$ and  a unique $\R$-equivariant $C^{k-4}$ bundle map
\begin{equation}\label{e:1.24}
\mathfrak{h}:{\bf H}^0(B)(\delta)\to {\bf H}^-(B)\dot{+}({\bf
H}^+(B)\cap XN{\cal O})
\end{equation}
whose restriction $\mathfrak{h}_\gamma$ on each fibre ${\bf
H}^0(B)(\delta)_\gamma$ is of class $C^{k-3}$, $\R_\gamma$-equivariant and satisfies
\begin{eqnarray}\label{e:1.25}
\mathfrak{h}_\gamma(0)=0,\quad ({\bf P}^+_\gamma +{\bf
P}^-_\gamma)\circ A_\gamma\bigl(v+ \mathfrak{h}_\gamma(v)\bigr)=0,
\quad\forall v\in {\bf H}^0(B)(\delta)_\gamma.
\end{eqnarray}
(Here the $C^{k-4}$ (resp. $C^{k-3}$) smoothness of $\mathfrak{h}$ (resp. $\mathfrak{h}_\gamma$)
 comes from the fact that the bundle map in (\ref{e:1.18}) (resp.  $A_\gamma$)
  is of class $C^{k-4}$ (resp. $C^{k-3}$)).
Moreover,  the $\R$-invariant functional
\begin{eqnarray}\label{e:1.26}
{\cal L}^\circ_\triangle:{\bf H}^0(B)(\delta)\ni (\gamma,v)\to{\cal
L}\circ{\rm EXP}_\gamma\bigl(v+
 \mathfrak{h}_\gamma(v)\bigr)\in\R
\end{eqnarray}
is $C^{k-4}$, and restricts to an $\R_\gamma$-invariant and $C^{k-3}$ functional in each fiber
${\bf H}^0(B)(\delta)_\gamma$, denoted by ${\cal L}^\circ_{\triangle
\gamma}$.

\begin{theorem}\label{th:1.7}
Under the  above notation there exist $\varrho\in (0, \delta)$, an
$\R$-invariant open neighborhood $U$ of the zero section of $N{\cal
O}$, an $\R$-equivariant $C^{k-4}$ fiber map $\mathfrak{h}$ given by
(\ref{e:1.24}), and an $\R$-equivariant fiber-preserving
homeomorphism $\Upsilon:N{\cal O}(\varrho)\to U$ such that for all
$(\gamma, u)\in N{\cal O}(\varrho)$,
$$
{\cal L}\circ{\rm EXP}\circ\Upsilon(\gamma, u)=\|{\bf P}_\gamma^+u\|^2_1-\|{\bf
P}_\gamma^-u\|^2_1+ {\cal L}^\circ_\triangle(\gamma, {\bf P}_\gamma^0u).
$$
 Moreover, $\Upsilon$ also satisfies the following properties:
\begin{description}
\item[(i)] for any $(\gamma,u^0\oplus u^+\oplus u^-)\in N{\cal
 O}(\varrho)\cap\bigl({\bf H}^0(B)\oplus {\bf
H}^-(B)\oplus {\bf H}^+(B)\bigr)$,
$\Upsilon_\gamma(0)=0$, $\Upsilon_\gamma(u^0\oplus
u^-\oplus u^+)-\mathfrak{h}_{\gamma}(u^0)-u^0\in {\bf H}^-(B)_\gamma\oplus {\bf
H}^+(B)_\gamma$,  $\Upsilon_\gamma(u^0\oplus u^-\oplus u^+)\in {\bf
H}^-(B)_\gamma$ if and only if $u^0=0$ and $u^+=0$, and
$\Upsilon_\gamma(u^0\oplus u^-\oplus u^+)\in {\bf H}^-(B)_\gamma\oplus {\bf
H}^+(B)_\gamma$ if and only if $u^0=0$;
\item[(ii)] $\Upsilon$ restricts to a $C^{k-4}$ bundle embedding from $N{\cal
 O}(\varrho)\cap{\bf H}^0(B)$ to $XN{\cal O}$ (and so to $N{\cal O}$) because $\Upsilon(\gamma, v)=
 (\gamma, v+\mathfrak{h}_\gamma(v))$ for $(\gamma, v)\in N{\cal O}(\varrho)\cap{\bf H}^0(B)$;

\item[(iii)]  $\Upsilon\left(({\bf P}^-+{\bf P}^0)N{\cal
O}(\varrho)\right)$ is contained in $XN{\cal O}$, and $\Upsilon$ is
also a homeomorphism from $({\bf P}^-+{\bf P}^0)N{\cal O}(\varrho)$
onto $\Upsilon\left(({\bf P}^-+{\bf P}^0)N{\cal O}(\varrho)\right)$
even if the topology on the latter is taken as the induced one by
$XN{\cal O}$. {\rm (}This implies that $N{\cal O}$ and $XN{\cal O}$ induce
the same topology in $\Upsilon\left(({\bf P}^-+{\bf P}^0)N{\cal
O}(\varrho)\right)$.{\rm )}
\end{description}
\end{theorem}

Since every curve in a small neighborhood
of ${\cal O}$ in $\Lambda(M,\mathbb{I}_g)$ has image near the compact subset
 $\gamma_0(\R)\subset M$, as showed in Remark~\ref{rm:1.3}, Theorem~\ref{th:1.7}
 only depends on the values of the metric $g$ near $\gamma_0(\R)$.

Under the assumptions of Theorem~\ref{th:1.7} let ${\bf
H}^{0-}(B):={\bf H}^0(B)\oplus{\bf H}^-(B)$ and ${\bf
H}^{0-}(B)(\epsilon):=({\bf H}^0(B)\oplus{\bf H}^-(B))\cap N{\cal
O}(\epsilon)$. Then ${\bf H}^{0-}(B)\subset XN{\cal O}$. Define
$\mathfrak{L}:{\bf H}^{0-}(B)(\epsilon)\to\R$ by
\begin{eqnarray}\label{e:1.27}
 \mathfrak{L}(\gamma,v)=-\|{\bf P}^-_\gamma v\|^2_1+ {\cal
L}^\circ_{\triangle \gamma}({\bf P}^0_\gamma v).
\end{eqnarray}
Let $C_\ast({\cal L},{\cal O};\K)$ and $C_\ast(\mathfrak{L},{\cal O};\K)$
denote the critical groups of ${\cal L}$ and $\mathfrak{L}$ at isolated critical orbits
${\cal O}$  with the coefficient group $\K$, respectively.
 By some  standard deformation arguments we may derive from
Theorem~\ref{th:1.7}
\begin{eqnarray}\label{e:1.28}
C_\ast({\cal L},{\cal O};\K)\cong C_\ast(\mathfrak{L},{\cal O};\K).
\end{eqnarray}
Define $\mathscr{N}({\cal O}):={\rm EXP}\circ\Upsilon\bigl(N{\cal
 O}(\varrho)\cap{\bf H}^0(B)\bigr)$. This is an $\R$-invariant $C^{k-4}$ (embedded) submanifold of ${\cal X}(M,\mathbb{I}_g)$ (and so of $\Lambda(M,\mathbb{I}_g)$)
which has dimension $1+ m^0({\cal L}, {\cal O})$ and contains ${\cal O}$ in its interior.
 But $\mathscr{N}({\cal O})_\gamma:={\rm EXP}_\gamma\circ\Upsilon_\gamma\bigl(N{\cal
 O}(\varrho)_\gamma\cap{\bf H}^0(B)_\gamma\bigr)$, for each $\gamma\in{\cal O}$,
  is a $C^{k-3}$ (embedded) submanifold of ${\cal X}(M,\mathbb{I}_g)$ (and so of $\Lambda(M,\mathbb{I}_g)$)
which has dimension $m^0({\cal L}, {\cal O})$ and contains $\gamma$ as an interior point.
We call $\mathscr{N}({\cal O})$ (resp. $\mathscr{N}({\cal O})_\gamma$)
 a \textsf{characteristic submanifold} of $(\Lambda(M,\mathbb{I}_g), \langle\cdot,\cdot\rangle_1)$
for ${\cal L}$ at ${\cal O}$ (resp. $\gamma\in{\cal O}$), and call
$$
\mathscr{H}_\ast^0({\cal L},\gamma;\K):=C_\ast({\cal L}|_{\mathscr{N}({\cal O})_\gamma},\gamma;\K)
$$
the \textsf{characteristic invariant} of ${\cal L}$ at $\gamma\in{\cal O}$.
Clearly, $\mathscr{H}_\ast^0({\cal L},\gamma;\K)\cong C_\ast({\cal L}^\circ_{\triangle
\gamma}, 0;\K)$ and\linebreak $\mathscr{H}_q^0({\cal L},\gamma;\K)=0$
 for all $q\ge 2\dim M-1$.
Theorem~\ref{th:1.7} also implies the following \textsf{shifting theorem},
\begin{eqnarray}\label{e:1.29}
C_{q+ m^-({\cal L},{\cal O})}({\cal F}_\gamma, 0;\K)\cong C_q({\cal L}^\circ_{\triangle
\gamma}, 0;\K)\cong \mathscr{H}_\ast^0({\cal L},\gamma;\K)\quad\forall q=0,1,\cdots.
\end{eqnarray}
By Meyer-Vietoris theorem (cf. \cite{Tan82}),  if the Abelian group $\K$ is a field we get the following
important inequalities
\begin{eqnarray}\label{e:1.30}
\dim C_q({\cal L},{\cal O};\K)\le 2\sum^1_{i=0}
\dim\mathscr{H}_{q- m^-({\cal L}, {\cal O})-i}^0({\cal L},\gamma;\K)\quad\forall q,
\end{eqnarray}
which are sufficient for most of applications.

For a curve $c:\R\to M$ and $\tau\in\R$ we define a curve $c^\tau:\R\to M$ by
\begin{eqnarray}\label{e:1.31}
c^\tau(t)=c(\tau\cdot t)\quad\forall t\in\R.
\end{eqnarray}
In particular, for each $m\in\N$ we have an iteration map
\begin{eqnarray}\label{e:1.32}
\varphi_m: \Lambda(M,\mathbb{I}_g)\to \Lambda(M,\mathbb{I}^m_g)
\end{eqnarray}
given by $\varphi_m(c)=c^m$.
For any closed $\mathbb{I}_g$-invariant $F$-geodesic $\gamma_0$ and every $m\in\N$,
$\gamma_0^m$ is a closed $\mathbb{I}^m_g$-invariant $F$-geodesic, and so a critical point
of the functional ${\cal L}$ defined by (\ref{e:1.2}) on $\Lambda(M,\mathbb{I}^m_g)$.
Let $T\Lambda(M,\mathbb{I}^m_g)$ be equipped with the equivalent Riemannian-Hilbert structure
\begin{equation}\label{e:1.33}
(\xi,\eta)_m=m^2\int^1_0\langle\xi(t),\eta(t)\rangle dt+
\int^1_0\langle\nabla^g_\gamma\xi(t),\nabla^g_\gamma\xi(t)\rangle dt.
\end{equation}
 Denote by $\|\cdot\|_m$ the norm of it. Let
 \begin{equation}\label{e:1.34}
 {\rm EXP}^m: T\Lambda(M,\mathbb{I}^m_g)(\varepsilon)=\{(\gamma,v)\in
T\Lambda(M,\mathbb{I}^m_g)\,|\,\|v\|_m\}<\varepsilon\}\to
\Lambda(M,\mathbb{I}^m_g)
\end{equation}
be defined as in (\ref{e:1.17}).
Suppose that  ${\cal
O}^m:=\R\cdot\gamma_0^m$ is also an isolated critical orbit of ${\cal L}$ in $\Lambda(M,\mathbb{I}^m_g)$.
Denote by $\hat{N}{\cal O}^m$ the normal bundle of ${\cal
O}^m$ in $\Lambda(M,\mathbb{I}^m_g)$ with respect to the
Riemannian-Hilbert structure
(\ref{e:1.33}), and by $X\hat{N}{\cal O}^m:=T_{{\cal O}^m}{\cal X}(M, \mathbb{I}^m_g)\cap
\hat{N}{\cal O}^m$ the Banach vector subbundle of
$T_{{\cal O}^m}{\cal X}(M, \mathbb{I}^m_g)$. Correspondingly, we have also
 $X\hat{N}{\cal O}^m(\varepsilon)$, $\hat{N}{\cal O}^m(\varepsilon)^X$ and $\hat{N}{\cal
O}^m(\varepsilon)$. Then
\begin{equation}\label{e:1.35}
\hat{\cal F}: \hat N{\cal O}^m(\varepsilon)\to
\R,\;(y,v)\mapsto {\cal L}\circ{\rm EXP}^m(y,v)
\end{equation}
is an  $\R$-invariant functional of class $C^{2-0}$, and also restricts to
a $C^{k-3}$ map on $X\hat N{\cal O}^m(\varepsilon)$ for $\varepsilon>0$ small enough.
For each $y\in{\cal O}^m$ let $\hat{\cal F}_y$
denote the restriction of  $\hat{\cal F}$ to $\hat N{\cal
O}^m(\varepsilon)_y$, and let
 $\hat{\cal F}^X_y$ be the restrictions of  $\hat{\cal F}_y$ to $\hat N{\cal
O}^m(\varepsilon)_y\cap X\hat{N}{\cal O}^m_y$.
 Denote by $\hat{\nabla}\hat{\cal F}_y$  the gradient of $\hat{\cal F}_y$ with
respect to the inner product in (\ref{e:1.33}) on $\hat
N{\cal O}^m_y$, and by  $\hat A_y$ be the restriction of
$\hat\nabla\hat{\cal F}_y$ to $\hat{N}{\cal
O}^m(\varepsilon)_y\cap X\hat{N}{\cal O}^m_y$. Clearly,
\begin{equation}\label{e:1.36}
\hat A_{s\cdot y}(s\cdot v)=s\cdot \hat A_y(v)\quad\forall s\in
\R,\; (y,v)\in \hat{N}{\cal O}^m(\varepsilon)\cap
X\hat{N}{\cal O}^m.
\end{equation}
Moreover, if $\delta\in (0, \varepsilon)$ is small enough we have:\\
(i) $X\hat{N}{\cal O}^m(\delta)\ni (y,v)\mapsto (y, \hat A_y(v))\in X\hat{N}{\cal O}^m$
is a $C^{k-4}$ bundle map; \\
 (ii) the map $\hat A_y$  is $C^{k-3}$-smooth from
$X\hat{N}{\cal O}^m(\delta)_y$ to $X\hat{N}{\cal
O}^m_y$ (and so $\hat{\cal F}^X_y$ is $C^{k-2}$ in
$X\hat{N}{\cal O}^m(\delta)_y$);\\
 (iii) the symmetric
bilinear form $d^2\hat{\cal F}^X_y(0)$ has a continuous extension
$\hat{B}_y$ on $\hat{N}{\cal O}^m_y$;\\
 (iv) ${\bf H}^-(\hat{B}_y)+ {\bf H}^0(\hat{B}_y)\subset X\hat{N}{\cal
O}^m_y$ is of finite dimension and there exist a  $C^{k-2}$ Hilbert vector bundle orthogonal
decomposition
\begin{eqnarray}\label{e:1.37}
\hat{N}{\cal O}^m={\bf H}^-(\hat{B})\hat\oplus {\bf
H}^0(\hat{B})\hat\oplus{\bf H}^+(\hat{B})
 \end{eqnarray}
with respect to the inner product in (\ref{e:1.33}) and an induced $C^{k-2}$ Banach
space bundle direct sum decomposition $X\hat{N}{\cal O}^m={\bf
H}^-(\hat{B}){\dot{+} } {\bf H}^0(\hat{B}){\dot{+} }({\bf
H}^+(\hat{B})\cap X\hat{N}{\cal O}^m)$. Let $m^-({\cal L},{\cal O}^m)$
 (resp. $m^0({\cal L}, {\cal O}^m)$) be the Morse index (resp. the nullity)
 defined by (\ref{e:1.20}) when ${\cal O}$ is replaced by ${\cal O}^m$.
(Precisely, $m^0({\cal L}, {\cal O}^m):=\dim{\bf H}^0(B_y)$ and $m^-({\cal L}, {\cal O}^m):=\dim{\bf H}^-(B_y)$ for $y\in {\cal O}^m$).  Then
\begin{eqnarray}\label{e:1.38}
m^-({\cal L}, {\cal O}^m)=\dim{\bf H}^-(\hat{B}_y)\quad\hbox{and}\quad
m^0({\cal L}, {\cal O}^m)=\dim{\bf H}^0(\hat{B}_y)
 \end{eqnarray}
 because two sides of the first (resp. second) equality are equal to the Morse index
 (resp. nullity minus one) of the symmetric
bilinear form $d^2({\cal L}|_{{\cal X}(M, \mathbb{I}^m_g)})(\gamma_0^m)$.
 Let $\hat{\bf P}^\star_y$ be the orthogonal projections from
$\hat{N}{\cal O}^m_y$ onto ${\bf H}^\star(\hat{B}_y)$ in
(\ref{e:1.37}), $\star=+,0,-$. Since $\dim{\bf H}^0(\hat{B}_y)$ is
finite we may take $\epsilon\in (0,\delta)$
so small that
$$
{\bf H}^0(\hat{B})(\epsilon)_y:={\bf
H}^0(\hat{B}_y)\cap \hat{N}{\cal
O}^m(\epsilon)_y\subset X\hat{N}{\cal
O}^m(\delta)_y
$$
 and use  the implicit function theorem to get a
unique  $C^{k-3}$ map
\begin{equation}\label{e:1.39}
 \hat{\mathfrak{h}}_y:{\bf H}^0(\hat{B})(\epsilon)_y\to
 {\bf H}^-(\hat{B}_y){\dot{+}}({\bf H}^+(\hat{B}_y)\cap
X\hat{N}{\cal O}^m_y)
\end{equation}
satisfying $\hat{\mathfrak{h}}_y(0)=0$, $d\hat{\mathfrak{h}}_y(0)=0$
and
\begin{eqnarray}\label{e:1.40}
 (\hat{\bf P}^+_y +\hat{\bf P}^-_y)\circ \hat{A}_y\bigl(v+ \hat{\mathfrak{h}}_y(v)\bigr)=0
\quad\forall v\in {\bf H}^0(\hat{B})(\epsilon)_y.
\end{eqnarray}
By (\ref{e:1.36}) the map $\hat{\mathfrak{h}}_y$ is also
$\R_y$-equivariant; moreover the bundle map
$$
\hat{\mathfrak{h}}:{\bf H}^0(\hat{B})(\epsilon)\to
 {\bf H}^-(\hat{B}){\dot{+}}({\bf H}^+(\hat{B})\cap
X\hat{N}{\cal O}^m), (y,v)\mapsto (y,\hat{\mathfrak{h}}_y)
$$
is  $C^{k-4}$.  Define the functional
\begin{eqnarray}\label{e:1.41}
\hat{\cal L}^\circ_{\triangle y}:{\bf
H}^0(\hat{B})(\epsilon)_y\ni v\to \hat{\cal F}_y\bigl(v+
 \hat{\mathfrak{h}}_y(v)\bigr)\in\R.
\end{eqnarray}
It is  $C^{k-3}$ and has the isolated critical point $0$.
On the other hand the functional 
$$
\hat{\cal L}^\circ_{\triangle}:{\bf H}^0(\hat{B})(\epsilon)\to\R
$$
defined by $\hat{\cal L}^\circ_{\triangle}(y, v)=\hat{\cal L}^\circ_{\triangle y}(v)$
is only  $C^{k-4}$. Define the map
$\hat{\mathscr{J}}_{y}:\overline{{\bf
 H}^0(\hat{B})(\epsilon)_{y}}\times \bigr({\bf
 H}^-(\hat{B})(\epsilon)_{y}\oplus
 {\bf  H}^+(\hat{B})(\epsilon)_{y}\bigl)\to\R$  by
\begin{equation}\label{e:1.42}
\hat{\mathscr{J}}_{y}(u, v)=\hat{\cal F}_{y}(u+
\hat{\mathfrak{h}}_{y}(u)+ v)-\hat{\cal F}_{y}(u+
\hat{\mathfrak{h}}_{y}(u)).
\end{equation}
Combining  the proof method of \cite[Claim~6.3]{Lu3} with that of Theorem~\ref{th:1.7}
we can show that $\hat{\mathscr{J}}_{y}$ satisfies an analogous result
to Proposition~\ref{prop:3.8}. So shrinking  $\epsilon>0$ (if necessary)
 as in the proof of \cite[Lemma 3.6]{Lu2}, we obtain
a $\R_{y}$-invariant origin-preserving homeomorphism
$\widehat\Upsilon_{y}$ from $\hat{N}{\cal O}^m(\epsilon)_{y}$ onto a
neighborhood of $0_{y}\in\hat{N}{\cal O}^m_{y}$ and satisfying
$$
\hat{\cal F}_{y}\circ\widehat\Upsilon_{y}(u)=\|\hat{\bf
P}^+_{y_0}u\|^2_m-\|\hat{\bf P}^-_{y}u\|^2_1+ \hat{\cal
L}^\circ_{\triangle y}(\hat{\bf P}^0_{y}u)\quad\forall
u\in \hat{N}{\cal O}^m(\epsilon)_{y}.
$$
Summing up, we have:

\begin{theorem}\label{th:1.8}
Under the  above notation,
by shrinking the above $\epsilon>0$ (if necessary)
 there exist  an
$\R$-invariant open neighborhood $U$ of the zero section of $\hat{N}{\cal
O}^m$, an $\R$-equivariant $C^{k-4}$ fiber map $\hat{\mathfrak{h}}$ given by
(\ref{e:1.39}), and an $\R$-equivariant fiber-preserving
homeomorphism $\widehat{\Upsilon}:\hat{N}{\cal O}^m(\epsilon)\to \widehat{U}$ such that for all
$(\gamma, u)\in \hat{N}{\cal O}^m(\epsilon)$,
$$
{\cal L}\circ{\rm EXP}^m\circ\widehat{\Upsilon}(y, u)=\|\hat{\bf P}_y^+u\|^2_m-\|\hat{\bf
P}_y^-u\|^2_m+ \hat{\cal L}^\circ_\triangle(y, \hat{\bf P}_y^0u).
$$
Moreover, $\widehat{\Upsilon}$ also satisfies similar properties to ``moreover" part of Theorem~\ref{th:1.7};
in particular $\widehat{\Upsilon}$ restricts to a $C^{k-4}$ bundle embedding from $\hat{N}{\cal O}^m(\epsilon)\cap{\bf H}^0(\hat{B})$ to $X\hat{N}{\cal O}^m$ (and so to $\hat{N}{\cal O}^m$) because $\widehat{\Upsilon}(y, u)=
 (y, u+ \hat{\mathfrak{h}}_y(u))$ for $(y, u)\in \hat{N}{\cal O}^m(\epsilon)\cap{\bf H}^0(\hat{B})$;
\end{theorem}

 $\hat{\mathscr{N}}({\cal O}^m):={\rm EXP}^m\circ\widehat{\Upsilon}\bigl(
\hat{N}{\cal O}^m(\epsilon)\cap{\bf H}^0(\hat{B})\bigr)$ is a $C^{k-4}$ (embedded) submanifold of ${\cal X}(M,\mathbb{I}^m_g)$ (and so of $\Lambda(M,\mathbb{I}^m_g)$)
which contains ${\cal O}^m$ in its interior. But for each $y\in{\cal O}^m$, $\hat{\mathscr{N}}({\cal O}^m)_y:={\rm EXP}^m_y\circ\widehat{\Upsilon}_y\bigl(
\hat{N}{\cal O}^m(\epsilon)_y\cap{\bf H}^0(\hat{B})_y\bigr)$ is a $C^{k-3}$ (embedded) submanifold of ${\cal X}(M,\mathbb{I}^m_g)$ (and so of $\Lambda(M,\mathbb{I}^m_g)$)
which has dimension $m^0({\cal L}, {\cal O}^m)$ and which contains $y$ as an interior point, and is called
 a \textsf{characteristic submanifold} of the Hilbert-Riemannian manifold $(\Lambda(M,\mathbb{I}^m_g), (\cdot,\cdot)_m)$ for ${\cal L}$ at $y\in{\cal O}^m$. We also call
$$
\hat{\mathscr{H}}_\ast^0({\cal L}, y;\K):=C_\ast({\cal L}|_{\hat{\mathscr{N}}({\cal O}^m)_y},y;\K)
$$
the \textsf{characteristic invariant} of ${\cal L}$ at $y\in{\cal O}^m$ because as in the proof of \cite[(6.31)]{Lu3}, we can show that
$$
C_\ast({\cal L}|_{\hat{\mathscr{N}}({\cal O}^m)_y},y;\K)=C_\ast({\cal L}|_{{\mathscr{N}}({\cal O}^m)_y},y;\K)
\quad\forall y\in{\cal O}^m,
$$
 where ${\mathscr{N}}({\cal O}^m)_y$  is the characteristic submanifold of the Hilbert-Riemannian manifold $(\Lambda(M,\mathbb{I}^m_g),\\ \langle\cdot,\cdot\rangle_1)$ for ${\cal L}$ at $y\in{\cal O}^m$.

As for Theorem~\ref{th:1.4} we have

\begin{theorem}\label{th:1.9}
Under the assumptions of Theorem~\ref{th:1.8}, let
$\tilde\Lambda$ be a Riemannian-Hilbert submanifold of $\bigl(\Lambda(M, \mathbb{I}^m_g),
(\cdot,\cdot)_m\bigr)$ containing ${\cal O}^m$
and satisfying the following conditions:
\begin{description}
\item[(i)]  $\tilde{\cal X}:={\cal X}(M, \mathbb{I}^m_g)\cap\tilde\Lambda$ is a Banach submanifold of
${\cal X}(M, \mathbb{I}^m_g)$; 

\item[(ii)]   ${\rm EXP}^m$ in (\ref{e:1.34}) restricts to a diffeomorphism
 $\widetilde{\rm EXP}^m$ (resp.  $\widetilde{\rm EXP}^{mX}$ )  from
$$
\tilde{N}{\cal O}^m(\varepsilon):=\hat{N}{\cal O}^m(\varepsilon)\cap T_{{\cal O}^m}\tilde\Lambda
\quad{\rm (resp.}\,
X\tilde{N}{\cal O}^m(\varepsilon):=X\hat{N}{\cal O}^m(\varepsilon)\cap T_{{\cal O}^m}\tilde{\cal L}\;
{\rm )}
$$
to an open neighborhood of ${\cal O}^m$ in $\tilde\Lambda$ (resp.  $\tilde{\cal X}$)
(by shrinking $\varepsilon>0$ if necessary);

\item[(iii)] $\hat\nabla({\cal L}\circ\widetilde{\rm EXP}_{y})(u)=\hat\nabla\hat{\cal F}_y(u)$ for any $(y, u)\in \tilde{N}{\cal O}^m(\varepsilon)$;

\item[(iv)] ${\bf H}^0\bigl(d^2({\cal L}|_{{\cal X}(M, \mathbb{I}^m_g)})(\gamma_0^m)\bigr)\subset
T_{\gamma_0^m}\tilde\Lambda$ (which is equivalent to the fact ${\bf H}^0\bigl(d^2\hat{\cal F}^X_{\gamma_0^m}(0)\bigr)\subset X\tilde{N}{\cal O}^m_{\gamma_0^m}$).
\end{description}
Then the splitting lemma in Theorem~\ref{th:1.8}
restricts to that of ${\cal L}\circ\widetilde{\rm EXP}$ around ${\cal O}^m$ in $\tilde\Lambda$,
precisely  $\widehat{\Upsilon}$ in Theorem~\ref{th:1.8} restricts to
the $\R$-equivariant fiber-preserving homeomorphism $\widetilde{\Upsilon}$ from
$\tilde{N}{\cal O}^m(\epsilon):=\hat{N}{\cal O}^m(\epsilon)\cap T_{{\cal O}^m}\tilde\Lambda$ to $\widehat{U}\cap\tilde\Lambda$ such that
$$
{\cal L}\circ\widetilde{\rm EXP}^m\circ\widetilde{\Upsilon}(y, u)=\|\hat{\bf P}_y^+u\|^2_m-\|\hat{\bf
P}_y^-u\|^2_m+ \hat{\cal L}^\circ_\triangle(y, \hat{\bf P}_y^0u)
$$
for all $(y, u)\in \tilde{N}{\cal O}^m(\epsilon)$ (by shrinking $\epsilon>0$ if necessary).
 Moreover the inclusion $\tilde\Lambda\hookrightarrow \Lambda(M, \mathbb{I}^m_g)$ induces isomorphisms
\begin{description}
\item[(a)] $\hat{\mathscr{H}}_\ast^0({\cal L}|_{\tilde\Lambda}, y;\K)\cong \hat{\mathscr{H}}_\ast^0({\cal L}, y;\K)\;\forall y\in{\cal O}^m$ if $m^0({\cal L}, {\cal O}^m)=m^0({\cal L}|_{\tilde\Lambda}, {\cal O}^m)$;

\item[(b)] $C_\ast({\cal L}|_{\tilde\Lambda}, {\cal O}^m;\K)\cong C_\ast({\cal L}, {\cal O}^m;\K)$
if $m^\star({\cal L}, {\cal O}^m)=m^\star({\cal L}|_{\tilde\Lambda}, {\cal O}^m)$, $\star=-,0$.
\end{description}
\end{theorem}

\begin{remark}\label{rm:1.10}
{\rm (I) If $F^2(v)=g(v,v)$ for all $v\in TM$, this result follows from \cite[Lemma~7]{GM1} immediately.
(II) Assume that there exists a Hilbert space isomorphism on $\hat{N}{\cal O}^m_{\gamma_0^m}$ which preserves
$X\hat{N}{\cal O}^m_{\gamma_0^m}$ and $\hat{\cal F}_{\gamma_0^m}$, such that $\tilde{N}{\cal O}^m_{\gamma_0^m}$
is the fixed point set of it. Then (iii) of Theorem~\ref{th:1.9} is satisfied, and (iv)
of Theorem~\ref{th:1.9} can also hold if $m^0({\cal L}, {\cal O}^m)=m^0({\cal L}|_{\tilde\Lambda}, {\cal O}^m)$. The following corollary may be easily seen from the proof of Theorem~\ref{th:1.12} below.}
\end{remark}

\begin{corollary}\label{cor:1.11}
 Under the assumptions of Theorem~\ref{th:1.8} suppose for some $k\in\N$ that
 $\mathbb{I}_g^k(\gamma_0(t))=\gamma_0(t)\;\forall t\in\R$. Let $S$ be
 the connected component of ${\rm Fix}(\mathbb{I}^k_g)$ containing $\gamma_0(\R)$,
 which is a totally geodesic submanifold of $M$, and let $\Lambda(S, \mathbb{I}_g)$
 be equipped with the Riemannian-Hilbert structure   defined by  $g|_S$.  Then the inclusion
 $\Lambda(S, \mathbb{I}_g)\hookrightarrow \Lambda(M, \mathbb{I}_g)$ induces isomorphisms
 $$
\mathscr{H}^0_\ast({\cal L}|_{\Lambda(S, \mathbb{I}_g)}, \gamma;\K)\cong\mathscr{H}^0_\ast({\cal L}, \gamma;\K)\quad
 \forall\gamma\in{\cal O}
 $$
provided that  $\mathbb{I}_g\in I(M,F)$ and $m^0({\cal L}, {\cal O})=m^0({\cal L}|_{\Lambda(S, \mathbb{I}_g)}, {\cal O})$.
\end{corollary}

In the Riemannian case this result also holds if $S$ is a
totally geodesic submanifold of $(M,g)$ with $\gamma_0(\R)\subset S=\mathbb{I}_g(S)$
(Proposition~3.5 in \cite{GroTa78}), which is very important for studies of Riemannian isometry-invariant geodesics \cite{GroTa78, Tan82}. When $F^2(\cdot)\ne g(\cdot,\cdot)$ we are not able to prove that the gradient $\nabla{\cal L}$ is tangent to
the submanifold $\Lambda(S, \mathbb{I}_g)$ at any $\gamma\in\Lambda(S, \mathbb{I}_g)$.

 From Theorems~\ref{th:1.8},~\ref{th:1.9} we may obtain

 \begin{theorem}\label{th:1.12}
Under the assumptions of Theorem~\ref{th:1.8},
suppose that  ${\cal O}^m:=\R\cdot\gamma_0^m$ is also an isolated critical submanifold of ${\cal L}$ in $\Lambda(M,\mathbb{I}^m_g)$. Then the iteration map in (\ref{e:1.32}) induces isomorphisms
\begin{equation}\label{e:1.43}
(\varphi_m)_\ast: \mathscr{H}^0_\ast({\cal L}, \gamma;\K)\to
\hat{\mathscr{H}}^0_\ast({\cal L}, \gamma^m;\K)
\quad \forall\gamma\in{\cal O}
\end{equation}
if $m^0({\cal L}, {\cal O})=m^0({\cal L}, {\cal O}^m)$, and
\begin{equation}\label{e:1.44}
(\varphi_m)_\ast: C_\ast({\cal L}, {\cal O};\K)\to
C_\ast({\cal L}, {\cal O}^m;\K)
\end{equation}
if $m^\star({\cal L}, {\cal O})=m^\star({\cal L}, {\cal O}^m)$, $\star=-,0$.
\end{theorem}

When $\mathbb{I}_g=id_M$,
$\Lambda(M, \mathbb{I}_g)$ is equal to $\Lambda M=W^{1,2}(S^1,M)$, where
$S^1=\R/\Z=[0,1]/\{0,1\}$. If  $\gamma\in \Lambda M$ is a closed $F$-geodesic such that its
orbit ${\cal O}=S^1\cdot\gamma$ of the $S^1$-action defined by
(\ref{e:1.16}) with $p=1$ is an isolated critical one, all corresponding results above
hold true when the $\R$-actions are replaced by the $S^1$-actions.

In applications we also need  variants of some of the results above.
For $\tau>0$, define a map ${\bf T}_\tau:W^{1,2}_{\rm loc}(\R, M)\to W^{1,2}_{\rm loc}(\R, M)$ by
$({\bf T}_\tau\xi)(t)=\xi(t+\tau)\;\forall t\in\R$.
Suppose that $\gamma_0\in\Lambda(M, \mathbb{I}_g)$  has the least irrational period $\alpha>0$.
Following \cite[Section~3]{Tan82}, for each $m\in\N\cup\{0\}$  we
introduce  the Hilbert manifolds
$$
\Lambda^{m\alpha+1}(M, \mathbb{I}_g)=\bigl\{x\in W^{1,2}_{\rm loc}(\R, M)\,|\,
\mathbb{I}_gx={\bf T}_{m\alpha+1}x\;\bigr\},
$$
with the Riemannian metric $\langle\cdot,\cdot\rangle_{m\alpha+1}$  defined by
$$
\langle\!\langle\xi,\eta\rangle\!\rangle_{m\alpha+1}=\frac{1}{m\alpha+1}\int^{m\alpha+1}_0\langle\xi(t),\eta(t)\rangle dt+
(m\alpha+1)\int^{m\alpha+1}_0\langle\nabla^g_x\xi(t),\nabla^g_x\xi(t)\rangle dt
$$
for $\xi, \eta\in T_x \Lambda^{m\alpha+1}(M, \mathbb{I}_g)$, and
$$
\Lambda^{m\alpha}(M)=\bigl\{x\in W^{1,2}_{\rm loc}(\R, M)\,|\, {\bf T}_{m\alpha}x=x\bigr\},
$$
with the Riemannian metric $\langle\!\langle\cdot,\cdot\rangle\!\rangle_{m\alpha}$ as above.
Correspondingly, we have also the Banach manifolds
${\cal X}^{m\alpha+1}(M, \mathbb{I}_g)=\bigl\{x\in C^1(\R, M)\,|\,
\mathbb{I}_gx={\bf T}_{m\alpha+1}x\;\bigr\}$
and 
$$
{\cal X}^{m\alpha}(M)=\bigl\{x\in C^1(\R, M)\,|\, {\bf T}_{m\alpha}x=x\bigr\}.
$$ 
For each $l\in\N$,
${\bf T}_{l\alpha}$ gives not only an isometry from
$\Lambda^{m\alpha+1}(M, \mathbb{I}_g)$ to itself, but also one from $\Lambda^{l\alpha}(M)$
to itself. Moreover the fixed point sets of both the isometries are $\Lambda^{m\alpha+1}(M, \mathbb{I}_g)\cap
\Lambda^{l\alpha}(M)$. This set can be understood as totally geodesic submanifolds of
both $\Lambda^{m\alpha+1}(M, \mathbb{I}_g)$ and  $\Lambda^{l\alpha}(M)$, and
these two manifolds are diffeomorphic. Define the energy
functionals
\begin{eqnarray*}
&&^{m\alpha+1}\!{\cal L}(x):=\frac{1}{2}\int^{m\alpha+1}_0
F^2(x,\dot x)dt,\quad\forall x\in \Lambda^{m\alpha+1}(M, \mathbb{I}_g),\\
&& ^{l\alpha}\!{\cal L}(x):=\frac{1}{2}\int^{l\alpha}_0
F^2(x,\dot x)dt,\quad\forall x\in\Lambda^{l\alpha}(M),\\
&&^{m\alpha+1}\!{\cal L}^{l\alpha}:=\,^{m\alpha+1}\!{\cal L}|_{\Lambda^{m\alpha+1}(M, \mathbb{I}_g)\cap
\Lambda^{l\alpha}(M)},\\
&&^{l\alpha}\!{\cal L}^{m\alpha+1}:=\,^{l\alpha}\!{\cal L}|_{\Lambda^{m\alpha+1}(M, \mathbb{I}_g)\cap
\Lambda^{l\alpha}(M)}.
\end{eqnarray*}
There exists a Riemannian isometry
$$\psi_m: (\Lambda(M, \mathbb{I}_g), \langle\cdot,\cdot\rangle_1)\to
(\Lambda^{m\alpha+1}(M, \mathbb{I}_g), \langle\!\langle\cdot,\cdot\rangle\!\rangle_{m\alpha+1})
$$
defined by $\psi_m(x)(t)=x(t/(m\alpha+1))\;\forall t\in\R$, whose inverse is given by
 $\psi_m^{-1}(y)=y^{m\alpha+1}$ for $y\in\Lambda^{m\alpha+1}(M, \mathbb{I}_g)$. It is clear that
$$
(m\alpha+1){\cal L}=^{m\alpha+1}\!{\cal L}\circ\psi_m\quad\hbox{on}\;\Lambda(M, \mathbb{I}_g).
$$
Note that $\gamma_0^{m\alpha+1}$ also sits in $\Lambda(M, \mathbb{I}_g)$ and that
$\psi_m(\gamma_0^{m\alpha+1})=\gamma_0$. When $\R\cdot\gamma_0^{m\alpha+1}$
is an isolated critical submanifold in $\Lambda(M, \mathbb{I}_g)$, so is $\R\cdot\gamma_0$
in $\Lambda^{m\alpha+1}(M, \mathbb{I}_g)$; and
\begin{equation}\label{e:1.45}
m^\ast({\cal L}, \R\cdot\gamma_0^{m\alpha+1})=m^\ast(\,^{m\alpha+1}\!{\cal L},\R\cdot\gamma_0),\quad
\ast=0,-.
\end{equation}
It follows that $\psi_m$ induces isomorphisms
\begin{equation}\label{e:1.46}
\mathscr{H}^0_\ast({\cal L}, \gamma_0^{m\alpha+1};\K)\cong \mathscr{H}^0_\ast(\,^{m\alpha+1}\!{\cal L},\gamma_0;\K).
\end{equation}
By Remark~\ref{rm:1.10}(II) and Theorem~\ref{th:1.9}
we can arrive at the following result which is analogous to \cite[Lemma~3.2]{Tan82}.

\begin{theorem}\label{th:1.13}
Suppose that $\R\cdot\gamma_0^{m\alpha+1}$
is  an isolated critical submanifold for ${\cal L}$ in $\Lambda(M, \mathbb{I}_g)$
{\rm (}and hence $\R\cdot\gamma_0$ is not only one for
$^{m\alpha+1}\!{\cal L}$ in $\Lambda^{m\alpha+1}(M, \mathbb{I}_g)$
but also  for $^{m\alpha+1}\!{\cal L}^{l\alpha}$ in\linebreak $\Lambda^{m\alpha+1}(M, \mathbb{I}_g)\cap
\Lambda^{l\alpha}(M)${\rm )}. Then
$$
\mathscr{H}^0_\ast(\,^{m\alpha+1}\!{\cal L}, \gamma_0;\K)\cong \mathscr{H}^0_\ast(\,^{m\alpha+1}\!{\cal L}^{l\alpha}, \gamma_0;\K)
$$
 if $m^0\bigl({\cal L}, \R\cdot\gamma_0^{m\alpha+1}\bigr)
=m^0\bigl(\,^{m\alpha+1}\!{\cal L}^{l\alpha}, \R\cdot\gamma_0\bigr)$.
\end{theorem}

\subsection{The case $N=G(\mathbb{I}_F)$}\label{sec:1.3}

By an isometry on a Finsler manifold $(M,F)$ we mean a diffeomorphism
of $M$ onto itself whose differential  preserves the Finsler metric $F$.
Clearly, such an isometry  preserves the (possibly nonsymmetric) Finsler distance of each pair
of points of $M$. The converse is the Finslerian version of the Myers-Steenrod theorem,
which states that a distance preserving mapping of $M$ onto itself is an
isometry. This result was proved by Brickell \cite{Bri} and
 Deng and Hou \cite{DenHo}. Recently,  Aradi and Kertesz \cite{ArKe} gave a
  simple proof of it, and also showed that the isometries preserve geodesics.
   As in the Riemannian case it follows from the above result that
   the compact-open topology turns the isometry group $I(M,F)$ of a connected Finsler manifold $(M, F)$
    into a locally compact transformation group of $M$ (and so being a Lie transformation group)
whose isotropy group $I_x(M,F)$ at each $x\in M$ is compact (see \cite{DenHo}).
In particular, if $M$ is compact so is $I(M,F)$.
Suppose that $\mathbb{I}\in{\rm Diff}(M)$ is  contained
in a compact subgroup of ${\rm Diff}(M)$, which implies that $\mathbb{I}$  is homotopic to
an element of ${\rm Diff}(M)$ of finite order.
(In particular, an isometry $\mathfrak{I}\in I(M,F)$  of \textsf{compact type}, i.e.,
the subgroup generated by $\mathfrak{I}$ is relatively compact in  $I(M,F)$, is
such a diffeomorphism.  Observe that an isometry on $(M, F)$ is of compact type if either it has fixed points
or $M$ is compact.)

\begin{theorem}\label{th:1.14}
 Let $\gamma_0\in \Lambda(M,\mathbb{I})$ be a closed
$\mathbb{I}$-invariant $F$-geodesic with ${\cal L}(\gamma_0)=c>0$
such that $\R_{\gamma_0}$ is infinite cyclic with generator $p>0$.
All results in Section~\ref{sec:1.2} hold if
there exists an $\mathbb{I}$-invariant Riemannian metric $g$ on $M$, i.e.,  $\mathbb{I}\in I(M,g)$,
and the associated Hilbert manifolds are endowed with the Hilbert-Riemannian structure induced by this
metric.
\end{theorem}

\subsection{An application}\label{sec:1.4}

For an isometry $\mathbb{I}$ of finite order (resp. any isometry $\mathbb{I}$)
on a simply connected and closed Riemannian mainfold $M$,
 as an extension of a theorem by  Gromoll and Meyer
it was proved by Grove and Tanaka \cite{GroTa78} (resp.  by Tanaka \cite{Tan82}) that
there exist infinitely many distinct $\mathbb{I}$-invariant closed geodesics
if $C^0(I, M)^{\mathbb{I}}:=\{x\in C^0([0,1], M)\,|\, \mathbb{I}(x(0))=x(1)\}$ with the
compact-open topology has an unbounded sequence of Betti numbers.
The following is a generalization of their result on Finsler manifolds.

\begin{theorem}\label{th:1.15}
 Assume that $F$ is a $C^k$ Finsler metric on a connected and closed
manifold $M$ of  dimension at least $2$ for $k\ge 5$. If an isometry $\mathbb{I}\in I(M, F)$ has  at most finitely many $\mathbb{I}$-invariant geodesics on $M$,
and there exists an $\mathbb{I}$-invariant Riemannian metric $g$ on $M$,  then
$\{\dim H_k(C^0(I, M)^{\mathbb{I}}; \K)\,|\,
k\ge 2\dim M\}$ is bounded for any field $\K$.
\end{theorem}

As in the Riemannian case, for example
 \cite{GroHaV78,GroHa82,GroHa91, Tan82, McCZ, VPS76} etc.,
 this theorem may lead to many interesting results.
 We list a few of them.

Because of the compactness of $I(M, F)$ every isometry $\mathbb{I}\in I(M, F)$ is isotopic to
a finite order $F$-isometry which commutes with $\mathbb{I}$ as in the proof of
\cite[Proposition~2.1(iii)]{GroHa91}, and hence rigid at $1$ by \cite[Theorem~3.9]{GroHaV78}.
Then as in \cite[Theorem~2.6]{GroHa91}, the second theorem of \cite{GroHaV78},
 Theorem~\ref{th:1.15} and the theory of minimal models lead to

\begin{corollary}\label{cor:1.16}
Under the assumptions of Theorem~\ref{th:1.15}, if $M$ is also simply connected then
\begin{description}
\item[(i)] $\dim\sum_k(\pi_k(M)\otimes\Q)<\infty$;
\item[(ii)] $\dim(\pi_{\rm even}(M)\otimes\Q)^{\mathbb{I}}\le\dim(\pi_{\rm odd}(M)\otimes\Q)^{\mathbb{I}}\le
1$,  where $(\pi_\ast(M)\otimes\Q)^{\mathbb{I}}$
  denote the $\mathbb{I}$-invariant part of rational homotopy.
 \end{description}
 \end{corollary}

\begin{corollary}\label{cor:1.17}
 Let $M$ be a  simply connected and closed manifold and let $F$ be a
$C^k$ Finsler metric on $M$ with $k\ge 5$. Then  every isometry $\mathbb{I}\in I(M, F)$ has
infinitely many $\mathbb{I}$-invariant geodesics  provided that
  the sequence of Betti numbers of the space $C^0(I, M)^{\mathbb{I}}$ for some field $\K$ is unbounded.
In particular, every
isometry $\mathbb{I}\in I(M, F)$ which is homotopic to ${\rm id}_M$ has infinitely many $\mathbb{I}$-invariant geodesics if one of
the following three equivalent conditions is satisfied:
\begin{description}
\item[(i)] $\sup_j\dim H_j(C^0(S^1, M); \Q)=\infty$;
\item[(ii)] the cohomology algebra $H^\ast(M;\Q)$ requires at least two generators;
\item[(iii)] $\dim(\pi_{\rm odd}(M)\otimes\Q)\ge 2$.
 \end{description}
 \end{corollary}

The part of ``In particular" is straightforward by
 Sullivan-Vigu\'e theorem \cite{VPS76} (cf. \cite[Proposition~5.14]{FOT08}), which covers
 all rational hyperbolic manifolds.
Recall that a simply connected and closed manifold $M$ with $\dim\sum_k(\pi_k(M)\otimes\Q)<\infty$ is called \textsf{rational elliptic},
and \textsf{rational hyperbolic} otherwise,  see \cite{FHT01, FOT08}. Actually,
the second class is vastly bigger than the first; for example,
every simply connected and closed $4$-manifold with second Betti number $k\ge 3$ is rational hyperbolic
\cite[Example~3.8]{FOT08}. But all compact simply connected homogeneous spaces
belong to the first class; in particular, if $M$ is one of  only four simply connected compact symmetric spaces
of rank one, $S^n$, $\CP^n$, $\mathbb{H}P^n$ and ${\C}aP^2$,
$\dim H_j(C^0(S^1, M); \K)$ is always bounded for any field $\K$.
  Fortunately, an important theorem by  McCleary and Ziller \cite{McCZ} and the first part of Corollary~\ref{cor:1.17}
  lead to

\begin{corollary}\label{cor:1.18}
 Let $M$ be a connected closed  manifold which has the same
homotopy type as a compact simply-connected homogeneous
space which is not diffeomorphic to a symmetric space of rank $1$, and let $F$ be a
$C^k$ Finsler metric on $M$ with $k\ge 5$. Then each isometry $\mathbb{I}$  on $(M,F)$ which is homotopic to ${\rm id}_M$ has infinitely many $\mathbb{I}$-invariant geodesics of constant speed.
\end{corollary}

 Very recently,  Jones and McCleary \cite{JoMcC} proved  that each manifold
 in Corollary~\ref{cor:1.18} satisfies: there exists a prime integer $p>0$
 such that  the algebra $H^\ast(M; \F_p)$ cannot be generated by one element and that
  $\sum_{i\le n}\dim H_i(\Omega M;\F_p)\le Cn^K$, $\forall n\in\N$, for some constant $C$ and an integer $K$,
 where  $\F_p$ is the finite field with $p$ elements and $\Omega M$ refers to the based loop space of $M$;
 and their proof of \cite[Theorem~1.2]{JoMcC} also shows that these two conditions  imply that
  the sequence $\dim H_j(C^0(S^1, M); \F_p)$ is unbounded and hence
Corollary~\ref{cor:1.18} holds true on such manifolds.

 Many conclusions on isometric invariant geodesics on Riemannian manifolds
 can be directly generalized to Finsler manifolds,  for instance, Proposition~2.7, Lemma~3.4, Theorem~3.7 and Corollary~3.8, Proposition~3.9 in \cite{Gro73I}, and \cite[Theorem~2.6]{Gro74} (with negative sectional curvature replaced  by negative flag curvature) and \cite[Theorem~B]{GroHa82}.
  Furthermore, it is not too hard to generalize some of results in
\cite{Hin88,Maz1, Maz2,PP05, Rad89} to Finsler manifolds.

Recently, the famous theorem by  Gromoll and Meyer was also generalized
from viewpoint of contact geometry by  Hryniewicz and Macarini \cite{HrMa} using contact homology
 and by Mclean \cite{Mcl} using  symplectic homology, respectively.
It is not  very hard to generalize their methods and to get the corresponding contact analogy of Theorem~\ref{th:1.15} for contactomorphisms of finite order. We shall consider such extensions in a general setting.

\subsection{Concluding remarks}\label{sec:1.5}

The results and methods in this paper may be extended toward different directions, for instance,
pseudo-Finsler manifolds considered in \cite{Jav13,Jav13+, Jav14}, and
 the Lagrange geometry introduced by J. Kern \cite{Ker} and developed
by R. Miron and M. Anastasiei \cite{MirA}. A {\it Lagrange
manifold} $L^n=(M, L(x,y))$ is a pair consisting of a smooth $n$-dimensional
manifold $M$ and a fundamental (or metric) function $L:TM\to\R$ which is a {\it regular Lagrangian}
 in the following sense: a) $L$ is of
class $C^\infty$ on $TM\setminus 0_M$ and continuous on $TM$, b) the matrix
with the entries $g_{ij}(x,y)=\frac{1}{2}\frac{\partial^2 L}{\partial y^i\partial y^j}$
 has rank $n$ on $TM\setminus 0_M$. It is a Finsler manifold if and only
 if its fundamental function $L(x,y)$ is positive and $2$-homogeneous with respect to the fiber
 variable of $TM\setminus 0_M\to M$.  For a closed  Lagrange manifold $L^n=(M, L(x,y))$,
  if $L:TM\to\R$ is of class $C^k$ ($k\ge 5$) and satisfies the conditions (L1)-(L2)
  on $TM\setminus 0_M$ below (\ref{e:2.4}), the method in this paper is still effective
  near every regular solution.  A class of important examples of such manifolds is
  the so-called {\it almost Finslerian Lagrange manifold} (AFL-manifold). This
 is a Lagrange manifold $(M, L(x,y))$ whose fundamental function is given by
$L(x,y)=F^2(x,y)+ A(x)(y)+ U(x)$, where $A$ is a smooth $1$-form on $M$ (called magnetic vector potential)
 and $U$ is a smooth function on $M$ (called potential function).
 For applications of our method to some variational problems of higher dimension, see \cite{Lu4}.

\vspace{2mm}

\noindent{\bf  Organization of this paper}. In Section~\ref{sec:2} we give the
proofs of Theorem~\ref{th:1.1},~\ref{th:1.4}. In Section~\ref{sec:3} we prove Theorem~\ref{th:1.7}.
 The proofs of Theorem~\ref{th:1.12} and Corollary~\ref{cor:1.11} are given in Section~\ref{sec:4}.
  In Section~\ref{sec:5} we give some propositions and the proofs of Lemmas~\ref{lem:2.2},~\ref{lem:2.3}.
Finally,  the proof of Theorem~\ref{th:1.15} is completed in Section~\ref{sec:6}.

 \noindent{\bf Acknowledgements}.  The author thanks an anonymous referee for helpful
comments, pointing out innumerable mistakes, raising several delicate points
 which I had overlooked, and suggesting ways to improve the exposition.


\section{Proofs of Theorem~\ref{th:1.1},~\ref{th:1.4}}\label{sec:2}
\setcounter{equation}{0}

\subsection{Proof of Theorem~\ref{th:1.1}}\label{sec:2.1}
\setcounter{equation}{0}

\subsubsection{The proof in a suitable chart}\label{sec:2.1.1}
 For conveniences of computations
we need to consider a coordinate chart around $\gamma_0$ different
from that of (\ref{e:1.8}), which appeared in \cite{AbF} and were
intensively used in the context of Finsler geometry in
\cite{Ca, CaJaMa2, Lu3}.
Since $\gamma_0$ is $C^k$, for  the
Riemannian metric $g$ on $M$ satisfying (\ref{e:1.7}) we may choose
a parallel orthogonal $C^k$ frame field along $\gamma_0$ with
respect to it, $I\ni t\to (e_1(t),\cdots, e_n(t))$. For some small
open ball ${\bf B}_{2\rho}(\R^n)$ in $\R^n$ centered at $0$ and of
radius $2\rho$ we get a $C^k$ map $\phi:I\times {\bf
B}_{2\rho}(\R^n)\to M$ given by
\begin{equation}\label{e:2.1}
\phi(t,x)=\exp_{\gamma_0(t)}\left(\sum^n_{i=1}x_ie_i(t)\right).
\end{equation}
Since $M_0$ (resp. $M_1$) is totally geodesic near $\gamma_0(0)$
(resp. $\gamma_0(1)$) with respect to $g$ there exist linear
subspaces $V_i\subset\R^n$, $i=0,1$, such that $x\in V_i$ if and
only if $\sum^n_{s=1}x_se_s(i)\in T_{\gamma_0(i)}M_i$, $i=0,1$. By
shrinking $\rho>0$ we get that $x\in V_i\cap {\bf
B}_{2\rho}(\R^n)$ if and only if $ \phi(i,x)\in M_i$, $i=0,1$. Set
$V:=V_0\times V_1$ and
\begin{eqnarray*}
&&H_V=W^{1,2}_{V}(I, \R^n):=\{\xi\in W^{1,2}(I, \R^n)\,|\,
(\xi(0),\xi(1))\in V\},\\
&& X_V=C^1_{V}(I, \R^n):=\{\xi\in C^1(I, \R^n)\,|\,
(\xi(0),\xi(1))\in V\}.
\end{eqnarray*}
They are equipped with norms
$\|\xi\|_{H_V}=\|\xi\|_1$ and $\|\xi\|_{X_V}=\|\xi\|_{\infty}+\|\dot{\xi}\|_{\infty}$, respectively.
Then $\|\xi\|_{H_V}\le\|\xi\|_{X_V}\;\forall\xi\in X_V$.
Let ${\bf B}_{2\rho}(H_V):=\{\xi\in H_V\,|\,\|\xi\|_{H_V}<2\rho\}$.
 For $\rho>0$ small enough   the map
\begin{equation}\label{e:2.2}
\Phi:{\bf B}_{2\rho}(H_V) \to\Lambda_N(M)
\end{equation}
defined by $\Phi(\xi)(t)=\phi(t,\xi(t))$, gives a $C^{k-3}$ coordinate
chart around $\gamma_0$ on $\Lambda_N(M)$ by \cite[Theorem~4.3]{PiTa}, and we can require that
all curves in ${\rm Im}(\Phi)$ have images contained in a compact neighborhood
of $\gamma_0([0,1])$.  Define $\tilde F:I\times {\bf
B}_{2\rho}(\R^n)\times\R^n\to\R$ by
$$
\tilde F(t, x, v)=F\bigr(\phi(t,x), d\phi(t,x)[(1,v)]\bigl).
$$
We have $\tilde F(t, 0, 0)=F(\phi(t,0),
\partial_t\phi(t, 0)[1])=F(\gamma_0(t),\dot\gamma_0(t))\equiv\sqrt{c}$.
Moreover the $C^{2-0}$ function $\tilde L:=\tilde F^2$ satisfies
$\tilde L(t, x, v)=L\bigr(\phi(t,x), d\phi(t,x)[(1,v)]\bigl)$, and it
is $C^k$ in $(I\times {\bf B}_{2\rho}(\R^n)\times\R^n)\setminus{\cal
Z}$, where
\begin{equation}\label{e:2.3}
{\cal Z}:=\bigl\{(t,x,v)\in I\times {\bf
B}_{2\rho}(\R^n)\times\R^n\,|\,
\partial_x\phi(t,x)[v]=-\partial_t\phi(t,x)\bigr\},
\end{equation}
a  closed subset in $I\times {\bf
B}_{2\rho}(\R^n)\times\R^n$, see \cite[p. 861]{CaJaMa2}. Furthermore, we have

\begin{claim}\label{cl:2.1}
The set ${\cal Z}$ is  a
submanifold of dimension $n+1$.
\end{claim}

\noindent{\bf Proof}.  Let $\pi_2:I\times M\to M$ be the natural projection.
Since the map $\phi$
in (\ref{e:2.1}) yields an embedding of codimension zero
$$
\overline{\phi}:I\times {\bf B}_{2\rho}(\R^n)\to I\times M
$$
given by $\overline{\phi}(t,x)=(t, \phi(t, x))$, we obtain a bundle
embedding of codimension zero covering $\overline{\phi}$,
\begin{equation}\label{e:2.4}
\overline{\phi}^\star:I\times {\bf
B}_{2\rho}(\R^n)\times\R^n\to \pi_2^\ast TM
\end{equation}
 given by
$\overline{\phi}^\star(t,x,v)=d\phi(t,x)[(1,v)]$. Observe that ${\cal Z}$ is
the inverse image of the zero section $0_{\pi_2^\ast TM}$ of
$\pi_2^\ast TM$ under $\overline{\phi}^\star$, i.e., ${\cal
Z}=(\overline{\phi}^\star)^{-1}(0_{\pi_2^\ast TM})$, which is not
necessarily the zero section of the trivial bundle $I\times {\bf
B}_{2\rho}(\R^n)\times\R^n\to I\times {\bf B}_{2\rho}(\R^n)$.
The claim follows.
\hfill$\Box$\vspace{2mm}

Since $F^2$ is fiberwise positively homogeneous of degree $2$ and
strongly convex there exist positive constants $C_1<C_2$  such that
for all $(t,x,v)\in I\times {\bf
B}_{2\rho}(\R^n)\times\R^n\setminus{\cal Z}$,
\begin{description}
\item[(L1)] $\sum_{ij}\frac{\partial^2}{\partial v_i\partial
v_j}\tilde{L}(t,x,v)u_iu_j\ge C_1|{\bf u}|^2\quad\forall {\bf
u}=(u_1,\cdots,u_n)\in\R^n$,

\item[(L2)] $\Bigl| \frac{\partial^2}{\partial x_i\partial
x_j}\tilde{L}(t,x,v)\Bigr|\le C_2(1+ |v|^2),\quad \Bigl|
\frac{\partial^2}{\partial x_i\partial v_j}\tilde{L}(t,x,v)\Bigr|\le
C_2(1+
|v|),\quad\hbox{and}\\
 \Bigl| \frac{\partial^2}{\partial v_i\partial
v_j}\tilde{L}(t,x,v)\Bigr|\le C_2$
\end{description}
(cf. \cite[(1) $\&$ (2)]{CaJaMa2}).

For every $m\in\N$, $C^m_N(I, M):=\{\gamma\in C^m(I,M)\,|\,
(\gamma(0), \gamma(1))\in N\}$ is a Banach manifold and dense in
${\cal X}$ and $\Lambda_N(M)$. In order to distinguish from the regularity of the curve
$\gamma\in\Lambda_N(M)$ defined below (\ref{e:1.2}),  a curve $\gamma\in C^m_N(I, M)$
is called {\bf strongly regular} if $\{t\in I\,|\,\dot\gamma(t)=0\}$ is empty.
Let $C^m_N(I, M)_{\rm reg}$ denote the subset of all strongly regular
curves in $C^m_N(I, M)$. Let
$C^m_{V}(I, \R^n):=\{\xi\in C^m(I,
\R^n)\,|\, (\xi(0),\xi(1))\in V\}$ and
\begin{eqnarray*}
C^m_{V}(I, {\bf B}_{2\rho}(\R^n)):=\{\xi\in C^m(I, {\bf
B}_{2\rho}(\R^n))\,|\, (\xi(0),\xi(1))\in V\}.
\end{eqnarray*}
Note that every  $\xi\in(\Phi|_{{\bf B}_{\sqrt{2}\rho}(H_V)})^{-1}(C^m_N(I,
M)_{\rm reg})$ sits in ${\bf B}_{\sqrt{2}\rho}(H_V)\cap
C^m_{V}(I, {\bf B}_{2\rho}(\R^n))$ because
$\|\xi\|_\infty\le\sqrt{2}\|\xi\|_{1}<2\rho$ by the integral mean
value theorem and the H\"older inequality. ({\it Here $\|\xi\|_1:=\sqrt{\|\xi\|^2_{L^2}+\|\dot{\xi}\|^2_{L^2}}$.
If $\xi\|_1:=\|\xi\|_{L^2}+\|\dot{\xi}\|_{L^2}$ it holds that $\|\xi\|_\infty\le\|\xi\|_1$}).
Hence it determines a map
\begin{equation}\label{e:2.5}
\check{\xi}:I\to I\times {\bf B}_{2\rho}(\R^n)\times\R^n,\quad
t\mapsto\check{\xi}(t)=(t,\xi(t),\dot\xi(t)).
\end{equation}
Clearly, $\check{\xi}\in C^{m-1}(I,I\times {\bf
B}_{2\rho}(\R^n)\times\R^n)$ and with $\gamma=\Phi(\xi)$ it holds
that
$$
 (\check{\xi})^{-1}({\cal
Z})=\{t\in I\,|\,(t, \xi(t ),\dot\xi(t ))\in {\cal Z}\}=\{t\in
I\,|\,\dot\gamma(t)=0\}.
$$
 This shows that
\begin{equation}\label{e:2.6}
(\Phi|_{{\bf B}_{\sqrt{2}\rho}(H_V)})^{-1}(C^m_N(I, M)_{\rm reg})=
C^m_{V}(I, {\bf B}_{2\rho}(\R^n))_{\rm reg}\cap {\bf
B}_{\sqrt{2}\rho}(H_V),
\end{equation}
where $C^m_{V}(I, {\bf B}_{2\rho}(\R^n))_{\rm reg}$ is the set of
all $\xi\in C^m_{V}(I, {\bf B}_{2\rho}(\R^n))$ such that
$(\check{\xi})^{-1}({\cal Z})=\{t\in I\,|\,(t, \xi(t ),\dot\xi(t
))\in {\cal Z}\}$ is empty, i.e., $\xi$ is an immersion. If $n>1$,
by Theorem~2.12 in \cite[page 53]{Hi} the set of all $C^m$
immersions from $I$ to ${\bf B}_{2\rho}(\R^n)$ is dense in $C^m(I,
{\bf B}_{2\rho}(\R^n))$. However we need a stronger result.

 Taking $m=k+1$ we have  $C^k$ maps
\begin{equation}\label{e:2.7}
{\bf E}:I\times C^{k+1}_{V}(I, {\bf B}_{2\rho}(\R^n))\to I\times
{\bf B}_{2\rho}(\R^n)\times\R^n,\;(t,\xi)\mapsto \check{\xi}(t)
\end{equation}
and ${\bf E}_\xi:I\to I\times {\bf B}_{2\rho}(\R^n)\times\R^n$ given
by ${\bf E}_\xi(t)={\bf E}(t,\xi)$ for $\xi\in C^{k+1}_{V}(I, {\bf
B}_{2\rho}(\R^n))$.

\begin{lemma}\label{lem:2.2}
 There exists a residual  subset
$C^{k+1}_{V}(I, {\bf B}_{2\rho}(\R^n))^\circ_{\rm reg}$  in
$C^{k+1}_{V}(I, {\bf B}_{2\rho}(\R^n))$, which contains
$C^{k+1}_{V}(I, {\bf B}_{2\rho}(\R^n))_{\rm reg}$ and is equal to
the latter  in case $n>1$, such that for every $\xi\in
C^{k+1}_{V}(I, {\bf B}_{2\rho}(\R^n))^\circ_{\rm reg}$ the set
$({\bf E}_\xi)^{-1}({\cal Z})$ is a $C^k$ submanifold of dimension
$1-n$, and so empty (in case $n>1$) and at most a finite set (in
case $n=1$).
 Moreover, if  $\xi_i\in C^{k+1}_{V}(I, {\bf
B}_{2\rho}(\R^n))^\circ_{\rm reg}$, $i=1,2$, then for a generic
$C^{k+1}$ path $\mathbbm{p}:[0, 1]\to C^{k+1}_{V}(I, {\bf
B}_{2\rho}(\R^n))$ connecting $\xi_0$ to $\xi_1$ the set $\{(s,t)\in
[0, 1]\times I\,|\, {\bf E}(t,\mathbbm{p}(s))\in{\cal Z}\}$ is a
$C^k$ submanifold of $[0, 1]\times I$ of dimension $2-n$, and so
empty for $n>2$, and a finite set set for $n=2$.
\end{lemma}

We postpone the proof of it  to Section~\ref{sec:5}.

 Since $\gamma_0$ is strongly regular, i.e.,
$\partial_t\phi(t,0)=\dot\gamma_0(t)\ne 0$ at each $t\in I$, and
$\partial_x\phi(t,x)$ is injective, we deduce that $(t, 0,
0)\notin{\cal Z}\;\forall t\in I$.  It follows that
\begin{equation}\label{e:2.8}
I\times {\bf B}_{2r}(\R^n)\times {\bf B}_{2r}(\R^n)\subset (I\times
{\bf B}_{2\rho}(\R^n)\times\R^n)\setminus{\cal Z}
\end{equation}
for some $0<r<\rho/2$.  Let ${\bf B}_{\tau}(X_V)=\{\xi\in
X_V\,|\, \|\xi\|_{X_V}<\tau\}$ for $\tau>0$. Then
\begin{equation}\label{e:2.9}
{\bf B}_{2r}(X_V)\subset {\bf B}_{2r}(H_V),\quad (\check{\xi})^{-1}({\cal
Z})=\emptyset\;\forall \xi\in {\bf
B}_{2r}(X_V)
\end{equation}
by (\ref{e:2.8}).  Define
the action functional
\begin{equation}\label{e:2.10}
\tilde{\cal L}: {\bf B}_{2r}(H_V)\to\R,\;\xi\mapsto \tilde{\cal
L}(\xi)=\int^1_0\tilde L(t,\xi(t),\dot\xi(t))dt,
\end{equation}
that is, $\tilde{\cal L}={\cal L}\circ\Phi$. Since $\Phi$ is of class $C^{k-3}$, $\tilde{\cal L}$
 is $C^{2-0}$ for $k\ge 5$, and
has $0\in {\bf B}_{2r}(H_V)$ as a
 critical point. (But $\tilde{\cal L}$ is $C^{k-2}$ on ${\bf B}_{2r}(X_V)$
 by the last claim of Lemma~\ref{lem:2.3}).
 For any $\zeta\in {\bf B}_{2r}(H_V)$ and $\xi\in
 H_V$, we have
\begin{eqnarray*}
d\tilde{\cal L}(\zeta)[\xi]
 \!\!\!&=&\!\!\!\int_0^{1}  \bigl[\partial_{q}\tilde
L(t,\zeta(t),\dot{\zeta}(t))\cdot\xi(t) + \partial_{v} \tilde
L(t,\zeta(t),\dot{\zeta}(t))\cdot\dot{\xi}(t) \bigr] dt
\end{eqnarray*}
(cf. \cite[\S3]{Lu1}).
  As in \cite[(4.14)]{Lu3} we can compute
 the gradient
\begin{eqnarray}\label{e:2.11}
\nabla\tilde{\cal L}(\zeta)(t)&=&e^t\int^t_0\left[
e^{-2s}\int^s_0e^{\tau}f(\tau)d\tau\right]ds  + c_1e^t+
c_2e^{-t}\nonumber\\
&&\qquad +\int^t_0 \partial_{v}\tilde{L}(s,\zeta(s),\dot{\zeta}(s))ds,
\end{eqnarray}
where $c_1, c_2\in\R^n$ are suitable constant vectors and
\begin{eqnarray}\label{e:2.12}
 f(t)&=&- \partial_{q} \tilde
L(t,\zeta(t),\dot{\zeta}(t))+
G(\zeta)(t)+ c_0t\nonumber\\
&=& - \partial_{q} \tilde L(t,\zeta(t),\dot{\zeta}(t))+ \int^t_0 \partial_{v}
\tilde L(s,\zeta(s),\dot{\zeta}(s))ds
\end{eqnarray}
because  $G(\zeta)(t):=\int^t_0\big[\partial_{v}
\tilde L(s,\zeta(s),\dot{\zeta}(s))ds-c_0\big]dt$.

\textsf{From now on we take $m=3$.} Then $\min\{k-3,m-1\}>1$ because
$k\ge 5$.  Let $U={\bf B}_{2r}(H_V)$, $U_X:={\bf B}_{2r}(H_V)\cap
X_V$  as an open subset of $X_V$ and
\begin{eqnarray}\label{e:2.13}
U_X^{\rm reg}:=U_X\cap C^3_{V}(I, {\bf B}_{2\rho}(\R^n))_{\rm
reg}={\bf B}_{2r}(H_V)\cap C^3_{V}(I, {\bf B}_{2\rho}(\R^n))_{\rm
reg}.
\end{eqnarray}
Clearly, ${\bf B}_{2r}(X_V)\subset U_X$.  By (\ref{e:2.6}), $U_X^{\rm reg}=(\Phi|_{{\bf B}_{2r}(H_V)})^{-1}(C^3_N(I, M)_{\rm
reg})$  since $2r<\rho<\sqrt{2}\rho$.
 Lemma~\ref{lem:2.2} implies that $U_X^{\rm reg}$  is
dense in $U_X$
 (and hence in $U$). Let $\tilde{\cal L}_X$ be the restriction of $\tilde{\cal L}$ to $U_X$.
 By (\ref{e:2.9}) we deduce that $\tilde{\cal L}_X$
 has with respect to the $C^1$-topology second Frech\'et derivative at each
$\zeta\in U_X^{\rm reg}\cup{\bf B}_{2r}(X_V)$,
\begin{eqnarray}\label{e:2.14}
 d^2 \tilde{\cal L}_X
  (\zeta)[\xi,\eta]
   = \int_0^{1} \Bigl(\!\! \!\!\!&&\!\!\!\!\!\partial_{vv}
    \tilde L\bigl(t,\zeta(t),\dot{\zeta}(t)\bigr)
\bigl[\dot{\xi}(t), \dot{\eta}(t)\bigr] \nonumber\\
&&+ \partial_{qv} \tilde
  L\bigl(t,\zeta(t), \dot{\zeta}(t)\bigr)
\bigl[\xi(t), \dot{\eta}(t)\bigr]\nonumber \\
&& + \partial_{vq} \tilde
  L\bigl(t,\zeta(t),\dot{\zeta}(t)\bigr)
\bigl[\dot{\xi}(t), \eta(t)\bigr] \nonumber\\
&&+  \partial_{qq} \tilde L\bigl(t,\zeta(t), \dot{\zeta}(t)\bigr)
\bigl[\xi(t), \eta(t)\bigr]\Bigr) \, dt
\end{eqnarray}
for any $\xi,\eta\in X_V$.
  Moreover,
 for every $\zeta\in U_X^{\rm reg}\cup{\bf B}_{2r}(X_V)$ it is easily seen
 from (\ref{e:2.11}) that the function
$\nabla\tilde{\cal L}(\zeta)(t)$ is continuous differentiable, and
\begin{eqnarray}\label{e:2.15}
\frac{d}{dt}\nabla\tilde{\cal L}(\zeta)(t)&=& e^t\int^t_0\left[
e^{-2s}\int^s_0e^{\tau}f(\tau)d\tau\right]ds
+e^{-t}\int^t_0e^{\tau}f(\tau)d\tau
\nonumber\\
&&\quad + c_1e^t -c_2e^{-t}+ \partial_{v} \tilde{L}\bigl(t,\zeta(t),\dot{\zeta}(t)\bigr)\quad\forall t\in I.
\end{eqnarray}
Let  $\tilde A$ be the restriction of the gradient
$\nabla\tilde{\cal L}$ to $U_X$.
By (\ref{e:2.11}) and (\ref{e:2.15}) we can  deduce that
$\tilde A(\zeta)\in X_V$ for  $\zeta\in
U_X^{\rm reg}\cup{\bf B}_{2r}(X_V)\subset U_X$, and that
$U_X\supseteq U_X^{\rm reg}\cup{\bf B}_{2r}(X_V)\ni\zeta\mapsto \tilde A(\zeta)\in X_{V}$ is continuous.
Furthermore we shall prove in Section~\ref{sec:5}.

\begin{lemma}\label{lem:2.3}
 The map $U_X\supseteq U_X^{\rm reg}\cup{\bf B}_{2r}(X_V)\ni\zeta\mapsto \tilde A(\zeta)\in X_{V}$ is Frech\'et differentiable, and $U_X^{\rm reg}\cup{\bf B}_{2r}(X_V)\ni\zeta\to d\tilde{A}(\zeta)\in\mathscr{L}(X_V)$ is continuous. In particular, $\tilde{A}$ restricts
to a $C^{k-3}$ map from ${\bf B}_{2r}(X_V)$ to $X_V$.
\end{lemma}

From (\ref{e:2.14}) and the conditions ({\bf L1}) and ({\bf L2}) below
(\ref{e:2.4}) it easily follows that
for any $\zeta\in U_X^{\rm reg}\cup{\bf B}_{2r}(X_V)$ there
exists a constant $C(\zeta)>0$ such that
$$
|d^2 {\tilde{\cal L}}_{X}
  (\zeta)[\xi,\eta]|\le
  C(\zeta)\|\xi\|_{H_V}\cdot\|\eta\|_{H_V}
\quad\forall \xi,\eta\in X_V;
$$
This shows that the right side of (\ref{e:2.14}) is also a bounded symmetric
bilinear form on $H_V$. As in \cite[\S3]{Lu1}, from these we obtain
a map $\tilde B: U_X^{\rm reg}\cup{\bf B}_{2r}(X_V)\to \mathscr{L}_s(H_V)$ such that
\begin{equation}\label{e:2.16}
\bigl(d\tilde A(\zeta)[\xi], \eta\bigr)_{H_V}=d^2 {\tilde{\cal L}}_{X}
  (\zeta)[\xi,\eta]=\bigl(\tilde B(\zeta)\xi,
  \eta\bigr)_{H_V}
\end{equation}
for  any $\zeta\in U_X^{\rm reg}\cup{\bf B}_{2r}(X_V)$ and $\xi, \eta\in X_V$. Moreover,
$\tilde{B}(0)(X_V)\subset X_V$, $(\tilde{B}(0))^{-1}(X_V)\subset X_V$, and
$\tilde B(0)\in\mathscr{L}_s(H_V)$ is a  Fredholm operator
with  finite dimensional negative definite and null spaces $H^-_V$ and $H^0_V$ (which are
 contained in $X_V$). The map
$\tilde B: U_X^{\rm reg}\cup{\bf B}_{2r}(X_V)\to \mathscr{L}_s(H_V)$
 has also a decomposition
$$
\tilde B(\zeta)=\tilde P(\zeta)+ \tilde Q(\zeta)\quad\forall
\zeta\in U_X^{\rm reg}\cup{\bf B}_{2r}(X_V),
$$
where $\tilde P(\zeta)\in\mathscr{L}_s(H_V)$ is a positive definitive
linear operator defined by
\begin{eqnarray*}
(\tilde{P}(\zeta)\xi, \eta)_{H_V}
   = \int_0^{1} \Bigl(\!\! \!\!\!&&\!\!\!\!\!\partial_{vv}
    \tilde L\bigl(t,\zeta(t),\dot{\zeta}(t)\bigr)
\bigl[\dot{\xi}(t), \dot{\eta}(t)\bigr]+
\bigl(\xi(t),\eta(t)\bigr)\Bigr) \, dt,
\end{eqnarray*}
 and $\tilde Q(\zeta)\in\mathscr{L}_s(H_V)$ is a compact
linear operator given by
\begin{eqnarray*}
(\tilde{Q}(\zeta)\xi, \eta)_{H_V}
    &=& \int_0^{1} \Bigl( \partial_{qv} \tilde
  L\bigl(t,\zeta(t), \dot{\zeta}(t)\bigr)
\bigl[\xi(t), \dot{\eta}(t)\bigr] + \partial_{vq} \tilde
  L\bigl(t,\zeta(t),\dot{\zeta}(t)\bigr)
\bigl[\dot{\xi}(t), \eta(t)\bigr] \nonumber\\
&&\hspace{10mm}+  \partial_{qq} \tilde L\bigl(t,\zeta(t), \dot{\zeta}(t)\bigr)
\bigl[\xi(t), \eta(t)\bigr]-\bigl(\xi(t), \eta(t)\bigr)\Bigr) \, dt.
\end{eqnarray*}

\begin{lemma}\label{lem:2.3+}
The operators $\tilde{P}(\zeta)$ and $\tilde{Q}(\zeta)$ have the following properties:
\begin{description}
\item[(P1)] For any sequence $(\zeta_k)\subset U_X^{\rm reg}\cup{\bf B}_{2r}(X_V)$ with $\|\zeta_k\|_{H_V}\to 0$
it holds that $\|\tilde P(\zeta_k)\xi-\tilde P(0)\xi\|_{H_V}\to 0$
for any $\xi\in H_V$;

\item[(P2)] For any sequence $(\zeta_k)\subset U_X^{\rm reg}\cup{\bf B}_{2r}(X_V)$ with $\|\zeta_k\|_{H_V}\to 0$
 there exist constants $C_0>0$ and $k_0\in\N$ such that
$(\tilde P(\zeta_k)\xi, \xi)_{H_V}\ge C_0\|\xi\|^2_{H_V}$ for all
$\xi\in H_V,\;k\ge n_0$;

\item[(Q)] The  map $\tilde Q:U_X^{\rm reg}\cup{\bf B}_{2r}(X_V)\to
\mathscr{L}_s(H_V)$ is continuous at $0$ with respect to the topology
induced from $H_V$ on $U_X^{\rm reg}\cup{\bf B}_{2r}(X_V)$.
\end{description}
\end{lemma}

\noindent{\bf Proof}.  The proof is the same as that of \cite[pages 568-569]{Lu1}.
We only prove ({\bf P1}) for the sake of clearness. ({\bf P2})
directly follows from ({\bf P1}) and positivity of $\tilde{P}(0)$ by ({\bf L2}).

For any $\xi,\eta\in H_V$ we have
\begin{eqnarray*}
&&(\tilde{P}(\zeta_k)\xi- \tilde{P}(0)\xi,
  \eta)_{H_V}\\
&=&\int_0^{1} \Bigl(\partial_{vv}
    \tilde L\bigl(t,\zeta_k(t),\dot{\zeta}_k(t)\bigr)
\bigl[\dot{\xi}(t), \dot{\eta}(t)\bigr]- \partial_{vv}
    \tilde L\left(t,0, 0\right)
\bigl[\dot{\xi}(t), \dot{\eta}(t)\bigr]\Bigr) \, dt,\\
&=&\int_0^{1}\sum^n_{j=1}\left[\sum^n_{i=1}\left(\frac{\partial^2\tilde
L}{\partial v_i\partial
v_j}\left(t,\zeta_k(t),\dot{\zeta}_k(t)\right)-
\frac{\partial^2\tilde L}{\partial v_i\partial
v_j}\left(t,0, 0\right)\right)
\dot{\xi}_i(t)\right]\cdot\dot{\eta}_j(t) \, dt,
\end{eqnarray*}
and hence
\begin{eqnarray*}
&&\|\tilde{P}(\zeta_k)\xi- \tilde{P}(0)\xi\|_{H_V}\\
&\le
&\sqrt{n}\left(\int_0^{1}\sum^n_{j=1}\left|\sum^n_{i=1}\left(\frac{\partial^2\tilde
L}{\partial v_i\partial
v_j}\left(t,\zeta_k(t),\dot{\zeta}_k(t)\right)-
\frac{\partial^2\tilde L}{\partial v_i\partial
v_j}\left(t, 0, 0\right)\right)
\dot{\xi}_i(t)\right|^2\, dt\right)^{1/2}.
\end{eqnarray*}
So it suffices to prove that for each fixed pair $(i,j)$,
\begin{eqnarray*}
\int_0^{1}\left|\left(\frac{\partial^2\tilde
L}{\partial v_i\partial
v_j}\left(t,\zeta_k(t),\dot{\zeta}_k(t)\right)-
\frac{\partial^2\tilde L}{\partial v_i\partial
v_j}\left(t, 0, 0\right)\right)
\dot{\xi}_i(t)\right|^2\, dt\to 0
\end{eqnarray*}
as $k\to\infty$.  By a
contradiction suppose that there exist $c_0>0$ and a subsequence of $(\zeta_k)$, still denoted by $(\zeta_k)$,
such that
$$
\int_0^{1}\left|\left(\frac{\partial^2\tilde
L}{\partial v_i\partial
v_j}\left(t,\zeta_k(t),\dot{\zeta}_k(t)\right)-
\frac{\partial^2\tilde L}{\partial v_i\partial
v_j}\left(t, 0, 0\right)\right)
\dot{\xi}_i(t)\right|^2 dt
\ge c_0,\quad\forall k=1,2,\cdots.
$$
Since $\|\zeta_k\|_{H_V}\to 0$, $(\zeta_k)$ converges uniformly to zero on $[0,1]$ by
 \cite[Prop.1.2]{MW}, and there exist a subsequence $(\zeta_{k_s})$ and a function $h\in L^2[0,1]$
 such that $\dot{\zeta}_{k_s}(t)\to 0$ a.e. on $[0,1]$, and $|\dot{\zeta}_{k_s}(t)|\le h(t)\;\forall s$ a.e. on $[0,1]$
  by \cite[Th.4.9]{Bre}. From this we deduce
  $$
  \left|\left(\frac{\partial^2\tilde
L}{\partial v_i\partial
v_j}\left(t,\zeta_{k_s}(t),\dot{\zeta}_{k_s}(t)\right)-
\frac{\partial^2\tilde L}{\partial v_i\partial
v_j}\left(t, 0, 0\right)\right)
\dot{\xi}_i(t)\right|^2\to 0\quad\hbox{a.e. on $[0,1]$}
$$
as $s\to \infty$. Moreover,  ({\bf L2}) below (\ref{e:2.4}) implies
$$
\left|\frac{\partial^2\tilde
L}{\partial v_i\partial
v_j}\left(t,\zeta_k(t),\dot{\zeta}_k(t)\right)-
\frac{\partial^2\tilde L}{\partial v_i\partial
v_j}\left(t, 0, 0\right)\right|\le 2C_2,\quad\forall t\in [0,1],\;\forall k\in\N.
$$
The dominated convergence theorem leads to
$$
\int_0^{1}\left|\left(\frac{\partial^2\tilde
L}{\partial v_i\partial
v_j}\left(t,\zeta_{k_s}(t),\dot{\zeta}_{k_s}(t)\right)-
\frac{\partial^2\tilde L}{\partial v_i\partial
v_j}\left(t, 0, 0\right)\right)
\dot{\xi}_i(t)\right|^2 dt
\to 0
$$
as $s\to\infty$.
 This contradiction affirms ({\bf P1}).
\hfill$\Box$\vspace{2mm}

Let $H^+_V$ be the positive definite space of $\tilde B(0)$. By  \cite{Lu1} there exists  $a_0>0$ such that
\begin{equation}\label{e:2.17}
(\tilde B(0)\xi, \xi)_{H_V}\ge
2a_0\|\xi\|^2_{H_V}\;\forall \xi\in H_V^+,\quad
(\tilde B(0)\xi, \xi)_{H_V}\le -2a_0\|\xi\|^2_{H_V}\;\forall
\xi\in H_V^-.
 \end{equation}
 As in the proof of \cite[Theorem~1.1]{Lu1} (or
 in the proof of Lemmas~3.3, 3.4 in \cite{Lu2}), we can use these to
 prove the following lemmas.

\begin{lemma}\label{lem:2.4}
 There exists a function $\omega:U_X^{\rm reg}\cup{\bf B}_{2r}(X_V)\to [0, \infty)$  such that
 $\omega(\zeta)\to 0$ as $\zeta\in U_X^{\rm reg}\cup{\bf B}_{2r}(X_V)$ and $\|\zeta\|_{H_V}\to
0$, and that
$$
|(\tilde B(\zeta)\xi, \eta)_{H_V}- (\tilde B(0)\xi, \eta)_{H_V} |\le
\omega(\zeta) \|\xi\|_{H_V}\cdot\|\eta\|_{H_V}
$$
for any $\zeta\in U_X^{\rm reg}\cup{\bf B}_{2r}(X_V)$,  $\xi\in H^0_V\oplus H^-_V$ and
$\eta\in H_V$.
\end{lemma}

\begin{lemma}\label{lem:2.5}
By shrinking $U$ (or equivalently $r>0$ in (\ref{e:2.8})), we can
obtain  a number $a_1\in (0, 2a_0]$ such that for any $\zeta\in
U_X^{\rm reg}\cup{\bf B}_{2r}(X_V)$, $\omega(\zeta)<\min\{a_0, a_1\}/2$ and
\begin{description}
\item[(i)] $(\tilde B(\zeta)\xi, \xi)_{H_V}\ge a_1\|\xi\|^2_{H_V}\;\forall \xi\in H^+_V$;
\item[(ii)] $|(\tilde B(\zeta)\xi,\eta)_{H_V}|\le\omega(\zeta)\|\xi\|_{H_V}\cdot\|\eta\|_{H_V}\;
\forall \xi\in H^+_V, \forall \eta\in H^-_V\oplus H^0_V$;
\item[(iii)] $(\tilde B(\zeta)\xi,\xi)_{H_V}\le-a_0\|\xi\|^2_{H_V}\;\forall \xi\in H^-_V$.
\end{description}
\end{lemma}

\textsf{From now on we assume  $\dim H_V^0>0$} (because the
nondegenerate case is simpler). Consider the $C^{k-3}$ map
 $\tilde A:{\bf B}_{2r}(X_V)\to X_V$. We have known that $d\tilde{A}(0)=\tilde{B}(0)|_{X_V}$
has kernel $X^0_V=H^0_V={\rm Ker}(\tilde B(0))$ as sets.
Observe that the orthogonal decomposition
$H_V=H^-_V\oplus H^0_V\oplus H^+_V$ induces
  a (topological) direct sum decomposition $X_V=X^-_V\dot{+} X^0_V\dot{+}
X^+_V$, where $X^-_V=H^-_V$ and $X^+_V=X_V\cap H^+_V$, and that $H_V$ and $X_V$
induce equivalent norms on $H_V^0=X_V^0$.
By the implicit function theorem we
get a $\tau\in (0, r]$ and a $C^{k-3}$-map $\tilde h: \overline{{\bf
B}_\tau(H_V^0)}\to X^-_V\dot{+}X^+_V$ with $\tilde h(0)=0$ and
$d\tilde h(0)=0$ such that for each $\xi\in \overline{{\bf
B}_\tau(H^0_V)}$,
\begin{equation}\label{e:2.18}
\xi+ \tilde h(\xi)\in {\bf B}_{r}(X_V)\quad\hbox{and}\quad
(I-P^0_V)\tilde A(\xi+ \tilde h(\xi))=0,
\end{equation}
where $\overline{{\bf
B}_\tau(H_V^0)}$ denotes the closure of ${\bf
B}_\tau(H_V^0)$ and $P^0_V:H_V\to H^0_V$ is the orthogonal projection. It is not
hard to prove that the  Morse index $m^-({\cal L}, \gamma_0)$ and nullity
$m^0({\cal L}, \gamma_0)$ of $\gamma_0$ are equal to  $\dim X^-_V$ and $\dim
X^0_V$, respectively. (See (\ref{e:2.28})-(\ref{e:2.31}) in Section~\ref{sec:2.1.2}).

Since the map $\mathfrak{M}:\overline{{\bf B}_\tau(H^0_V)}\oplus
H_V^\pm\to H^0_V\oplus H_V^\pm,\;\xi+\zeta\mapsto \xi+ \tilde
h(\xi)+\zeta$, satisfies  $d\mathfrak{M}(0)={\rm id}$, we may shrink the above $\tau>0$ such that
$$
\mathfrak{M}\left(\overline{{\bf B}_\tau(H^0_V)}\oplus {\bf
B}_\tau(H^-_V)\oplus{\bf B}_\tau(H^+_V)\right)\subset {\bf
B}_{2r}(H_V)=U
$$
and that $\mathfrak{M}$ is a $C^{k-3}$ diffeomorphism from
$\overline{{\bf B}_\tau(H^0_V)}\oplus {\bf B}_\tau(H^-_V)\oplus{\bf
B}_\tau(H^+_V)$ onto its image. Recall that
$U_X^{\rm reg}$ is dense in
$U_X$ (and hence in $U$) by (\ref{e:2.13}) and Lemma~\ref{lem:2.2}. We deduce that
$W_X^{\rm reg}:=\mathfrak{M}^{-1}(U_X^{\rm reg})$
is a dense subset in $\overline{{\bf B}_\tau(H^0_V)}\oplus {\bf
B}_\tau(H^-_V)\oplus{\bf B}_\tau(H^+_V)$.

 Consider the functional ${\cal J}:\overline{{\bf B}_\tau(H^0_V)}\times \bigl({\bf B}_\tau(H^-_V)\oplus
 {\bf B}_\tau(H^+_V)\bigr)\to\R$ defined by
\begin{equation}\label{e:2.19}
{\cal J}(\xi, \zeta)=\tilde{\cal L}(\xi+ \tilde h(\xi)+
\zeta)-\tilde{\cal L}(\xi+ \tilde h(\xi))
\end{equation}
for $\xi\in \overline{{\bf B}_\tau(H^0_V)}$ and $\zeta\in {\bf
B}_\tau(H^-_V)\oplus {\bf B}_\tau(H^+_V)$. It is $C^{2-0}$,  and  it
holds that
\begin{eqnarray}
&&D_2{\cal J}(\xi,\zeta)[\eta] =((I-P^0_V)\nabla\tilde{\cal L}(\xi+ \tilde h(\xi)+ \zeta), \eta)_{H_V},\label{e:2.20}\\
&& {\cal J}(\xi, 0)=0\quad\hbox{and}\quad D_2{\cal J}(\xi,
0)[\eta]=0\label{e:2.21}
\end{eqnarray}
for every $(\xi,\zeta)\in \overline{{\bf B}_\tau(H^0_V)}\times \bigl({\bf
B}_\tau(H^-_V)\oplus
 {\bf B}_\tau(H^+_V)\bigr)$ and $\eta\in H^\pm_V$.
 We can require that $\overline{{\bf B}_\tau(H^0_V)}+ \tilde{h}(\overline{{\bf B}_\tau(H^0_V)})$
 is contained in ${\bf B}_{2r}(X_V)$. Then the
functional
\begin{equation}\label{e:2.22}
\tilde{\cal L}^\circ: \bar{\bf B}_\tau(H^0_V)\to\R,\; \xi\mapsto
\tilde{\cal L}(\xi+ \tilde h(\xi))
\end{equation}
is $C^{k-3}$, and as in the proof Lemma~3.1 of \cite{Lu2} we may deduce that  
$$d\tilde{\cal L}^\circ(\xi)[\zeta]=(\tilde A(\xi+
\tilde h(\xi)), \zeta)_{H_V}\;\forall \zeta\in H^0_V\quad\hbox{and}\quad
d^2\tilde{\cal L}^\circ(0)=0.
$$

Let us prove that Theorem~A.1 in \cite{Lu2} can be applied to the
functional ${\cal J}$. By (\ref{e:2.20}) and (\ref{e:2.21})  we only
need to prove

\begin{proposition}\label{prop:2.6}
\begin{enumerate}
\item[{\bf (i)}] For any $(\xi, \zeta)\in \overline{{\bf
B}_\tau(H^0_V)}\times {\bf B}_\tau(H^+_V)$, $\eta_1, \eta_2\in {\bf
B}_\tau(H^-_V)$,
$$
\bigl(D_2{\cal J}(\xi, \zeta+ \eta_2)-D_2{\cal J}(\xi, \zeta+
\eta_1)\bigr)[\eta_2-\eta_1]\le-\frac{a_0}{2}\|\eta_2-\eta_1\|^2_{H_V}.
$$

\item[{\bf (ii)}] $D_2{\cal J}(\xi, \zeta+\eta)[\zeta-\eta]\ge
\frac{a_1}{2}\|\zeta\|^2_{H_V}+ \frac{a_0}{2}\|\eta\|^2_{H_V}$ for
any $(\xi, \zeta, \eta)\in
 \overline{{\bf B}_\tau(H^0_V)}\times {\bf B}_\tau(H^+_V)\times {\bf
 B}_\tau(H^-_V)$; in particular we have
 $D_2{\cal J}(\xi, \zeta)[\zeta]\ge \frac{a_1}{2}\|\zeta\|^2_{H_V}$
for any $(\xi, \zeta)\in \overline{{\bf B}_\tau(H^0_V)}\times {\bf
B}_\tau(H^+_V)$.
\end{enumerate}
\end{proposition}

The last inequality in (ii) shows that the non-decreasing function
$p:(0, \tau]\to (0, \infty)$ in Theorem~A.1(iv) of \cite{Lu2} can be
chosen as $p(t)=\frac{a_1}{2}t^2$.

\noindent{\bf Proof of Proposition~\ref{prop:2.6}}.\quad Recalling that $n>1$,
Lemma~\ref{lem:2.2} implies
$$
U_X^{\rm reg}=U_X\cap C^{k+1}_{V}(I, {\bf
B}_{2\rho}(\R^n))_{\rm reg}=U_X\cap C^{k+1}_{V}(I, {\bf
B}_{2\rho}(\R^n))^\circ_{\rm reg}.
$$

\noindent{\bf Step 1}. \textsf{Proving {\rm (i)}}.\quad Assume $\eta_1\ne\eta_2$. By (\ref{e:2.20}),
\begin{eqnarray}\label{e:2.23}
&&\bigl(D_2{\cal J}(\xi, \zeta+ \eta_2)-D_2{\cal J}(\xi, \zeta+
\eta_1)\bigr)[\eta_2-\eta_1]\\
&=&(\nabla\tilde{\cal L}(\xi+ \tilde h(\xi)+ \zeta+ \eta_2),
\eta_2-\eta_1)_{H_V} - (\nabla\tilde{\cal L}(\xi+ \tilde h(\xi)+ \zeta+ \eta_1),
\eta_2-\eta_1)_{H_V}.\nonumber
\end{eqnarray}

Suppose that  $(\xi, \zeta, \eta_i)\in W_X^{\rm reg}$,
$i=1,2$, i.e.,
$$
\xi_i:=\xi+ \tilde h(\xi)+ \zeta+ \eta_i\in U_X^{\rm
reg}=U_X\cap C^{k+1}_{V}(I, {\bf
B}_{2\rho}(\R^n))^\circ_{\rm reg},\quad i=1,2.
$$
 Since $U_X\cap C^{k+1}_{V}(I, {\bf
B}_{2\rho}(\R^n))$ is convex the smooth curve
$$
[0,1]\ni s\to \mathbbm{p}(s):=\xi+ \tilde h(\xi)+ \zeta+
(1-s)\eta_1+ s\eta_2\in C^{k+1}_{V}(I, {\bf
B}_{2\rho}(\R^n))
$$
also sits in $U_X$.  By Lemma~\ref{lem:2.2} we have a $C^{k+1}$ path
$\mathbbm{q}:[0, 1]\to C^{k+1}_{V}(I, {\bf B}_{2\rho}(\R^n))$ connecting
$\mathbbm{p}(0)$ to $\mathbbm{p}(1)$, which also takes values in $U_X$, such
that
\begin{equation}\label{e:2.24}
\sum^{k+1}_{j=0}\sup_{s\in [0,1]}\left\|\frac{d^j}{ds^j}\mathbbm{p}(s)-\frac{d^j}{ds^j}\mathbbm{q}(s)
\right\|_{C^{k+1}_{V}(I, {\bf B}_{2\rho}(\R^n))}<\|\eta_2-\eta_1\|_{H_V},
\end{equation}
 and that the set
$\{(s,t)\in [0, 1]\times I\,|\, {\bf E}(t,\mathbbm{q}(s))\in{\cal
Z}\}$ is at most finite . The final claim implies that there exists a partition of $[0, 1]$,
$0=s_0<s_1<\cdots<s_m=1$,
 such that $\mathbbm{q}(s)\in U_X^{\rm
reg}$ for each $s\in [0,
1]\setminus\{s_0,s_1,\cdots,s_m\}$. It follows that the map $\tilde A$ is Frech\'et differentiable
at every $\mathbbm{q}(s)$ with $s\in [0,
1]\setminus\{s_0,s_1,\cdots,s_m\}$, and that the map $[0,
1]\setminus\{s_0,s_1,\cdots,s_m\}\ni
s\to d\tilde A(\mathbbm{q}(s))\in\mathscr{L}(X_V)$ is continuous by Lemma~\ref{lem:2.3}.
This implies that the continuous function
$$
\Gamma:[0, 1]\to\R,\; s\mapsto \bigl(\tilde
A(\mathbbm{q}(s)),\eta_2-\eta_1\bigr)_{H_V}.
$$
 is  differentiable in $[0, 1]\setminus\{s_0,s_1,\cdots,s_m\}$.  By the
Newton-Leibniz formula we have
\begin{eqnarray*}
&&\!\!\!\!\!\!\bigl(\tilde A(\xi+ \tilde h(\xi)+ \zeta+ \eta_2),
\eta_2-\eta_1\bigr)_{H_V} - \bigl(\tilde A(\xi+ \tilde h(\xi)+ \zeta+ \eta_1),
\eta_2-\eta_1\bigr)_{H_V}\\
&=&\!\!\!\!\!\!\Gamma(1)-\Gamma(0)=\sum^{m-1}_{j=0}(\Gamma(s_{j+1})-\Gamma(s_j))=
\sum^{m-1}_{j=0}\int^{s_{j+1}}_{s_j}\Gamma'(s)ds\\
&=&\!\!\!\!\!\!\sum^{m-1}_{j=0}\int^{s_{j+1}}_{s_j}\Bigl(d\tilde
A(\mathbbm{q}(s))\bigl[\mathbbm{q}'(s)\bigr],
\eta_2-\eta_1\Bigr)_{H_V}ds\\
&\stackrel{(\ref{e:2.16})}{=}&\!\!\!\!\!\!
\sum^{m-1}_{j=0}\int^{s_{j+1}}_{s_j}\Bigl(\tilde
B(\mathbbm{q}(s))\mathbbm{q}'(s),
\eta_2-\eta_1\Bigr)_{H_V}ds\\
&=&\!\!\!\!\!\! \sum^{m-1}_{j=0}\int^{s_{j+1}}_{s_j}\Bigl(\tilde
B(\mathbbm{q}(s))(\eta_2-\eta_1), \eta_2-\eta_1\Bigr)_{H_V}ds+\\
&+&\sum^{m-1}_{j=0}\int^{s_{j+1}}_{s_j} \Bigl(\tilde
B(\mathbbm{q}(s))\bigl(\mathbbm{q}'(s)-(
\eta_2-\eta_1)\bigr), \eta_2-\eta_1\Bigr)_{H_V}ds\\
&\le&\!\!\!\!\!\!
-a_0\|\eta_2-\eta_1\|^2_{H_V}+\sum^{m-1}_{j=0}\int^{s_{j+1}}_{s_j}\Bigl(\tilde
B(\mathbbm{q}(s))\bigl(\mathbbm{q}'(s)-(
\eta_2-\eta_1)\bigr), \eta_2-\eta_1\Bigr)_{H_V}ds.
\end{eqnarray*}
Moreover, since $\mathbbm{p}'(s)=\eta_2-\eta_1$,
Lemma~\ref{lem:2.5}(ii) and (\ref{e:2.24}) leads to
\begin{eqnarray*}
&&\left|\Bigl(\tilde
B(\mathbbm{q}(s))\bigl(\mathbbm{q}'(s)-(
\eta_2-\eta_1)\bigr), \eta_2-\eta_1\Bigr)_{H_V}\right|\\
&\le
&\omega(\mathbbm{q}(s)))\left\|\frac{d}{ds}(\mathbbm{q}(s)-\mathbbm{p}(s))\right\|_{H_V}
\cdot\|\eta_2-\eta_1\|_{H_V}\\
&\le&\frac{a_0}{2}\left\|\frac{d}{ds}(\mathbbm{q}(s)-\mathbbm{p}(s))\right\|_{X_V}
\cdot\|\eta_2-\eta_1\|_{H_V}\\
&\le&\frac{a_0}{2}\left\|\frac{d}{ds}(\mathbbm{q}(s)-\mathbbm{p}(s))\right\|_{C^3_{V}(I, {\bf B}_{2\rho}(\R^n))}\cdot\|\eta_2-\eta_1\|_{H_V}
\le\frac{a_0}{2}\|\eta_2-\eta_1\|^2_{H_V}
\end{eqnarray*}
for each $s\in [0, 1]\setminus\{s_0,s_1,\cdots,s_m\}$, and hence
$$
\left|\sum^{m-1}_{j=0}\int^{s_{j+1}}_{s_j}\Bigl(\tilde
B(\mathbbm{q}(s))\bigl(\mathbbm{q}'(s)-(
\eta_2-\eta_1)\bigr), \eta_2-\eta_1\Bigl)_{H_V}ds\right|\le
\frac{a_0}{2}\|\eta_2-\eta_1\|^2_{H_V}.
$$
 From this and (\ref{e:2.23}) It follows that
\begin{eqnarray*}
&&\left|\bigl(D_2{\cal J}(\xi, \zeta+ \eta_2)-D_2{\cal J}(\xi, \zeta+
\eta_1)\bigr)[\eta_2-\eta_1]\right|\\
&=&\left|\bigl(\tilde A(\xi+ \tilde h(\xi)+ \zeta+ \eta_2),
\eta_2-\eta_1\bigr)_{H_V} - \bigl(\tilde A(\xi+ \tilde h(\xi)+ \zeta+ \eta_1),
\eta_2-\eta_1\bigr)_{H_V}\right|\\
&\le&-\frac{a_0}{2}\|\eta_2-\eta_1\|^2_{H_V}\quad\forall (\xi,
\zeta, \eta_i)\in W_X^{\rm reg}, i=1,2.
\end{eqnarray*}
By the density of $W_X^{\rm reg}$ in $\overline{{\bf
B}_\tau(H^0_V)}\oplus {\bf B}_\tau(H^-_V)\oplus{\bf B}_\tau(H^+_V)$
this  yields (i).

\noindent{\bf Step 2}. \textsf{Proving {\rm (ii)}}.\quad
 By (\ref{e:2.20}) and (\ref{e:2.21}) we have
\begin{eqnarray}\label{e:2.25}
&&D_2{\cal J}(\xi, \zeta+ \eta)[\zeta-\eta]\nonumber\\
&=&\bigl(\nabla\tilde{\cal L}(\xi+ \tilde h(\xi)+ \zeta+ \eta), \zeta-\eta\bigr)_{H_V} -
\bigl(\nabla\tilde{\cal L}(\xi+ \tilde h(\xi)), \zeta-\eta\bigr)_{H_V}.
\end{eqnarray}
For a given $(\xi, \zeta, \eta)\in W_X^{\rm reg}$ with
$(\zeta,\eta)\ne (0, 0)$, i.e.,
$$\xi+ \tilde h(\xi)+ \zeta+ \eta\in  U_X^{\rm
reg}\quad\hbox{and}\quad\|\zeta\|_{H_V}+\|\eta\|_{H_V}\ne 0,
$$
since $\xi+ \tilde h(\xi)\in {\bf B}_{2r}(X_V)\subset
U_X$ and $U_X^{\rm reg}$ is dense in $U_X$ there
exists a sequence $(\zeta_l)\subset X_V$ with $\|\zeta_l\|_{X_V}\to
0$ as $l\to\infty$, such that $\xi+ \tilde h(\xi)+\zeta_l\in
U_X^{\rm reg}$ for all $l\in\N$.
By the inequalities that $\|\zeta_l\|_{H_V}\le\|\zeta_l\|_{X_V}\;\forall l$ we can
find an integer $l_0>0$ such that
\begin{equation}\label{e:2.26}
\|\zeta_l\|_{H_V}<
(\|\zeta\|_{H_V}+\|\eta\|_{H_V})/4\quad\forall l>l_0.
\end{equation}
\textbf{Fix a $l\in\N$ with $l>l_0$.}  Consider the smooth curves
$$
[0,1]\ni s\to \mathbbm{p}_l(s):=\xi+ \tilde h(\xi)+ (1-s)\zeta_l+
s(\zeta+ \eta)\in U_X.
$$
By Lemma~\ref{lem:2.2} we can choose a $C^{k+1}$ path
$\mathbbm{q}_l:[0, 1]\to C^{k+1}_{V}(I, {\bf B}_{2\rho}(\R^n))$
connecting $\mathbbm{p}_l(0)$ to $\mathbbm{p}_l(1)$, which also takes values
in $U_X$, such that
\begin{equation}\label{e:2.27}
\sum^{k+1}_{j=0}\sup_{s\in [0,1]}\left\|\frac{d^j}{ds^j}\mathbbm{p}_l(s)-\frac{d^j}{ds^j}\mathbbm{q}_l(s)
\right\|_{C^{k+1}_{V}(I, {\bf B}_{2\rho}(\R^n))}<\frac{1}{4}(\|\zeta\|_{H_V}+\|\eta\|_{H_V}),
\end{equation}
and that each set $\{(s,t)\in [0, 1]\times I\,|\, {\bf
E}(t,\mathbbm{q}_l(s))\in{\cal Z}\}$ is at most finite.
So we have a partition of $[0, 1]$,
$0=s_0<s_1<\cdots<s_m=1$,
 such that $\mathbbm{q}_l(s)\in U_X^{\rm
reg}$ for each $s\in [0,
1]\setminus\{s_0,s_1,\cdots,s_m\}$. It follows that
 the continuous function
$$
\Gamma_l:[0,1]\to\R,\; s\mapsto\bigl(\tilde
A(\mathbbm{q}_l(s)),\zeta-\eta\bigr)_{H_V}
$$
is differentiable at each point $s\in [0,
1]\setminus\{s_0,s_1,\cdots,s_m\}$. As in Step 1 we have
\begin{eqnarray*}
&&\bigl(\tilde A(\xi+ \tilde h(\xi)+ \zeta+ \eta), \zeta-\eta\bigr)_{H_V} -
\bigl(\tilde A(\xi+ \tilde h(\xi)+\zeta_l), \zeta-\eta\bigr)_{H_V}\\
&=&\!\!\!\!\!\!\Gamma_l(1)-\Gamma_l(0)=\sum^{m-1}_{j=0}\int^{s_{j+1}}_{s_j}\Gamma'_l(s)ds\\
&=&\!\!\!\!\!\!\sum^{m-1}_{j=0}\int^{s_{j+1}}_{s_j}\Bigl(d\tilde
A(\mathbbm{q}_l(s))[\mathbbm{q}'_l(s)],
\zeta-\eta\Bigr)_{H_V}ds\\
&\stackrel{(\ref{e:2.16})}{=}&\!\!\!\!\!\!
\sum^{m-1}_{j=0}\int^{s_{j+1}}_{s_j}\bigl(\tilde
B(\mathbbm{q}_l(s))\mathbbm{q}'_l(s),
\zeta-\eta\bigr)_{H_V}ds\\
&=&\!\!\!\!\!\! \sum^{m-1}_{j=0}\int^{s_{j+1}}_{s_j}\Bigl(\tilde
B(\mathbbm{q}_l(s))(\zeta+\eta), \zeta-\eta\Bigr)_{H_V}ds+\\
&&\quad +\sum^{m-1}_{j=0}\int^{s_{j+1}}_{s_j} \Bigl(\tilde
B(\mathbbm{q}_l(s))\bigl(\mathbbm{q}'_l(s)-(
\zeta+\eta)\bigr), \zeta-\eta\Bigr)_{H_V}ds\\
&=&\!\!\!\!\!\! \sum^{m-1}_{j=0}\int^{s_{j+1}}_{s_j}\bigl(\tilde
B(\mathbbm{q}_l(s))\zeta, \zeta\bigr)_{H_V}ds-
\sum^{m-1}_{j=0}\int^{s_{j+1}}_{s_j}\bigl(\tilde
B(\mathbbm{q}_l(s))\eta, \eta\bigr)_{H_V}ds\\
&&\quad +\sum^{m-1}_{j=0}\int^{s_{j+1}}_{s_j} \Bigl(\tilde
B(\mathbbm{q}_l(s))\bigl(\mathbbm{q}'_l(s)-(
\zeta+\eta)\bigr), \zeta-\eta\Bigr)_{H_V}ds\\
&\ge&\!\!\!\!\!\! a_1\|\zeta\|^2_{H_V}+
a_0\|\eta\|^2_{H_V}+\sum^{m-1}_{j=0}\int^{s_{j+1}}_{s_j}\Bigl(\tilde
B(\mathbbm{q}_l(s))\bigl(\mathbbm{q}'_l(s)-(
\zeta+\eta)\bigr), \zeta-\eta\Bigr)_{H_V}ds
\end{eqnarray*}
by Lemma~\ref{lem:2.5}. As in Step 1 we can derive from Lemma~\ref{lem:2.5},
 (\ref{e:2.26}) and (\ref{e:2.27}) that
\begin{eqnarray*}
&&\left|\sum^{m-1}_{j=0}\int^{s_{j+1}}_{s_j}\Bigl(\tilde
B(\mathbbm{q}_l(s))\bigl(\mathbbm{q}'_l(s)-(
\zeta+\eta)\bigr),
\zeta-\eta\Bigr)_{H_V}ds\right|\\
&\le&\sum^{m-1}_{j=0}\int^{s_{j+1}}_{s_j}\omega(\mathbbm{q}_l(s)))\left\|\mathbbm{q}'_l(s)-
\mathbbm{p}'_l(s)+ \zeta_l\right\|_{H_V}\cdot\|\zeta-\eta\|_{H_V}ds\\
 &\le& \frac{\min\{a_0,
 a_1\}}{2}\left(\frac{\|\zeta\|_{H_V}+\|\eta\|_{H_V}}{4}+ \|\zeta_l\|_{H_V}\right)\|\zeta-\eta\|_{H_V}\\
 &\le&
\frac{\min\{a_0,
 a_1\}}{2}(\|\zeta\|^2_{H_V}+\|\eta\|^2_{H_V})
\end{eqnarray*}
and so
\begin{eqnarray*}
&&\bigl(\tilde A(\xi+ \tilde h(\xi)+ \zeta+ \eta), \zeta-\eta\bigr)_{H_V} -
\bigl(\tilde A(\xi+ \tilde h(\xi)+\zeta_l), \zeta-\eta\bigr)_{H_V}\\
&&\ge \frac{\min\{a_0,
 a_1\}}{2}(\|\zeta\|^2_{H_V}+\|\eta\|^2_{H_V}).
\end{eqnarray*}
Recall that  $l>l_0$ is arbitrary. Let $l\to\infty$ we obtain
\begin{eqnarray*}
&&\bigl(\tilde A(\xi+ \tilde h(\xi)+ \zeta+ \eta), \zeta-\eta\bigr)_{H_V} -
\bigl(\tilde A(\xi+ \tilde h(\xi)), \zeta-\eta\bigr)_{H_V}\\
&&\ge \frac{\min\{a_0,
 a_1\}}{2}(\|\zeta\|^2_{H_V}+\|\eta\|^2_{H_V})
\end{eqnarray*}
and hence
$$
D_2{\cal J}(\xi, \zeta+ \eta)[\zeta-\eta]
\ge \frac{\min\{a_0,
 a_1\}}{2}(\|\zeta\|^2_{H_V}+\|\eta\|^2_{H_V})
$$
by (\ref{e:2.25}).
This can also hold for any $(\xi, \zeta,
\eta)$ in $\overline{{\bf B}_\tau(H^0_V)}\times {\bf B}_\tau(H^+_V)\times {\bf
 B}_\tau(H^-_V)$ because
 $\{(\xi, \zeta, \eta)\in W_X^{\rm reg}\,|\,
(\zeta,\eta)\ne (0, 0)\}$ is dense in
$ \overline{{\bf B}_\tau(H^0_V)}\times {\bf B}_\tau(H^+_V)\times {\bf
 B}_\tau(H^-_V)$.   (ii) is proved.
 \hfill$\Box$\vspace{2mm}

Then repeating the arguments in Step 4 of proof of
\cite[Theorem~1.1]{Lu1} (or the proof of Lemma~3.6 in  \cite{Lu2})
 we obtain the following
splitting theorem.

\begin{theorem}\label{th:2.7}
Under the notation above, there exists   an origin-preserving local
homeomorphism $\tilde\psi$ from ${\bf B}_\tau(H_V)$ to an open
neighborhood of $0\in H_V$  such that
$$
\tilde{\cal L}\circ\tilde\psi(\xi)=\|P_V^+\xi\|^2_{H_V}-
\|P^-_V\xi\|^2_{H_V} + \tilde{\cal
L}^{\circ}(P^0_V\xi)\quad\forall\xi\in{\bf B}_\tau(H_V).
$$
\end{theorem}

\subsubsection{Completing the proof  of
Theorem~\ref{th:1.1}}\label{sec:2.1.2}

The completion of the proof of Theorem~\ref{th:1.1} is done by repeating of Step~2 in \cite[\S4]{Lu3}.
Let us outline it.
The differential at $0$ of the chart $\Phi$ in
(\ref{e:2.2}), $d\Phi(0):H_V\to T_{\gamma_0}\Lambda_N(M)$, given by $d\Phi(0)[\xi]=\sum^n_{i=1}\xi_ie_i$,
 is a Hilbert space isomorphism and  ${\rm EXP}_{\gamma_0}\circ
d\Phi(0)=\Phi$ on ${\bf B}_{2r}(H_V):=\{\xi\in
H_V\,|\,\|\xi\|_1<2r\}$.
  Since $\tilde{\cal L}={\cal
L}\circ\Phi$ on ${\bf B}_{2r}(H_V)$ by (\ref{e:2.10}), we get
\begin{equation}\label{e:2.28}
{\cal L}\circ{\rm EXP}_{\gamma_0}\circ d\Phi(0)={\cal
L}\circ\Phi=\tilde{\cal L}\quad\hbox{on}\; {\bf B}_{2r}(H_V)
\end{equation}
and so $\nabla({\cal L}\circ{\rm
EXP}_{\gamma_0})(d\Phi(0)[\xi])=d\Phi(0)[\nabla\tilde{\cal L}(\xi)]$ for
all $\xi\in {\bf B}_{2r}(H_{V})$. It follows that
\begin{eqnarray}
&&{\cal A}(d\Phi(0)\xi)=d\Phi(0)\tilde A(\xi)\quad\forall\xi\in {\bf
B}_{2r}(H_{V})\cap X_{V},\label{e:2.29}\\
&&d{\cal A}(0)\circ d\Phi(0)=d\Phi(0)\circ d\tilde
A(0)\label{e:2.30}
\end{eqnarray}
because the Hilbert space isomorphism $d\Phi(0):H_{V}\to
T_{\gamma_0}\Lambda_N(M)$  induces a Banach space isomorphism from
$X_{V}$ to $T_{\gamma_0}{\cal X}=T_{\gamma_0}C^1_N(I,M)$.
(\ref{e:2.30}) implies that $d\Phi(0)(H_V^\star)= {\bf
H}^\star(d^2{\cal L}|_{{\cal X}_N}(\gamma_0))$ for $\star=-,0,+$, and therefore
\begin{eqnarray}\label{e:2.31}
d\Phi(0)\circ P_V^\star=P^\star\circ d\Phi(0),\quad \star=-,0,+.
\end{eqnarray}
Shrinking $\delta>0$ above (\ref{e:1.12}) so that $\delta<\tau$,
as below  \cite[(4.30)]{Lu3} we obtain
\begin{eqnarray*}
0=(I-P^0)\circ {\cal A}\bigl(\zeta+ d\Phi(0)\circ\tilde
h(d\Phi(0)^{-1}\zeta)\bigr)\quad\forall \zeta\in {\bf B}_\delta\bigl({\bf H}^0(d^2{\cal L}|_{{\cal
X}_N}(\gamma_0))\bigr)
\end{eqnarray*}
by (\ref{e:2.31}), (\ref{e:2.30}) and (\ref{e:2.18}).
This and (\ref{e:1.11}) yield  $h(\zeta)=d\Phi(0)\circ\tilde
h(d\Phi(0)^{-1}\zeta)$ by the uniqueness of $h$.
 From it, (\ref{e:1.12}), (\ref{e:2.28}) and the definition of $\tilde{\cal
L}^\circ$ in (\ref{e:2.22}) we derive
\begin{eqnarray}\label{e:2.32}
{\cal L}^\circ(d\Phi(0)\xi)
=\tilde{\cal L}\bigl(\xi+ \tilde h(\xi)\bigr)=\tilde{\cal
L}^\circ(\xi),\quad\forall\xi\in {\bf B}_\delta(H_V^0)
\end{eqnarray}
as in \cite[(4.31)]{Lu3}.
 Since  $d\Phi(0)\left({\bf
B}_\delta(H_V)\right)={\bf B}_\delta(T_{\gamma_0}\Lambda_N(M))$,
(\ref{e:2.31}) implies that
\begin{eqnarray}\label{e:2.33}
(P_V^\star\xi, P_V^\star\xi)_{H_V}=\langle
d\Phi(0)\circ P_V^\star\xi, d\Phi(0)\circ
P_V^\star\xi\rangle_{1}
=\langle P^\star\circ d\Phi(0)\xi, P^\star\circ
d\Phi(0)\xi\rangle_{1}
\end{eqnarray}
for $\xi\in {\bf B}_\delta(H_V)\subset {\bf B}_\tau(H_V)$ and
$\star=+,-$.
Define $\psi: {\bf B}_\delta(T_{\gamma_0}\Lambda_N(M))\to
T_{\gamma_0}\Lambda_N(M)$ by
$\psi=d\Phi(0)\circ\tilde\psi\circ[d\Phi(0)]^{-1}$. For $\zeta\in
{\bf B}_\delta(T_{\gamma_0}\Lambda_N(M))$ and
$\xi=[d\Phi(0)]^{-1}\zeta$,  we obtain
\begin{eqnarray*}
&&{\cal L}\circ{\rm EXP}_{\gamma_0}\circ\psi(\zeta)={\cal
L}\circ{\rm EXP}_{\gamma_0}\circ
d\Phi(0)\circ\tilde\psi\circ[d\Phi(0)]^{-1}\zeta=\tilde{\cal
L}\circ\tilde\psi(\xi),\\
&&{\cal L}^\circ(P^0\zeta)={\cal L}^\circ(P^0\circ
d\Phi(0)\xi)={\cal L}^\circ(d\Phi(0)\circ P^0_V\xi)=\tilde{\cal
L}^\circ(P^0_V\xi)
\end{eqnarray*}
by  (\ref{e:2.28}) and (\ref{e:2.32}). These, (\ref{e:2.33}) and Theorem~\ref{th:2.7}
give Theorem~\ref{th:1.1}. \hfill$\Box$\vspace{2mm}

\subsection{Proof  of
Theorem~\ref{th:1.4}}\label{sec:2.2}

Define the functional
$$
\mathscr{J}:\overline{{\bf B}_\tau\bigl({\bf H}^0(d^2{\cal L}|_{{\cal
X}_N}(\gamma_0))\bigr)}\times\Bigr({\bf B}_\tau\bigl({\bf H}^-(d^2{\cal L}|_{{\cal
X}_N}(\gamma_0))\bigr)\oplus{\bf B}_\tau\bigl({\bf H}^+(d^2{\cal L}|_{{\cal
X}_N}(\gamma_0))\bigr)\Bigl)\to\R
$$
by $\mathscr{J}(x,y)={\cal L}\circ{\rm EXP}_{\gamma_0}(x+ h(x)+ y)-
{\cal L}\circ{\rm EXP}_{\gamma_0}(x+ h(x))$, where $h$ is as in (\ref{e:2.32}),
i.e., $h(x)=d\Phi(0)\circ\tilde
h(d\Phi(0)^{-1}x)$. Note that the Hilbert space isomorphism $d\Phi(0)$
maps ${\bf B}_\tau(H^\ast_V)$ onto ${\bf B}_\tau\bigl({\bf H}^\ast(d^2{\cal L}|_{{\cal
X}_N}(\gamma_0))\bigr)$, $\ast=0,+,-$. We have
$$
{\cal J}(\xi, \zeta)=\mathscr{J}(d\Phi(0)\xi, d\Phi(0)\zeta)\quad
\forall\xi\in \overline{{\bf B}_\tau(H^0_V)}\quad\&\quad\forall\zeta\in {\bf
B}_\tau(H^-_V)\oplus {\bf B}_\tau(H^+_V),
$$
 where ${\cal J}$ is given by (\ref{e:2.19}). Moreover, (\ref{e:2.20}) and (\ref{e:2.21}) imply
\begin{eqnarray}
&&D_2\mathscr{J}(\xi,\zeta)[\eta] =((I-P^0)\nabla({\cal L}\circ{\rm EXP}_{\gamma_0})(\xi+ h(\xi)+ \zeta), \eta)_{H_V},\label{e:2.39}\\
&& \mathscr{J}(\xi, 0)=0\quad\hbox{and}\quad D_2\mathscr{J}(\xi,
0)[\eta]=0\label{e:2.40}
\end{eqnarray}
for any $\xi\in \overline{{\bf B}_\tau\bigl({\bf H}^0(d^2{\cal L}|_{{\cal
X}_N}(\gamma_0))\bigr)}$, $\zeta\in{\bf B}_\tau\bigl({\bf H}^-(d^2{\cal L}|_{{\cal
X}_N}(\gamma_0))\bigr)\oplus{\bf B}_\tau\bigl({\bf H}^+(d^2{\cal L}|_{{\cal
X}_N}(\gamma_0))\bigr)$ and $\eta\in {\bf H}^-(d^2{\cal L}|_{{\cal
X}_N}(\gamma_0))\oplus {\bf H}^+(d^2{\cal L}|_{{\cal
X}_N}(\gamma_0))$.
From Proposition~\ref{prop:2.6} we deduce

\begin{proposition}\label{prop:2.9}
\begin{enumerate}
\item[{\bf (i)}] For any $(\xi, \zeta)\in \overline{{\bf B}_\tau\bigl({\bf H}^0(d^2{\cal L}|_{{\cal
X}_N}(\gamma_0))\bigr)}\times{\bf B}_\tau\bigl({\bf H}^+(d^2{\cal L}|_{{\cal
X}_N}(\gamma_0))\bigr)$ and  $\eta_1, \eta_2\in {\bf B}_\tau\bigl({\bf H}^+(d^2{\cal L}|_{{\cal
X}_N}(\gamma_0))\bigr)$,
$$
\bigl(D_2\mathscr{J}(\xi, \zeta+ \eta_2)-D_2\mathscr{J}(\xi, \zeta+
\eta_1)\bigr)[\eta_2-\eta_1]\le-\frac{a_0}{2}\|\eta_2-\eta_1\|^2_{1}.
$$

\item[{\bf (ii)}]  For any $\xi\in\overline{{\bf B}_\tau\bigl({\bf H}^0(d^2{\cal L}|_{{\cal
X}_N}(\gamma_0))\bigr)}$, $\zeta\in{\bf B}_\tau\bigl({\bf H}^+(d^2{\cal L}|_{{\cal
X}_N}(\gamma_0))\bigr)$ and $\eta\in{\bf B}_\tau\bigl({\bf H}^+(d^2{\cal L}|_{{\cal
X}_N}(\gamma_0))\bigr)$ it holds that
$$
D_2\mathscr{J}(\xi, \zeta+\eta)[\zeta-\eta]\ge
\frac{a_1}{2}\|\zeta\|^2_{1}+ \frac{a_0}{2}\|\eta\|^2_{1};
$$
 in particular we have
 $$
 D_2\mathscr{J}(\xi, \zeta)[\zeta]\ge \frac{a_1}{2}\|\zeta\|^2_{1}
 $$
for any $(\xi, \zeta)\in \overline{{\bf B}_\tau\bigl({\bf H}^0(d^2{\cal L}|_{{\cal
X}_N}(\gamma_0))\bigr)}\times{\bf B}_\tau\bigl({\bf H}^+(d^2{\cal L}|_{{\cal
X}_N}(\gamma_0))\bigr)$.
\end{enumerate}
\end{proposition}

(Theorem~\ref{th:1.1} may also follow from this as in the proof of Theorem~\ref{th:2.7}.)

Let $\tilde{\cal A}$ be the restriction of
$\nabla({\cal L}\circ\widetilde{\rm EXP}_{\gamma_0})$ to ${\bf
B}_{2\rho}(T_{\gamma_0}\tilde\Lambda)\cap T_{\gamma_0}\tilde{\cal X}_N$. Then
\begin{equation}\label{e:2.41}
\tilde{\cal A}(x)={\cal A}(x)\quad\forall x\in {\bf
B}_{2\rho}(T_{\gamma_0}\tilde\Lambda)\cap T_{\gamma_0}\tilde{\cal X}_N
\end{equation}
by the assumption (iii), where ${\cal A}$ is the map  below (\ref{e:1.7}).
So $\tilde{\cal A}:{\bf B}_{2\rho}(T_{\gamma_0}\tilde\Lambda)\cap T_{\gamma_0}\tilde{\cal X}_N
\to T_{\gamma_0}\tilde{\cal X}_N$  is  $C^{k-3}$ near $0\in T_{\gamma_0}\tilde{\cal X}_N$.
(\ref{e:2.41}) implies
$\langle\tilde{\cal A}(tx), y\rangle_1=\langle{\cal A}(tx), y\rangle_1$
for any $y\in T_{\gamma_0}\tilde\Lambda$ and $t\in\R$ with small $|t|$.
On both sides taking derivatives at $t=0$ yields
\begin{equation}\label{e:2.42}
d^2({\cal L}|_{\tilde{\cal X}_N})(\gamma_0)[x,y]=\langle d\tilde{\cal A}(0)x, y\rangle_1=\langle d{\cal A}(0)x, y\rangle_1=d^2({\cal L}|_{{\cal X}_N})(\gamma_0)[x,y].
\end{equation}
As for $d^2({\cal L}|_{{\cal X}_N})(\gamma_0)$ we also use $d^2({\cal L}|_{\tilde{\cal X}_N})(\gamma_0)$
to denote the continuous symmetric bilinear extension of it on $T_{\gamma_0}\tilde\Lambda$.
Then (\ref{e:2.42}) also holds for all $x\in T_{\gamma_0}\tilde\Lambda$
by the density of $T_{\gamma_0}\tilde{\cal X}_N$ in
$T_{\gamma_0}\tilde\Lambda$. Note that the conditions (i) and (iv) imply
 ${\bf H}^0\bigl(d^2({\cal L}|_{{\cal X}_N})(\gamma_0)\bigr)\subset T_{\gamma_0}\tilde{\cal X}_N$.
 It follows that
\begin{equation}\label{e:2.43}
{\bf H}^\ast\bigl(d^2({\cal L}|_{\tilde{\cal X}_N})(\gamma_0)\bigr)
\subset{\bf H}^\ast\bigl(d^2({\cal L}|_{{\cal X}_N})(\gamma_0)\bigr),\quad \ast=0,+,-.
\end{equation}
(Actually, for $\ast=+,-$ these are obvious by (\ref{e:2.42}).)
 Shrinking $\delta>0$ in (\ref{e:1.10}) (if necessary) and using
the implicit function theorem yield a unique $C^{k-3}$ map
$$
 \tilde{h}: {\bf B}_\delta\bigl({\bf H}^0(d^2{\cal L}|_{\tilde{\cal
X}_N}(\gamma_0))\bigr)\to {\bf H}^-(d^2{\cal L}|_{\tilde{\cal
X}_N}(\gamma_0))\dot{+}\bigl({\bf H}^+(d^2{\cal L}|_{\tilde{\cal
X}_N}(\gamma_0))\cap T_{\gamma_0}\tilde{\cal X}_N\bigr)
$$
such that $\tilde{h}(0)=0$, $d\tilde{h}(0)=0$ and
\begin{equation}\label{e:2.44}
(\tilde{P}^+ +\tilde{P}^-)\tilde{\cal A}(\xi+
\tilde{h}(\xi))=0\;\forall\xi\in {\bf B}_\delta\bigl({\bf H}^0(d^2{\cal
L}|_{\tilde{\cal X}_N}(\gamma_0))\bigr),
\end{equation}
where $\tilde{P}^\star$ ($\star=-,0,+$) are the orthogonal
projections,
$$
T_{\gamma_0}\tilde\Lambda={\bf H}^0(d^2{\cal
L}|_{\tilde{\cal X}_N}(\gamma_0))\oplus{\bf H}^+(d^2{\cal
L}|_{\tilde{\cal X}_N}(\gamma_0))\oplus{\bf H}^-(d^2{\cal
L}|_{\tilde{\cal X}_N}(\gamma_0))\to {\bf H}^\star(d^2{\cal
L}|_{\tilde{\cal X}_N}(\gamma_0)).
$$
By (\ref{e:2.41}) and (\ref{e:2.43}) we deduce that (\ref{e:2.44}) is equal to the equality
$$
({P}^+ +{P}^-){\cal A}(\xi+
\tilde{h}(\xi))=0\;\forall\xi\in {\bf B}_\delta\bigl({\bf H}^0(d^2{\cal
L}|_{\tilde{\cal X}_N}(\gamma_0))\bigr).
$$
By uniqueness, this and (\ref{e:1.11}) lead to
$$
\tilde{h}(\xi)=h(\xi)\quad\forall \xi\in {\bf B}_\delta\bigl({\bf H}^0(d^2{\cal
L}|_{\tilde{\cal X}_N}(\gamma_0))\bigr),
$$
which implies that for  ${\cal L}^\circ$ in (\ref{e:1.12}) and  $\tilde{\cal L}^\circ$ in
(\ref{e:2.22}) we have
\begin{equation}\label{e:2.45}
{\cal L}^\circ(\xi)=
\tilde{\cal L}^\circ(\xi)= {\cal L}\bigl({\rm EXP}_{\gamma_0}(\xi+  \tilde{h}(\xi))\bigr)\quad\forall \xi\in {\bf B}_\delta\bigl({\bf H}^0(d^2{\cal
L}|_{\tilde{\cal X}_N}(\gamma_0))\bigr).
\end{equation}
By (\ref{e:2.22}), $\tilde{\cal L}^\circ$ is $C^{k-3}$, has an isolated
critical point $0$, and $d^2\tilde{\cal L}^\circ(0)=0$.
Define
$$
\tilde{\mathscr{J}}:\overline{{\bf B}_\tau\bigl({\bf H}^0(d^2{\cal L}|_{\tilde{\cal
X}_N}(\gamma_0))\bigr)}\times\Bigr({\bf B}_\tau\bigl({\bf H}^-(d^2{\cal L}|_{\tilde{\cal
X}_N}(\gamma_0))\bigr)\oplus{\bf B}_\tau\bigl({\bf H}^+(d^2{\cal L}|_{\tilde{\cal
X}_N}(\gamma_0))\bigr)\Bigl)\to\R
$$
by $\tilde{\mathscr{J}}(x,y)={\cal L}\circ\widetilde{{\rm EXP}}_{\gamma_0}(x+ \tilde{h}(x)+ y)-
{\cal L}\circ\widetilde{{\rm EXP}}_{\gamma_0}(x+ \tilde{h}(x))$. It is easily checked that
$\tilde{\mathscr{J}}(x,y)={\mathscr{J}}(x,y)$ and hence that
Proposition~\ref{prop:2.9} also holds if
$\overline{{\bf B}_\tau\bigl({\bf H}^0(d^2{\cal L}|_{{\cal
X}_N}(\gamma_0))\bigr)}$, ${\bf B}_\tau\bigl({\bf H}^\star(d^2{\cal L}|_{{\cal
X}_N}(\gamma_0))\bigr)$ are replaced by
$\overline{{\bf B}_\tau\bigl({\bf H}^0(d^2{\cal L}|_{\tilde{\cal
X}_N}(\gamma_0))\bigr)}$ and ${\bf B}_\tau\bigl({\bf H}^\star(d^2{\cal L}|_{\tilde{\cal
X}_N}(\gamma_0))\bigr)$, $\star=+,-$, respectively.
Checking Step 4 of proof of \cite[Theorem~1.1]{Lu1} (or the proof of Lemma~3.6 in  \cite{Lu2})
the first claim is easily obtained. As a consequence
we have
\begin{equation}\label{e:2.46}
C_q(\tilde{\cal L}, \gamma_0;\K)\cong C_{q-m^-(\tilde{\cal L},\gamma_0)}(\tilde{\cal L}^\circ,
0;\K)\quad\forall q=0,1,\cdots,
\end{equation}
where $m^-(\tilde{\cal L},\gamma_0)=\dim{\bf H}^-(d^2{\cal L}|_{\tilde{\cal
X}_N}(\gamma_0))$.

If $m^0({\cal L}|_{\tilde\Lambda},\gamma_0)=m^0({\cal L},\gamma_0)$, (\ref{e:2.43}) leads to
${\bf H}^0\bigl(d^2({\cal L}|_{\tilde{\cal X}_N})(\gamma_0)\bigr)
={\bf H}^0\bigl(d^2({\cal L}|_{{\cal X}_N})(\gamma_0)\bigr)$ and so
${\cal L}^\circ=\tilde{\cal L}^\circ$ by (\ref{e:2.45}). The latter implies the claim in Theorem~\ref{th:1.4}(a).

Furthermore, if $m^-({\cal L}|_{\tilde\Lambda},\gamma_0)$ is also equal to $m^-({\cal L},\gamma_0)$, Theorem~\ref{th:1.2} and (\ref{e:2.46}) give
the claim in Theorem~\ref{th:1.4}(b).
\hfill$\Box$\vspace{2mm}

\section{Proof of Theorem~\ref{th:1.7}}\label{sec:3}
\setcounter{equation}{0}

The basic idea is similar to that of Theorem~\ref{th:1.1} in Section~\ref{sec:2},
but more complex. We shall complete the proof  in three steps.

\noindent{\bf Step 1}.  Since  $\gamma_0$ is a $C^k$-map to $M$ with
$k\ge 6$, starting with a unit orthogonal
 frame at $T_{\gamma(0)}M$ and using the parallel
transport along $\gamma_0$ with respect to the Levi-Civita
connection of the Riemannian metric $g$ we get a unit orthogonal
parallel $C^k$ frame field $\R\to \gamma_0^\ast TM,\;t\mapsto
(e_1(t),\cdots, e_n(t))$.
 Note that  there exists a unique orthogonal matrix
 $E_{\gamma_0}$ such that
$(e_1(1),\cdots, e_n(1))=(\mathbb{I}_{g\ast}e_1(0),\cdots,
\mathbb{I}_{g\ast}e_n(0))E_{\gamma_0}$. (\textsf{All vectors in
$\R^n$ will be understood as row vectors.}) By the elementary matrix theory there exists an
orthogonal matrix $\Xi$ such that
$$
\Xi^{-1}E_{\gamma_0}\Xi={\rm
diag}(S_1,\cdots,S_{\sigma})\in\R^{n\times n}\quad\hbox{with}\;{\rm
ord}(S_1)\ge\cdots\ge{\rm ord}(S_{\sigma}),
$$
where each $S_j$ is either $1$, or $-1$, or
${\scriptscriptstyle\left(\begin{array}{cc}
\cos\theta_j& \sin\theta_j\\
-\sin\theta_j& \cos\theta_j\end{array}\right)}$, $0<\theta_j<\pi$.
 So replacing $(e_1,\cdots,e_n)$ by $(e_1,\cdots,e_n)\Xi$ we may
assume
\begin{eqnarray}\label{e:3.1}
E_{\gamma_0}={\rm diag}(S_1,\cdots,S_{\sigma})\in\R^{n\times n}.
\end{eqnarray}
Since $\gamma_0(t+1)=\mathbb{I}_g(\gamma_0(t))\;\forall t\in\R$ and
$\R\ni t\mapsto (\mathbb{I}_{g\ast}(e_1(t)),\cdots,
\mathbb{I}_{g\ast}(e_n(t)))E_{\gamma_0}$ is also a unit orthogonal
parallel $C^k$ frame field along $\mathbb{I}_g\circ\gamma_0$, it is
easily proved that
\begin{eqnarray}\label{e:3.2}
(e_1(t+1),\cdots, e_n(t+1))=(\mathbb{I}_{g\ast}(e_1(t)),\cdots,
\mathbb{I}_{g\ast}(e_n(t)))E_{\gamma_0}\quad\forall t\in\R.
\end{eqnarray}
Consider the subspaces of $C^1(\R,\R^n)$ and
$W^{1,2}_{loc}(\R,\R^n)$,
\begin{eqnarray*}
&&X_{{\gamma_0}}=C^1_{{\gamma_0}}(\R,\R^n)=\{\xi\in C^1(\R,
\R^n)\,|\,\xi(t+1)^T=E_{{\gamma_0}}\xi(t)^T\;\forall t\},\\
&&H_{{\gamma_0}}=W^{1,2}_{{\gamma_0}}(\R,\R^n)=\{\xi\in
W^{1,2}_{loc}(\R,
\R^n)\,|\,\xi(t+1)^T=E_{{\gamma_0}}\xi(t)^T\;\forall t\}.
\end{eqnarray*}
 Here $\xi(t)^T$ denotes  the transpose of the
matrix $\xi(t)$ as usual.  The former is a Banach space according to
the usual $C^1$-norm because $E_{{\gamma_0}}\in O(n)$, and the latter is
 a Hilbert space with respect to the inner product
\begin{eqnarray}\label{e:3.3}
\langle\xi,\eta\rangle_{1,2}=\int_0^1[(\xi(t),\eta(t))_{\R^n}+
(\dot\xi(t),\dot\eta(t))_{\R^n}]dt.
\end{eqnarray}
 Let $\|\cdot\|_{1,2}$ be the induced norm. Clearly, $X_{{\gamma_0}}$ is dense in
$H_{{\gamma_0}}$, and there exists  a Hilbert space isomorphism
 $$
\mathfrak{I}_{{\gamma_0}}: ({H}_{{\gamma_0}},
\langle\cdot,\cdot\rangle_{1,2})\to
(T_{\gamma_0}\Lambda(\mathbb{I}_g, M),
\langle\cdot,\cdot\rangle_1),\quad\xi\mapsto \sum^n_{j=1}\xi_j{e}_j.
$$
 Moreover we have a small open ball $B^n_{2\rho}(0)\subset\R^n$ and a $C^k$
map
 \begin{eqnarray}\label{e:3.4}
\phi_{{\gamma_0}}:\R\times B^n_{2\rho}(0)\to
M,\;(t,x)\mapsto\exp_{\gamma_0(t)}\left(\sum^n_{i=1}x_i
e_i(t)\right),
 \end{eqnarray}
which satisfies
 $\phi_{{\gamma_0}}(t+1,x)=\phi_{{\gamma_0}}(t, (E_{{\gamma_0}}
x^T)^T)$ and
\begin{eqnarray*}
d\phi_{{\gamma_0}}(t+1,x)[(1,v)]=d\phi_{{\gamma_0}}(t,
(E_{{\gamma_0}} x^T)^T)[(1, (E_{{\gamma_0}} v^T)^T)]
\end{eqnarray*}
for $(t,x)\in \R\times B^n_{2\rho}(0)$. This yields  a
coordinate chart around $\gamma_0$ on $\Lambda(M,\mathbb{I}_g)$,
\begin{eqnarray}\label{e:3.5}
\Phi_{{\gamma_0}}:{\bf B}_{2\rho}({H}_{{\gamma_0}}):=\{\xi\in
 H_{{\gamma_0}}\,|\,\|\xi\|_{1,2}<2\rho\} \to \Lambda(M,\mathbb{I}_g)
\end{eqnarray}
given by $\Phi_{{\gamma_0}}(\xi)(t)=\phi_{{\gamma_0}}(t,\xi(t))$ for
$\xi\in {\bf B}_{2\rho}({H}_{{\gamma_0}})$. $\Phi_{{\gamma_0}}$ is only $C^{k-3}$ by Lemma~\ref{lem:5.2},
and $d\Phi_{{\gamma_0}}(0)={\mathfrak{I}}_{{\gamma_0}}$.

Shrink $\varepsilon>0$ in (\ref{e:1.17}) so that
$\varepsilon<\sqrt{2}\rho$. Then
$$
{\rm
EXP}_{\gamma_0}\bigl(T\Lambda(M,\mathbb{I}_g)(\varepsilon)_{\gamma_0}\bigr)\subset\Phi_{{\gamma_0}}\bigl({\bf
B}_{2\rho}({H}_{{\gamma_0}})\bigr).
$$
(Indeed, for $v=\sum^n_{i=1}v_ie_i\in
T\Lambda(M,\mathbb{I}_g)(\varepsilon)_{\gamma_0})$ write
 ${\bf v}=(v_1,\cdots,v_n)^T$ then $\|{\bf v}\|_\infty\le\sqrt{2}\|{\bf
v}\|_{1,2}=\sqrt{2}\|v\|_1<\sqrt{2}\varepsilon<2\rho$.) Define
$L_{{\gamma_0}}:\R\times B^n_{2\rho}(0)\times\R^n\to\R$ by
\begin{eqnarray}\label{e:3.6}
 L_{{\gamma_0}}(t, x, v)=L\bigr(\phi_{{\gamma_0}}(t,x),
d\phi_{{\gamma_0}}(t,x)[(1,v)]\bigl).
\end{eqnarray}
Then $L_{{\gamma_0}}(t,0,0)\equiv c\;\forall t\in \R$,
\begin{eqnarray*}
L_{{\gamma_0}}(t+1, x, v)&=&L\bigr(\phi_{{\gamma_0}}(t+1,x),
d\phi_{{\gamma_0}}(t+1,x)[(1,v)]\bigl)\\
&=&L\bigr(\phi_{{\gamma_0}}(t, (E_{{\gamma_0}} x^T)^T),
d\phi_{{\gamma_0}}(t, (E_{{\gamma_0}} x^T)^T)[(1, (E_{{\gamma_0}} v^T)^T)]\bigl)\\
&=&L_{{\gamma_0}}\bigl(t, (E_{{\gamma_0}} x^T)^T, (E_{{\gamma_0}}
v^T)^T\bigr)\quad\hbox{and}\\
\partial_x L_{{\gamma_0}}(t+1, x, v)&=&\partial_x L_{{\gamma_0}}\bigl(t,
(E_{{\gamma_0}} x^T)^T,
(E_{{\gamma_0}} v^T)^T\bigr)E_{{\gamma_0}},\\
\partial_v L_{{\gamma_0}}\bigl(t+1, x, v)&=&\partial_v L_{{\gamma_0}}\bigl(t,
(E_{{\gamma_0}} x^T)^T, (E_{{\gamma_0}} v^T)^T\bigr)E_{{\gamma_0}}.
\end{eqnarray*}
So for any $\xi\in {\bf B}_{2\rho}({H}_{{\gamma_0}})$, since
$\xi(t+1)^T=E_{{\gamma_0}}\xi(t)^T$  and
$\dot\xi(t+1)^T=E_{{\gamma_0}}\dot\xi(t)^T$ a.e. in $\mathbb{R}$, we get
\begin{eqnarray}
 L_{{\gamma_0}}(t+1,\xi(t+1),\dot\xi(t+1))&=&
L_{{\gamma_0}}(t,
(E_{{\gamma_0}} \xi(t+1)^T)^T, (E_{{\gamma_0}} \dot\xi(t+1)^T)^T)\nonumber\\
&=& L_{{\gamma_0}}(t,\xi(t),\dot\xi(t)),\label{e:3.7}\\
\partial_x L_{{\gamma_0}}(t+1,\xi(t+1),\dot\xi(t+1))&=&\partial_x
L_{{\gamma_0}}(t,\xi(t),\dot\xi(t))E_{{\gamma_0}},\label{e:3.8}\\
\partial_v L_{{\gamma_0}}(t+1,\xi(t+1),\dot\xi(t+1))&=&\partial_v
L_{{\gamma_0}}(t,\xi(t),\dot\xi(t))E_{{\gamma_0}}.\label{e:3.9}
\end{eqnarray}
 Define  the action
functional ${\cal L}_{{\gamma_0}}: {\bf
B}_{2\rho}({H}_{{\gamma_0}})\to\R$ by
$$
{\cal L}_{{\gamma_0}}(\xi)=\int^1_0
L_{{\gamma_0}}(t,\xi(t),\dot\xi(t))dt.
$$
 Then
${\cal L}_{{\gamma_0}}={\cal L}\circ\Phi_{{\gamma_0}}$ on ${\bf
B}_{2\rho}({H}_{{\gamma_0}})$. It is $C^{2-0}$, and its gradient is
given by (B.34) and (B.35) in \cite{Lu3}, that is,
\begin{eqnarray}\label{e:3.10}
\nabla{\cal L}_{{\gamma_0}}(\xi)(t)
 &=&\frac{1}{2}\int^\infty_{t}
e^{t-s}\left[\partial_x
L_{{\gamma_0}}\bigl(s,\xi(s),\dot{\xi}(s)\bigr)-\mathfrak{R}^{\xi}(s)\right]\,ds
\nonumber\\
&+&\frac{1}{2}\int_{-\infty}^{t} e^{s-t}\left[\partial_x
L_{{\gamma_0}}\bigl(s,\xi(s),\dot{\xi}(s)\bigr)-\mathfrak{R}^{\xi}(s)\right]\,ds+\mathfrak{R}^{\xi}(t),
\end{eqnarray}
where as $2={\rm ord}(S_p)>{\rm ord}(S_{p+1})$ for some
$p\in\{0,\cdots,\sigma\}$,
\begin{eqnarray}\label{e:3.11}
\mathfrak{R}^\xi(t)= \int^t_0\partial_v
L_{\gamma_0}\bigl(s,\xi(s),\dot{\xi}(s)\bigr)ds+
\int^1_0\partial_v L_{\gamma_0}\bigl(s,\xi(s),\dot{\xi}(s)\bigr)ds\times &&\nonumber\\
\times \left(\oplus_{l\le
p}\frac{\sin\theta_l}{2-2\cos\theta_l}{\scriptscriptstyle\left(\begin{array}{cc}
0& -1\\
1& 0\end{array}\right)}-\frac{1}{2}I_{2p}\right)\oplus{\rm
diag}(a_{p+1}(t), \cdots, a_{\sigma}(t))&&
\end{eqnarray}
with $a_j(t)=\frac{2t{S}_j+ 2t+1-{S}_j}{4}$,
$j=p+1,\cdots,{\sigma}=n-p$. (As usual $p=0$ or $p={\sigma}$ means
there is no the first or second term in (\ref{e:3.11}).)

Corresponding to (\ref{e:2.3}) let
\begin{equation}\label{e:3.12}
{\cal Z}_{{\gamma_0}}:=\bigl\{(t,x,v)\in \R\times {\bf
B}_{2\rho}(\R^n)\times\R^n\,|\,
\partial_x\phi_{{\gamma_0}}(t,x)[v]=-\partial_t\phi_{{\gamma_0}}(t,x)\bigr\}.
\end{equation}
It  is  a submanifold in $\R\times {\bf B}_{2\rho}(\R^n)\times\R^n$
of dimension $n+1$, and also a closed subset. Note that
$L_{{\gamma_0}}$ is $C^k$ in $(\R\times {\bf
B}_{2\rho}(\R^n)\times\R^n)\setminus{\cal Z}_{{\gamma_0}}$.
Moreover,  there exist positive constants $C_1<C_2$  such that
for all $(t,x,v)\in \R\times {\bf
B}_{2\rho}(\R^n)\times\R^n\setminus{\cal Z}_{\gamma_0}$,
\begin{description}
\item[(LL1)] $\sum_{ij}\frac{\partial^2}{\partial v_i\partial
v_j}L_{{\gamma_0}}(t,x,v)u_iu_j\ge C_1|{\bf u}|^2\quad\forall {\bf
u}=(u_1,\cdots,u_n)\in\R^n$,

\item[(LL2)] $\Bigl| \frac{\partial^2}{\partial x_i\partial
x_j}L_{{\gamma_0}}(t,x,v)\Bigr|\le C_2(1+ |v|^2),\quad \Bigl|
\frac{\partial^2}{\partial x_i\partial v_j}L_{{\gamma_0}}(t,x,v)\Bigr|\le
C_2(1+
|v|),\quad\hbox{and}\\
 \Bigl| \frac{\partial^2}{\partial v_i\partial
v_j}L_{{\gamma_0}}(t,x,v)\Bigr|\le C_2$
\end{description}
(cf. \cite[(1) $\&$ (2)]{CaJaMa2}).

 Given
$m\in\N$ we call a curve
$$
\gamma\in C^m(\R, M, \mathbb{I}_g):=\{\gamma\in C^m(\R, M)\,|\,
\gamma(t+1)=\mathbb{I}_g(\gamma(t)),\;\forall t\in\R\}
$$
{\bf strongly regular} if $\{t\in \R\,|\, \dot\gamma(t)=0\}$ is empty.
Let $C^m(\R, M, \mathbb{I}_g)_{\rm reg}$ denote the subset of all
strongly regular curves in $C^m(\R, M, \mathbb{I}_g)$, and let
$$
C^m_{{\gamma_0}}(\R,\R^n):=\{\xi\in C^m(\R,
\R^n)\,|\,\xi(t+1)^T=E_{{\gamma_0}}\xi(t)^T,\;\forall t\in\R\},
$$
which is a Banach space with the usual $C^m$-norm
because $E_{{\gamma_0}}\in O(n)$.  Then
$$
\Phi_{{\gamma_0}}^{-1}\bigl(C^m(\R, M, \mathbb{I}_g)_{\rm
reg}\bigr)=\left\{\xi\in{\bf B}_{2\rho}({H}_{{\gamma_0}})\cap
C^m_{{\gamma_0}}(\R,\R^n)\,|\, (\check{\xi})^{-1}({\cal Z}_{{\gamma_0}})=\emptyset\right\},
$$
where $\check{\xi}:\R\to \R\times {\bf B}_{2\rho}(\R^n)\times\R^n$ is defined by
$\check{\xi}(t)=(t,\xi(t),\dot\xi(t))$. As in
(\ref{e:2.8}), we can choose $0<r<\rho/2$ so small that
\begin{equation}\label{e:3.13}
\R\times {\bf B}_{2r}(\R^n)\times {\bf B}_{2r}(\R^n)\subset \R\times
{\bf B}_{2\rho}(\R^n)\times\R^n\setminus{\cal Z}_{{\gamma_0}}.
\end{equation}
 Put $U={\bf B}_{2r}(H_{\gamma_0})$,
  $U_X:={\bf B}_{2r}(H_{\gamma_0})\cap
X_{\gamma_0}$ as an open subset of $X_{\gamma_0}$, and
\begin{eqnarray}\label{e:3.14}
U_{X,k}^{\rm reg}:=U_X\cap (\Phi_{{\gamma_0}}|_{{\bf B}_{2r}({H}_{{\gamma_0}})})^{-1}\bigl(C^{k-1}(\R, M,
\mathbb{I}_g)_{\rm reg}\bigr).
\end{eqnarray}
Then ${\bf B}_{2r}(X_{\gamma_0})\subset U_X$, and
$(\check{\xi})^{-1}({\cal Z}_{{\gamma_0}})=\emptyset,\;\forall\xi\in
U_{X,k}^{\rm reg}\cup{\bf B}_{2r}(X_{\gamma_0})$
 by (\ref{e:3.13}) and (\ref{e:3.14}).
  For $\xi\in U_{X,k}^{\rm reg}\cup{\bf B}_{2r}(X_{\gamma_0})\subset U_X$, from (\ref{e:3.10}) and (\ref{e:3.11}) we see that $t\mapsto\nabla{\cal
L}_{{\gamma_0}}(\xi)(t)$ is continuous differentiable, and
\begin{eqnarray*}
\frac{d}{dt}\nabla{\cal
L}_{{\gamma_0}}(\xi)(t)&=&\frac{e^{t}}{2}\int^\infty_te^{-s}\left(
\partial_x
L_{{\gamma_0}}\bigl(s,\xi(s),\dot{\xi}(s)\bigr)-\mathfrak{R}^{\xi}(s) \right)\, ds\nonumber\\
&-&\frac{e^{-t}}{2}\int^t_{-\infty}e^{s}\left(\partial_x
L_{{\gamma_0}}\bigl(s,\xi(s),\dot{\xi}(s)\bigr)-\mathfrak{R}^{\xi})(s)
\right)\, ds \nonumber\\
&+&\partial_v L_{{\gamma_0}}\bigl(t,\xi(t),\dot{\xi}(t)\bigr)-
\left(\int^{1}_0\partial_v
L_{{\gamma_0}}\bigl(s,\xi(s),\dot{\xi}(s)\bigr)ds\right)\times\nonumber\\
&&\times \biggl[\frac{d}{dt}\left(\oplus_{l\le
p}\frac{\sin\theta_l}{2-2\cos\theta_l}{\scriptscriptstyle\left(\begin{array}{cc}
0& -1\\
1&
0\end{array}\right)}-\frac{1}{2}I_{2p}\right)\nonumber\\
&&\qquad\oplus\, {\rm diag}(\dot{a}_{p+1}(t), \cdots,
\dot{a}_\sigma(t))\biggr]
\end{eqnarray*}
provided $2={\rm ord}(S_p)>{\rm ord}(S_{p+1})$ for some
$p\in\{0,\cdots,\sigma\}$, where $\mathfrak{R}^{\xi}$ is given by
(\ref{e:3.11}). Let  $\mathbb{A}_{\gamma_0}$ be the restriction of the gradient
  $\nabla{\cal L}_{{\gamma_0}}$ to $U_X$. Then
$\mathbb{A}_{\gamma_0}(\xi)\in X_{\gamma_0}$ for $\xi\in U_{X,k}^{\rm reg}\cup{\bf B}_{2r}(X_{\gamma_0})\subset U_X$, and $U_X\supseteq U_{X,k}^{\rm reg}\cup{\bf B}_{2r}(X_{\gamma_0})
\ni \xi\mapsto \mathbb{A}_{\gamma_0}(\xi)\in X_{\gamma_0}$
is continuous with respect to the $C^1$-topology.  Corresponding to Lemma~\ref{lem:2.3} we have

\begin{lemma}\label{lem:3.1}
 The map $U_X\supseteq U_{X,k}^{\rm reg}\cup{\bf B}_{2r}(X_{\gamma_0})
\ni \xi\mapsto \mathbb{A}_{\gamma_0}(\xi)\in X_{\gamma_0}$ is Frech\'et differentiable, and
$U_X\supseteq U_{X,k}^{\rm reg}\cup{\bf B}_{2r}(X_{\gamma_0})
\ni \xi\mapsto d\mathbb{A}_{\gamma_0}(\xi)\in \mathscr{L}(X_{\gamma_0})$
 is continuous. In particular, $\mathbb{A}_{\gamma_0}$ restricts
to a $C^{k-3}$ map from ${\bf B}_{2r}(X_{\gamma_0})$ to $X_{\gamma_0}$.
\end{lemma}

Hence the restriction  of ${\cal L}_{{\gamma_0}}$
to $U_X$, denoted by ${\cal L}^{X}_{{\gamma_0}}$,  has second Frech\'et derivative at each $\zeta\in
U_{X,k}^{\rm reg}\cup {\bf B}_{2r}(X_{\gamma_0})$ given by
\begin{eqnarray*}
 d^2 {\cal L}^{X}_{{\gamma_0}}
  (\zeta)[\xi,\eta]
   = \int_0^{1} \Bigl(\!\! \!\!\!&&\!\!\!\!\!\partial_{vv}
     L_{{\gamma_0}}\bigl(t,\zeta(t),\dot{\zeta}(t)\bigr)
\bigl[\dot{\xi}(t), \dot{\eta}(t)\bigr]
+ \partial_{qv}
  L_{{\gamma_0}}\bigl(t,\zeta(t), \dot{\zeta}(t)\bigr)
\bigl[\xi(t), \dot{\eta}(t)\bigr]\nonumber \\
&& + \partial_{vq}
  L_{{\gamma_0}}\bigl(t,\zeta(t),\dot{\zeta}(t)\bigr)
\bigl[\dot{\xi}(t), \eta(t)\bigr] \nonumber\\
&&+  \partial_{qq} L_{{\gamma_0}}\bigl(t,\zeta(t),
\dot{\zeta}(t)\bigr) \bigl[\xi(t), \eta(t)\bigr]\Bigr) \, dt
\end{eqnarray*}
for any  $\xi,\eta\in X_{{\gamma_0}}$. From the last expression
and ({\bf LL1})-({\bf LL2}) below (\ref{e:3.12}), it
easily follows that for any $\zeta\in U_{X,k}^{\rm reg}\cup {\bf B}_{2r}(X_{\gamma_0})$ there
exists $C(\zeta)>0$ such that
\begin{equation}\label{e:3.15}
\big|d^2 {\cal L}^{X}_{{\gamma_0}} (\zeta)[\xi,\eta]\big|\le
  C(\zeta)\|\xi\|_{1,2}\cdot\|\eta\|_{1,2}
\quad\forall \xi,\eta\in X_{{\gamma_0}}.
\end{equation}
It follows that there exists a   map $\mathbb{B}_{\gamma_0}:
 U_{X,k}^{\rm reg}\cup {\bf B}_{2r}(X_{\gamma_0})\to
\mathscr{L}_s({H}_{{\gamma_0}})$ such that
\begin{equation}\label{e:3.16}
\langle d\mathbb{A}_{\gamma_0}(\zeta)[\xi], \eta\rangle_{1,2}=d^2
{\cal L}^{X}_{{\gamma_0}}
  (\zeta)[\xi,\eta]=\langle\mathbb{B}_{\gamma_0}(\zeta)\xi,
  \eta\rangle_{1,2}\quad\forall \xi, \eta\in X_{{\gamma_0}}.
\end{equation}
  Moreover,  by \cite{Lu1} $\mathbb{B}_{\gamma_0}(0)\in\mathscr{L}(H_{\gamma_0})$ is a self-adjoint Fredholm
operator, and the map $\mathbb{B}_{\gamma_0}$  has a decomposition
$\mathbb{B}_{\gamma_0}(\zeta)=P_{{\gamma_0}}(\zeta)+
Q_{{\gamma_0}}(\zeta)\;\forall \zeta\in U_{X,k}^{\rm reg}\cup {\bf B}_{2r}(X_{\gamma_0})$, where
$P_{{\gamma_0}}(\zeta)\in\mathscr{L}_s(H_{\gamma_0})$ is a positive
definitive linear operator defined by
\begin{eqnarray}\label{e:3.17}
\langle P_{{\gamma_0}}(\zeta)\xi,
  \eta\rangle_{1,2}
   = \int_0^1 \left(\partial_{vv}L_{{\gamma_0}}\bigl(t,\zeta(t),\dot\zeta(t)\bigr)
[\dot\xi(t), \dot\eta(t)]+ \bigl(\xi(t), \eta(t)\bigr)_{\R^n}\right)
\, dt,
\end{eqnarray}
and $Q_{{\gamma_0}}(\zeta)\in\mathscr{L}_s(H_{\gamma_0})$ is a compact
linear operator. As in the proof of Lemma~\ref{lem:2.3+}, from this and ({\bf LL1})-({\bf LL2})
we may deduce
\begin{lemma}\label{lem:3.2}
The operators $P_{{\gamma_0}}$ and $Q_{{\gamma_0}}$ have the following properties:
\begin{description}
\item[($\mathfrak{P}1$)] For any
sequence $(\xi_j)\subset U_{X,k}^{\rm reg}\cup {\bf B}_{2r}(X_{\gamma_0})$ with $\|\xi_j\|_{1,2}\to
0$ it holds that $\|P_{{\gamma_0}}(\xi_j)\eta-
P_{{\gamma_0}}(0)\eta\|_{1,2}\to 0$ for any
$\eta\in H_{\gamma_0}$;

\item[($\mathfrak{P}2$)] There exist constants $\varepsilon_{\gamma_0}$ and
$C_{\gamma_0}>0$ such that
$$
\langle P_{{\gamma_0}}(\xi)\eta, \eta\rangle_{1,2}\ge
C_{\gamma_0}\|\eta\|^2_{1,2}\quad\forall \eta\in H_{\gamma_0}
$$
for any  $\xi\in U_{X,k}^{\rm reg}\cup {\bf B}_{2r}(X_{\gamma_0})$ with $\|\xi\|_{1,2}\le
\varepsilon_{\gamma_0}$;
\item[($\mathfrak{Q}$)] $\|Q_{{\gamma_0}}(\xi_j)-
Q_{{\gamma_0}}(0)\|_{\mathscr{L}_s(H_{{\gamma_0}})}\to 0$ for any
sequence $(\xi_j)\subset U_{X,k}^{\rm reg}\cup {\bf B}_{2r}(X_{\gamma_0})$ with $\|\xi_j\|_{1,2}\to
0$.
\end{description}
\end{lemma}

Let  $H^+_{{\gamma_0}}$, $H^-_{{\gamma_0}}$ and $H^0_{{\gamma_0}}$
denote the positive definite, negative definite and null spaces of
$\mathbb{B}_{\gamma_0}(0)$. Note that $H^-_{{\gamma_0}}\subset
X_{{\gamma_0}}$ and $H^0_{{\gamma_0}}\subset X_{{\gamma_0}}$  as Banach
subspaces of $X_{{\gamma_0}}$, they are denoted by
 $X^-_{{\gamma_0}}$ and $X^0_{{\gamma_0}}$. Set
$X^+_{{\gamma_0}}:=X_{{\gamma_0}}\cap H^+_{{\gamma_0}}$. We arrive
at a direct sum decomposition of Banach spaces,
$X_{{\gamma_0}}=X^0_{{\gamma_0}}\dot{+} X^+_{{\gamma_0}}\dot{+}
X^-_{{\gamma_0}}$. (From (\ref{e:3.26}) it follows that the Morse index
$m^-({\cal L}, {\cal O})$ and nullity $m^0({\cal L}, {\cal O})$ of ${\cal O}$ are equal
to  $\dim X^-_{{\gamma_0}}$ and $\dim X^0_{{\gamma_0}}-1$,
respectively). As in (\ref{e:2.17})  we have  a positive number
$a_{\gamma_0}>0$ such that  the set $[-2a_{\gamma_0},
2a_{\gamma_0}]\setminus\{0\}$ is disjoint with the spectrum of
$\mathbb{B}_{\gamma_0}(0)$, which implies
 \begin{equation}\label{e:3.18}
\left.\begin{array}{ll}
 &\langle \mathbb{B}_{\gamma_0}(0)\xi, \xi\rangle_{1,2}\ge
2a_{\gamma_0}\|\xi\|^2_{1,2}\quad\forall \xi\in H_{{\gamma_0}}^+,\\
&\langle \mathbb{B}_{\gamma_0}(0)\xi, \xi\rangle_{1,2}\le
-2a_{\gamma_0}\|\xi\|^2_{1,2}\quad\forall \xi\in H_{{\gamma_0}}^-.
\end{array}\right\}
 \end{equation}

Let $\mathbb{P}^\ast_{{\gamma_0}}$ be the orthogonal projections
from $H_{{\gamma_0}}=H^-_{{\gamma_0}}\oplus H^0_{{\gamma_0}}\oplus
H^+_{{\gamma_0}}$ onto $H^\ast_{{\gamma_0}}$, $\ast=-,0,+$. Choose a
basis of $H^-_{{\gamma_0}}\oplus H^0_{{\gamma_0}}$, $\mathbbm{e}_i$,
$i=1,\cdots, 1+ m^-({\cal L}, {\cal O})+ m^0({\cal L}, {\cal O})$.   As in the proofs of
\cite[Lemmas~3.3,3.4]{Lu2} we may use  (\ref{e:3.18}), Lemma~\ref{lem:3.2}
 to derive the following results corresponding to
Lemmas~\ref{lem:2.4},~\ref{lem:2.5}.

\begin{lemma}\label{lem:3.3}
 There exists a function
 $\omega:{U}_{X,k}^{\rm reg}\cup {\bf B}_{2r}(X_{\gamma_0})\to [0, \infty)$
  such that  $\omega(\zeta)$ uniformly converges to $0$ as
  $\zeta\in {U}_{X,k}^{\rm reg}\cup {\bf B}_{2r}(X_{\gamma_0})$ and $\|\zeta\|_{1,2}\to 0$, and that
$$
|\langle \mathbb{B}_{\gamma_0}(\zeta)\xi, \eta\rangle_{1,2}- \langle
\mathbb{B}_{\gamma_0}(0)\xi, \eta\rangle_{1,2} |\le \omega(\zeta)
\|\xi\|_{1,2}\cdot\|\eta\|_{1,2}
$$
for any $\zeta\in {U}_{X,k}^{\rm reg}\cup {\bf B}_{2r}(X_{\gamma_0})$,  $\xi\in
H^0_{{\gamma_0}}\oplus H^-_{{\gamma_0}}$ and $\eta\in
H_{{\gamma_0}}$.
\end{lemma}

\begin{lemma}\label{lem:3.4}
By shrinking  ${U}={\bf B}_{2r}(H_{{\gamma_0}})$ (or equivalently
$r>0$), we may obtain  a number $a_{\gamma_0}'\in (0,
2a_{\gamma_0}]$ such that for any $\zeta\in {U}_{X,k}^{\rm reg}\cup {\bf B}_{2r}(X_{\gamma_0})$ we
have $\omega(\zeta)<\min\{a_{\gamma_0}', a_{\gamma_0}\}/2$ and
\begin{description}
\item[(i)] $\langle \mathbb{B}_{\gamma_0}(\zeta)\xi, \xi\rangle_{1,2}\ge a_{\gamma_0}'\|\xi\|^2_{1,2}\;\forall
\xi\in H^+_{{\gamma_0}}$,
\item[(ii)] $|\langle \mathbb{B}_{\gamma_0}(\zeta)\xi,\eta\rangle_{1,2}|\le
\omega(\zeta)\|\xi\|_{1,2}\cdot\|\eta\|_{1,2}\; \forall \xi\in
H^+_{{\gamma_0}}, \forall \eta\in H^-_{{\gamma_0}}\oplus
H^0_{{\gamma_0}}$,
\item[(iii)] $\langle \mathbb{B}_{\gamma_0}(\zeta)\xi,\xi\rangle_{1,2}\le-a_{\gamma_0}\|\xi\|^2_{1,2}\;
\forall \xi\in H^-_{{\gamma_0}}$.
\end{description}
\end{lemma}

Note that $\mathbb{A}_{\gamma_0}$ is $C^{k-3}$ in ${\bf
B}_{2r}(X_{{\gamma_0}})\subset {U}_{X}$ and that
$H_{{\gamma_0}}$ and $X_{{\gamma_0}}$ induce equivalent norms on
$X^0_{{\gamma_0}}={\rm Ker}(d\mathbb{A}_{\gamma_0}(0))$. It follows
from the
 implicit function theorem  that there exist
$\tau\in (0, r/2)$ and a $C^{k-3}$-map $h: \overline{{\bf
B}_\tau(H^0_{{\gamma_0}})}\to X^-_{{\gamma_0}}\oplus
X^+_{{\gamma_0}}$ with $h(0)=0$ and $dh(0)=0$, such that for all
$\xi\in \overline{{\bf B}_\tau(H^0_{{\gamma_0}})}$,
\begin{equation}\label{e:3.19}
\xi+  h(\xi)\in {\bf B}_{r}(X_{{\gamma_0}})\quad\hbox{and}\quad
(I-\mathbb{P}^0_{{\gamma_0}}) \mathbb{A}_{\gamma_0}(\xi+ h(\xi))=0.
\end{equation}

Define  ${\bf J}:\overline{{\bf
B}_\tau(H^0_{{\gamma_0}})}\times \bigl({\bf
B}_\tau(H^-_{{\gamma_0}})\oplus
 {\bf B}_\tau(H^+_{{\gamma_0}})\bigr)\to\R$ by
\begin{equation}\label{e:3.20}
{\bf J}(\xi, \zeta)={\cal L}_{{\gamma_0}}(\xi+ h(\xi)+ \zeta)-{\cal
L}_{{\gamma_0}}(\xi+  h(\xi))
\end{equation}
for  $\zeta\in {\bf B}_\tau(H^-_{{\gamma_0}})\oplus {\bf
B}_\tau(H^+_{{\gamma_0}})$ and $\xi\in \overline{{\bf
B}_\tau(H^0_{{\gamma_0}})}$. It is $C^{2-0}$, and satisfies
\begin{eqnarray}
&&D_2{\bf J}(\xi,\zeta)[\eta] =\langle
(I-\mathbb{P}^0_{{\gamma_0}})\nabla{\cal L}_{{\gamma_0}}(\xi+ h(\xi)+
\zeta), \eta\rangle_{1,2},
\label{e:3.21}\\
&& {\bf J}(\xi, 0)=0\quad\hbox{and}\quad D_2{\bf J}(\xi,
0)[\eta]=0\label{e:3.22}
\end{eqnarray}
for every $(\xi,\zeta)\in \overline{{\bf
B}_\tau(H^0_{{\gamma_0}})}\times \bigl({\bf
B}_\tau(H^-_{{\gamma_0}})\oplus
 {\bf B}_\tau(H^+_{{\gamma_0}})\bigr)$ and $\eta\in H^+_{{\gamma_0}}+H^-_{{\gamma_0}}$.
Moreover, by the last claim of Lemma~\ref{lem:3.1} the functional
\begin{equation}\label{e:3.23}
{\cal L}^\circ_{{\gamma_0}}: \overline{{\bf
B}_\tau(H^0_{{\gamma_0}})}\to\R,\; \xi\mapsto {\cal
L}_{{\gamma_0}}(\xi+ h(\xi))
\end{equation}
is $C^{k-3}$, and $d{\cal L}^\circ_{{\gamma_0}}(\xi)[\eta]=\langle
\mathbb{A}_{\gamma_0}(\xi+
 h(\xi)), \eta\rangle_{1,2}\;\forall \eta\in H^0_{{\gamma_0}}$ and
$d^2{\cal L}^\circ_{{\gamma_0}}(0)=0$.

Set ${U}_{X,k}:={\bf B}_{2r}(H_{{\gamma_0}})\cap
C^{k-1}_{{\gamma_0}}(\R,\R^n)$ with $2r<\rho$, and  consider the $C^{k-2}$
evaluation map
\begin{equation}\label{e:3.24}
\mathfrak{E}:\R\times {U}_{X,k}\to \R\times {\bf
B}_{2\rho}(\R^n)\times\R^n,\quad (t,\xi)\mapsto
(t,\xi(t),\dot\xi(t)).
\end{equation}
 Let $\mathfrak{E}_\xi:\R\to \R\times {\bf
B}_{2\rho}(\R^n)\times\R^n$ be defined by
$\mathfrak{E}_\xi(t)=\mathfrak{E}(t,\xi)$ for $\xi\in {U}_{X,k}$. Note that
$U_{X,k}^{\rm reg}=U_{X,k}\cap (\Phi_{{\gamma_0}}|_{{\bf B}_{2r}({H}_{{\gamma_0}})})^{-1}\bigl(C^{k-1}(\R, M,
\mathbb{I}_g)_{\rm reg}\bigr)$ by (\ref{e:3.14}).
By the same proof as that of Lemma~\ref{lem:2.2}  we may obtain the corresponding result.

\begin{lemma}\label{lem:3.5}
 There exists a residual subset
 $C^{k-1}_{{\gamma_0}}(\R,\R^n)^{\circ}_{\rm reg}$ in  $C^{k-1}_{{\gamma_0}}(\R,\R^n)$
 such that\linebreak ${\bf B}_{2r}(H_{{\gamma_0}})\cap C^{k-1}_{{\gamma_0}}(\R,\R^n)^{\circ}_{\rm reg}\supseteq
 {U}_{X,k}^{\rm reg}$ with equality in case $n>1$, and
  that  for every $\xi\in C^{k-1}_{{\gamma_0}}(\R,\R^n)^{\circ}_{\rm reg}$ the set
$(\mathfrak{E}_{\xi})^{-1}({\cal Z}_{{\gamma_0}})$ is a $C^{k-2}$ submanifold of
dimension $1-n$, and so empty (in case $n>1$) and at most a finite
set (in case $n=1$).
 Moreover, if $\xi_i\in C^{k-1}_{{\gamma_0}}(\R,\R^n)^{\circ}_{\rm reg}$, $i=1,2$, then for a generic $C^{k-1}$ path
$\mathbbm{p}:[0,1]\to C^{k-1}_{{\gamma_0}}(\R,\R^n)$ connecting $\xi_0$ to $\xi_1$ the
set $\{(s,t)\in [0, 1]\times \R\,|\,
\bigr(\mathfrak{E}(\mathbbm{p}(s)\bigl)(t)\in{\cal
Z}_{{\gamma_0}}\}$ is a $C^{k-2}$ submanifold of $[0, 1]\times \R$ of
dimension $2-n$, and so empty for $n>2$, and at most a finite set
for $n=2$.
\end{lemma}

Since the $C^{k-3}$-map $h: \overline{{\bf B}_\tau(H^0_{{\gamma_0}})}\to
X^-_{{\gamma_0}}\oplus X^+_{{\gamma_0}}$ satisfies $h(0)=0$ and
$dh(0)=0$, shrinking $\tau>0$ we can guarantee that the map
$$
\mathfrak{M}:\overline{{\bf B}_\tau(H^0_{{\gamma_0}})}\oplus
\bigl(H_{{\gamma_0}}^-\oplus H^+_{{\gamma_0}}\bigr)\to H_{{\gamma_0}}^-\oplus
H^+_{{\gamma_0}}
$$
given by $\mathfrak{M}(\xi+\zeta)= \xi+ h(\xi)+\zeta$ for $\xi\in
\overline{{\bf B}_\tau(H^0_{{\gamma_0}})}$ and $\zeta\in
H_{{\gamma_0}}^-\oplus H^+_{{\gamma_0}}$,
 is a $C^{k-3}$ diffeomorphism onto a neighborhood of $0$ in
${U}={\bf B}_{2r}(H_{{\gamma_0}})$. So
${W}_{X,k}^{\rm
reg}:=(\mathfrak{M})^{-1}({U}_{X,k}^{\rm reg})$
 is a residual subset
in $\overline{{\bf B}_\tau(H^0_{{\gamma_0}})}\oplus
{\bf B}_\tau(H_{{\gamma_0}}^+)\oplus {\bf B}_\tau(H^0_{{\gamma_0}})$. As in the proof of
Proposition~\ref{prop:2.6} we can use this and
Lemmas~\ref{lem:3.3},~\ref{lem:3.4} to prove

\begin{proposition}\label{prop:3.6}
\begin{enumerate}
\item[{\bf (i)}]  For any $(\xi, \zeta)\in \overline{{\bf
B}_\tau(H^0_{{\gamma_0}})}\times {\bf B}_\tau(H^+_{{\gamma_0}})$,
$\eta_1, \eta_2\in {\bf B}_\tau(H^-_{{\gamma_0}})$,
 $$
 (D_2{\bf J}(\xi,
\zeta+ \eta_2)-D_2{\bf J}(\xi, \zeta+
\eta_1))[\eta_2-\eta_1]\le-\frac{a_{\gamma_0}}{2}
\|\eta_2-\eta_1\|_{1,2};
$$
\item[{\bf (ii)}]  For any $(\xi, \zeta, \eta)\in
 \overline{{\bf B}_\tau(H^0_{{\gamma_0}})}\times {\bf B}_\tau(H^+_{{\gamma_0}})
 \times {\bf B}_\tau(H^-_{{\gamma_0}})$ we have
$$
D_2{\bf J}(\xi,
\zeta+\eta)[\zeta-\eta]\ge\frac{a_{\gamma_0}'}{2}\|\zeta\|^2_{1,2}+
\frac{a_{\gamma_0}}{2}\|\eta\|^2_{1,2};
$$
in particular $D_2{\bf J}(\xi,
\zeta)[\zeta]\ge\frac{a_{\gamma_0}'}{2}\|\zeta\|^2_{1,2}$
 for any $(\xi, \zeta)\in
\overline{{\bf B}_\tau(H^0_{{\gamma_0}})}\times {\bf
B}_\tau(H^+_{{\gamma_0}})$.
\end{enumerate}
\end{proposition}

If $\dim H_{{\gamma_0}}^0=0$ this result also holds with
 ${\bf J}(0, \zeta)={\cal
L}_{{\gamma_0}}(\zeta)$ for $\zeta\in {\bf
B}_\tau(H^-_{{\gamma_0}})\oplus
 {\bf B}_\tau(H^+_{{\gamma_0}})$.

\noindent{\bf Step 2}. Observe that the critical set of ${\cal
L}_{\gamma_0}$ is an one-dimensional $C^{k-1}$ critical manifold $S:=
\Phi^{-1}_{{\gamma_0}}\bigl({\cal O}\cap{\rm
Im}(\Phi_{\gamma_0})\bigr)$ containing $0$ as an interior point, and
that $T_{\gamma_0}{\cal O}=\dot\gamma_0\R\subset
T_{\gamma_0}\Lambda(M,\mathbb{I}_g)$. Since
$d\Phi_{\gamma_0}(0)=\mathfrak{I}_{{\gamma_0}}$
  is an isomorphism there exists a
unique $\zeta_0\in H_{\gamma_0}$ satisfying
${\mathfrak{I}}_{{\gamma_0}}(\zeta_0)=\dot\gamma_{0}$,
  that is, $\dot\gamma_0(t)=\sum^n_{j=1}\zeta_{0j}(t)e_j(t),\;\forall t\in\R$.
Hence $\zeta_{0j}(t)=g(\dot\gamma_0(t), e_j(t))\;\forall t\in \R$
and $j=1,\cdots,n$. Clearly,  $\zeta_0\in
C_{\gamma_0}^{k-1}(\R,\R^n)\subset X_{\gamma_0}$  and
$T_0S=\R\zeta_0$. The normal space of $S$ at $0\in S$ is  the
orthogonal complementary of $\R\zeta_0$ in the Hilbert space
$H_{\gamma_0}$, denoted by $H_{{\gamma_0},0}$. Set
$X_{{\gamma_0},0}=H_{{\gamma_0},0}\cap X_{\gamma_0}$ and
\begin{eqnarray*}
\mathscr{U}={\bf B}_{2r}(H_{{\gamma_0},0})={\bf
B}_{2r}(H_{\gamma_0})\cap H_{{\gamma_0},0},
\end{eqnarray*}
moreover, define the open subset of $X_{{\gamma_0},0}$
\begin{eqnarray*}
 \mathscr{U}_{X}:={\bf B}_{2r}(H_{\gamma_0})\cap
X_{{\gamma_0},0}.
\end{eqnarray*}
Let ${\cal L}_{{\gamma_0},0}$
be the restriction of ${\cal L}_{{\gamma_0}}$ to $\mathscr{U}$. Then
$$
\nabla{\cal L}_{{\gamma_0},0}(\xi)=\nabla{\cal
L}_{{\gamma_0}}(\xi)-\frac{\langle \nabla{\cal L}_{{\gamma_0}}(\xi),
\zeta_0\rangle_{1,2}}{\| \zeta_0\|_{1,2}^2}
\zeta_0,\quad\forall\xi\in \mathscr{U},
$$
and so the restriction of $\nabla{\cal L}_{{\gamma_0},0}$ to
$\mathscr{U}_{X}$, denoted by $\mathbb{A}_{\gamma_0,0}$, is given by
\begin{equation}\label{e:3.25}
\mathbb{A}_{\gamma_0,0}(\xi)=\mathbb{A}_{\gamma_0}(\xi)-\frac{\langle
 \mathbb{A}_{\gamma_0}(\xi),
\zeta_0\rangle_{1,2}}{\| \zeta_0\|_{1,2}^2} \zeta_0,\quad\forall
\xi\in \mathscr{U}_{X}.
\end{equation}
By Lemma~\ref{lem:3.1}, $\mathbb{A}_{\gamma_0,0}$ restricts to a $C^{k-3}$ map from ${\bf
B}_{2r}(X_{{\gamma_0},0})\subset{\bf
B}_{2r}(X_{{\gamma_0}})\cap\mathscr{U}_{X}$ to $X_{{\gamma_0},0}$.
 Let ${\cal L}^{X}_{{\gamma_0},0}$ be the restriction of ${\cal L}^{X}_{{\gamma_0}}$ to
 $\mathscr{U}_{X}$. It restricts to a $C^{k-2}$ functional on
 ${\bf B}_{2r}(X_{{\gamma_0},0})$.
 For  each $\xi\in {\bf B}_{2r}(X_{{\gamma_0},0})$,  let
$\mathbb{B}_{\gamma_0,0}(\xi)$  denote both the continuous symmetric
bilinear extension  form of $d^2({\cal
L}^{X}_{{\gamma_0},0})(\xi)$  on ${H}_{{\gamma_0},0}$ and the
corresponding self-adjoint operator. Then
$$
\mathbb{B}_{\gamma_0,0}(\xi)\eta=\mathbb{B}_{\gamma_0}(\xi)\eta-\frac{\langle
\mathbb{B}_{\gamma_0}(\xi)\eta,
\zeta_0\rangle_{1,2}}{\|\zeta_0\|_{1,2}^2}\zeta_0\quad\forall\eta\in
{H}_{{\gamma_0},0}.
$$
From this, (\ref{e:3.16}) and (\ref{e:3.25}) it follows that for $\xi\in
{\bf B}_{2r}(X_{{\gamma_0},0})$ and $\zeta,\eta\in
X_{{\gamma_0},0}$,
$$
\langle d\mathbb{A}_{\gamma_0,0}(\xi)[\zeta],
\eta\rangle_{1,2}=d^2({\cal
L}^{X}_{{\gamma_0},0})(\xi)[\zeta,\eta]=\langle\mathbb{B}_{\gamma_0,0}(\xi)\zeta,
  \eta\rangle_{1,2}.
$$
 Let $\Pi_{{\gamma_0}}:H_{\gamma_0}\to
H_{{\gamma_0},0}$ be the orthogonal projection given by
$$
\Pi_{{\gamma_0}}\eta=\eta-\frac{\langle \eta,
\zeta_0\rangle_{1,2}}{\|\zeta_0\|_{1,2}^2}\zeta_0.
$$
Then
$\mathbb{B}_{{\gamma_0},0}(\xi)=\Pi_{{\gamma_0}}\circ\mathbb{B}_{\gamma_0}(\xi)|_{H_{{\gamma_0},0}},\;\forall
\xi\in {\bf B}_{2r}(X_{{\gamma_0},0})$. Let
$H^+_{{\gamma_0},0}$, $H^-_{{\gamma_0},0}$ and
 $H^0_{{\gamma_0},0}$ be the
positive definite, negative definite and null spaces of
 $\mathbb{B}_{\gamma_0,0}(0)$. We have
$$
H^+_{{\gamma_0}}=H^+_{{\gamma_0},0},\quad
H^-_{{\gamma_0}}=H^-_{{\gamma_0},0}\quad\hbox{and}\quad
H^0_{{\gamma_0}}=H^0_{{\gamma_0},0}\oplus{\R}\zeta_0.
$$
 Note that $H^-_{{\gamma_0},0}\subset
X_{{\gamma_0},0}$ and $H^0_{{\gamma_0},0}\subset X_{{\gamma_0},0}$.
They are denoted by $X^-_{{\gamma_0},0}$ and $X^0_{{\gamma_0},0}$
as Banach subspaces of $X_{{\gamma_0},0}$. Set
$X^+_{{\gamma_0},0}:=X_{{\gamma_0}}\cap H^+_{{\gamma_0},0}=X_{{\gamma_0},0}\cap H^+_{{\gamma_0},0}$. We
obtain a direct sum decomposition of Banach spaces,
$X_{{\gamma_0},0}=X^0_{{\gamma_0},0}\dot{+} X^+_{{\gamma_0},0}\dot{+}
X^-_{{\gamma_0},0}$. Moreover, (\ref{e:3.36}) and (\ref{e:3.37}) imply
\begin{equation}\label{e:3.26}
m^-({\cal L}, {\cal O})=\dim
X^-_{{\gamma_0},0}\quad\hbox{and}\quad
m^0({\cal L}, {\cal O})=\dim X^0_{{\gamma_0},0}.
\end{equation}
Let $h_{0}$ be the restriction to $\overline{{\bf
B}_\tau(H^0_{{\gamma_0},0})}\subset \overline{{\bf
B}_\tau(H^0_{{\gamma_0}})}$ of the map $h$ in (\ref{e:3.19}). It is
a $C^{k-3}$ map to $X^-_{{\gamma_0},0}\oplus
X^+_{{\gamma_0},0}=X^-_{{\gamma_0}}\oplus X^+_{{\gamma_0}}$ and
satisfies $h_{0}(0)=0$, $dh_{0}(0)=0$ and
\begin{equation}\label{e:3.27}
\xi+  h_{0}(\xi)\in {\bf
B}_{r}(X_{{\gamma_0},0})\quad\hbox{and}\quad
(I-\mathbb{P}^0_{{\gamma_0},0}) \mathbb{A}_{\gamma_0,0}(\xi+
h_{0}(\xi))=0
\end{equation}
for each $\xi\in \overline{{\bf B}_\tau(H^0_{{\gamma_0},0})}$ by
(\ref{e:3.19}), where
$\mathbb{P}^\ast_{{\gamma_0},0}:H_{{\gamma_0},0}\to
H^\ast_{{\gamma_0},0}$, $\ast=0,-,+$, are the orthogonal
projections.
 Functional ${\bf J}$ in (\ref{e:3.20}) restricts to
\begin{equation}\label{e:3.28}
{\bf J}_{0}:\overline{{\bf B}_\tau(H^0_{{\gamma_0},0})}\times \left({\bf
B}_\tau(H^-_{{\gamma_0},0})\oplus
 {\bf B}_\tau(H^+_{{\gamma_0},0})\right)\to\R,
 \end{equation}
 that is, ${\bf J}_{0}(\xi, \zeta)={\cal L}_{{\gamma_0},0}(\xi+
h_{0}(\xi)+ \zeta)-{\cal L}_{{\gamma_0},0}(\xi+ h_{0}(\xi))$ for
$\xi\in \overline{{\bf B}_\tau(H^0_{{\gamma_0},0})}$ and $\zeta\in
{\bf B}_\tau(H^-_{{\gamma_0},0})\oplus {\bf
B}_\tau(H^+_{{\gamma_0},0})$. It is $C^{2-0}$, and
(\ref{e:3.21})-(\ref{e:3.21}) imply
\begin{eqnarray}
&&D_2{\bf J}_{0}(\xi,\zeta)[\eta] =\langle
(I-\mathbb{P}^0_{{\gamma_0},0})\nabla{\cal L}_{{\gamma_0},0}(\xi+
h_{0}(\xi)+ \zeta), \eta\rangle_{1,2},
\label{e:3.29}\\
&& {\bf J}_{0}(\xi, 0)=0\quad\hbox{and}\quad D_2{\bf J}_{0}(\xi,
0)[\eta]=0\label{e:3.30}
\end{eqnarray}
for every $(\xi,\zeta)\in \overline{{\bf
B}_\tau(H^0_{{\gamma_0},0})}\times {\bf
B}_\tau(H^-_{{\gamma_0},0})\oplus
 {\bf B}_\tau(H^+_{{\gamma_0},0})$ and $\eta\in H^+_{{\gamma_0},0}+H^-_{{\gamma_0},0}$.
Moreover the functional ${\cal L}^\circ_{{\gamma_0}}$ in
(\ref{e:3.23}) restricts to a $C^{k-2}$ one
\begin{equation}\label{e:3.31}
{\cal L}^\circ_{{\gamma_0},0}: \overline{{\bf
B}_\tau(H^0_{{\gamma_0},0})}\to\R,\; \xi\mapsto {\cal
L}_{{\gamma_0}}(\xi+ h_{0}(\xi)),
\end{equation}
which satisfies $d{\cal L}^\circ_{{\gamma_0},0}(\xi)[\eta]=\langle
\mathbb{A}_{\gamma_0,0}(\xi+
 h_{0}(\xi)), \eta\rangle_{1,2}\;\forall \eta\in H^0_{{\gamma_0},0}$, and
$d^2{\cal L}^\circ_{{\gamma_0},0}(0)=0$.
From Proposition~\ref{prop:3.6} we immediately obtain

\begin{proposition}\label{prop:3.7}
\begin{enumerate}
\item[{\bf (i)}]  For any $(\xi, \zeta)\in \overline{{\bf
B}_\tau(H^0_{{\gamma_0},0})}\times {\bf B}_\tau(H^+_{{\gamma_0},0})$,
$\eta_1, \eta_2\in {\bf B}_\tau(H^-_{{\gamma_0},0})$,
 $$
 (D_2{\bf J}_0(\xi,
\zeta+ \eta_2)-D_2{\bf J}_0(\xi, \zeta+
\eta_1))[\eta_2-\eta_1]\le-\frac{a_{\gamma_0}}{2}
\|\eta_2-\eta_1\|_{1,2};
$$
\item[{\bf (ii)}]  For any $(\xi, \zeta, \eta)\in
 \overline{{\bf B}_\tau(H^0_{{\gamma_0},0})}\times {\bf B}_\tau(H^+_{{\gamma_0},0})
 \times {\bf B}_\tau(H^-_{{\gamma_0},0})$ we have
$$
D_2{\bf J}_0(\xi,
\zeta+\eta)[\zeta-\eta]\ge\frac{a_{\gamma_0}'}{2}\|\zeta\|^2_{1,2}+
\frac{a_{\gamma_0}}{2}\|\eta\|^2_{1,2};
$$
in particular $D_2{\bf J}_0(\xi,
\zeta)[\zeta]\ge\frac{a_{\gamma_0}'}{2}\|\zeta\|^2_{1,2}$
 for any $(\xi, \zeta)\in
\overline{{\bf B}_\tau(H^0_{{\gamma_0},0})}\times {\bf
B}_\tau(H^+_{{\gamma_0},0})$.
\end{enumerate}
\end{proposition}

\noindent{\bf Step 3}. For $\gamma\in{\cal O}$ let ${\cal
F}_{\gamma}={\cal F}|_{N{\cal O}(\varepsilon)_{\gamma}}$. Note that
$d\Phi_{{\gamma_0}}(0)=\mathfrak{I}_{{\gamma_0}}$ restricts to a
Hilbert space isomorphism $\mathfrak{I}_{{\gamma_0},0}:H_{{\gamma_0},0}\to N{\cal
O}_{\gamma_0}$. Take a positive number
$\epsilon<\min\{\tau,\delta/4\}$ such that  the map $\mathfrak{h}$
in (\ref{e:1.24}) satisfies
\begin{equation}\label{e:3.32}
 u+\mathfrak{h}_{\gamma_0}(u)\in XN{\cal O}(\delta/2)_{\gamma_0},\quad\forall u\in {\bf
 H}^0(B)(\epsilon)_{\gamma_0}.
\end{equation}
As in Section~\ref{sec:2.2} we have
\begin{eqnarray}\label{e:3.33}
&& \mathfrak{h}_{\gamma_0}(u)=\mathfrak{I}_{{\gamma_0},0}\circ
h_{0}(\mathfrak{I}_{{\gamma_0},0}^{-1}u)\quad\forall u\in{\bf H}^0(B)(\epsilon)_{\gamma_0},\\
&&{\cal L}^\circ_{\triangle \gamma_0}(u)={\cal
L}^\circ_{{\gamma_0},0}(\mathfrak{I}_{{\gamma_0},0}^{-1}u)\quad\forall
u\in {\bf H}^0(B)(\epsilon)_{\gamma_0}.\nonumber
\end{eqnarray}
For each $\gamma\in {\cal O}$ define the map
$\mathscr{J}_{\gamma}:\overline{{\bf
 H}^0(B)(\epsilon)_{\gamma}}\times \left({\bf
 H}^-(B)(\epsilon)_{\gamma}\oplus
 {\bf  H}^+(B)(\epsilon)_{\gamma}\right)\to\R$  by
\begin{equation}\label{e:3.34}
\mathscr{J}_{\gamma}(u, v)={\cal F}_{\gamma}(u+
\mathfrak{h}_{\gamma}(u)+ v)-{\cal F}_{\gamma}(u+
\mathfrak{h}_{\gamma}(u)).
\end{equation}
As in (\ref{e:2.28})-(\ref{e:2.31}) we have
\begin{eqnarray}
&&{\cal F}_{\gamma_0}\circ \mathfrak{I}_{{\gamma_0},0}={\cal
L}_{{\gamma_0},0}\quad\hbox{on}\;
 {\bf B}_{2r}(H_{{\gamma_0},0}),\nonumber\\
&&A_{\gamma_0}(\mathfrak{I}_{{\gamma_0},0}\xi)=\mathfrak{I}_{{\gamma_0},0}
\mathbb{A}_{\gamma_0,0}(\xi),\quad\forall\xi\in {\bf
B}_{2r}(H_{{\gamma_0},0})\cap X_{{\gamma_0},0},\label{e:3.35}\\
&&dA_{\gamma_0}(0)\circ \mathfrak{I}_{{\gamma_0},0}=\mathfrak{I}_{{\gamma_0},0}\circ d\mathbb{A}_{\gamma_0,0}(0),\label{e:3.36}\\
&&\mathfrak{I}_{{\gamma_0},0}\circ
\mathbb{P}_{{\gamma_0},0}^\star={\bf P}^\star_{\gamma_0}\circ
\mathfrak{I}_{{\gamma_0},0},\quad \star=-,0,+.\label{e:3.37}
\end{eqnarray}
It follows from these, (\ref{e:3.28}) and (\ref{e:3.33}) that
$\mathscr{J}_{\gamma_0}(u, v)={\bf
J}_{0}(\mathfrak{I}_{{\gamma_0},0}u, \mathfrak{I}_{{\gamma_0},0}v)$.
Hence $\mathscr{J}_{\gamma_0}$ is $C^{2-0}$ and
(\ref{e:3.29})-(\ref{e:3.30}) imply
\begin{eqnarray*}
&&D_2\mathscr{J}_{\gamma_0}(u, v)[w] =\langle ({\bf
P}^+_{\gamma_0}+{\bf P}^-_{\gamma_0})\nabla{\cal F}_{\gamma_0}(u+
\mathfrak{h}_{\gamma_0}(u)+ v), w\rangle_{1},
\\
&& \mathscr{J}_{\gamma_0}(u, 0)=0\quad\hbox{and}\quad
D_2\mathscr{J}_{\gamma_0}(u, 0)[w]=0
\end{eqnarray*}
for every $(u, v)\in \overline{{\bf
 H}^0(B)(\epsilon)_{\gamma_0}}\times\left({\bf
 H}^-(B)(\epsilon)_{\gamma_0}\oplus
 {\bf  H}^+(B)(\epsilon)_{\gamma_0}\right)$ and $w\in {\bf  H}^-(B)_{\gamma_0}+  {\bf  H}^+(B)_{\gamma_0}$.
These and Proposition~\ref{prop:3.7} immediately lead to

\begin{proposition}\label{prop:3.8}
\begin{enumerate}
\item[{\bf (i)}]  For any $(u, v, w_i)\in \overline{{\bf
 H}^0(B)(\epsilon)_{\gamma_0}}\times {\bf
 H}^+(B)(\epsilon)_{\gamma_0}\times {\bf
 H}^-(B)_{\gamma_0}$, $i=1,2$,
 $$
 (D_2\mathscr{J}_{\gamma_0}(u,
v+ w_2)-D_2\mathscr{J}_{\gamma_0}(u, v+
w_1))[w_2-w_1]\le-\frac{a_{\gamma_0}}{2} \|w_2-w_1\|_{1}.
$$
\item[{\bf (ii)}]  For any $(u, v, w)\in
 \overline{{\bf
 H}^0(B)(\epsilon)_{\gamma_0}}\times {\bf
 H}^+(B)(\epsilon)_{\gamma_0}\times {\bf
 H}^-(B)(\epsilon)_{\gamma_0}$,
$$
D_2\mathscr{J}_{\gamma_0}(u, v+
w)[v-w]\ge\frac{a_{\gamma_0}'}{2}\|v\|^2_{1}+
\frac{a_{\gamma_0}}{2}\|w\|^2_{1};
$$
in particular
 $D_2\mathscr{J}_{\gamma_0}(u, v)[v]\ge\frac{a_{\gamma_0}'}{2}\|v\|^2_{1}$ for any $(u, v)\in
\overline{{\bf
 H}^0(B)(\epsilon)_{\gamma_0}}\times {\bf
 H}^+(B)(\epsilon)_{\gamma_0}$.
\end{enumerate}
\end{proposition}

 By Proposition~\ref{prop:3.8} the functional
$\mathscr{J}_{\gamma_0}$  satisfies the conditions of Theorem~A.1 in
\cite{Lu2}. Then as in the proof of \cite[Lemma 3.6]{Lu2} we obtain
a $\R_{\gamma_0}$-invariant origin-preserving homeomorphism
$\Upsilon_{\gamma_0}$ from $N{\cal O}(\epsilon)_{\gamma_0}$ onto a
neighborhood of $0_{\gamma_0}\in N{\cal O}_{\gamma_0}$ and satisfies
$$
{\cal F}_{\gamma_0}\circ\Upsilon_{\gamma_0}(u)=\|{\bf
P}^+_{\gamma_0}u\|^2_1-\|{\bf P}^-_{\gamma_0}u\|^2_1+ {\cal
L}^\circ_{\triangle \gamma_0}({\bf P}^0_{\gamma_0}u),\quad\forall
u\in N{\cal O}(\epsilon)_{\gamma_0}.
$$
Since each $\gamma\in {\cal O}$  may be written as
$\gamma=s\cdot\gamma_0$ with some $s\in\R$, we have
$$
((-s)\cdot
u,(-s)\cdot v)\in\overline{{\bf
 H}^0(B)(\epsilon)_{\gamma_0}}\times \left({\bf
 H}^-(B)(\epsilon)_{\gamma_0}\oplus
 {\bf  H}^+(B)(\epsilon)_{\gamma_0}\right)
 $$
for $(u,v)\in\overline{{\bf
 H}^0(B)(\epsilon)_{\gamma}}\times \left({\bf
 H}^-(B)(\epsilon)_{\gamma}\oplus
 {\bf  H}^+(B)(\epsilon)_{\gamma}\right)$.
 From (\ref{e:3.34}) we derive
$\mathscr{J}_{\gamma}(u,v)=\mathscr{J}_{\gamma_0}((-s)\cdot
u,(-s)\cdot v)$ because the map $\mathfrak{h}$ in (\ref{e:1.24}) is
 $\R$-equivariant.
Define $\Upsilon:N{\cal O}(\epsilon)\to N{\cal O}$ by
$\Upsilon_{s\cdot\gamma_0}(u)=s\cdot\Upsilon_{\gamma_0}((-s)\cdot u)$ for any $u\in
N{\cal O}(\epsilon)_{s\cdot \gamma_0}$ and any $s\in \R$.
It is $\R$-equivariant fiber-preserving homeomorphism onto an open
neighborhood of ${\cal O}$ in $N{\cal O}$ and satisfies
${\cal F}\circ\Upsilon(u)=\|{\bf P}^+u\|^2_1-\|{\bf P}^-u\|^2_1+
{\cal L}^\circ_\triangle({\bf P}^0u)$ for all $u\in N{\cal
O}(\epsilon)$.
\hfill$\Box$\vspace{2mm}

\section{Proofs of Corollary~\ref{cor:1.11} and Theorem~\ref{th:1.12}}\label{sec:4}
\setcounter{equation}{0}

\noindent{\bf Proof of Corollary~\ref{cor:1.11}}.\quad
We only need to check that the conditions (i)-(iv) of Theorem~\ref{th:1.9}
are satisfied with $\tilde\Lambda=\Lambda(S, \mathbb{I}_g)$ and
  $m=1$.

Observe that $\tilde{\cal X}:={\cal X}(M, \mathbb{I}_g)\cap\tilde\Lambda={\cal X}(S, \mathbb{I}_g)$. So
 it is a Banach submanifold of ${\cal X}(M, \mathbb{I}_g)$.
Since $S$ is a totally geodesic submanifold of $(M, g)$ it is easily seen that the condition (ii) holds.

Let us  prove that the condition (iii) is satisfied, i.e.,
 \begin{equation}\label{e:4.1}
 \nabla({\cal L}\circ\widetilde{\rm EXP}_{\gamma})(\xi)=\nabla{\cal F}_\gamma(\xi)\quad\forall (\gamma, \xi)\in \tilde{N}{\cal O}(\varepsilon).
 \end{equation}
Define a map $\tilde{\mathbb{I}}_g:\Lambda(M, \mathbb{I}_g)\to \Lambda(M, \mathbb{I}_g)$ by
$\tilde{\mathbb{I}}_g(\beta)(t)=\mathbb{I}_g(\beta(t)),\;\forall t\in\R$, for $\beta\in
\Lambda(M, \mathbb{I}_g)$. We claim that it is a Riemannian isometry on $\Lambda(M, \mathbb{I}_g)$,
that is,
\begin{equation}\label{e:4.2}
\langle d\tilde{\mathbb{I}}_g(\beta)[\xi], d\tilde{\mathbb{I}}_g(\beta)[\eta]\rangle_1=
\langle \xi, \eta\rangle_1,\quad\forall\beta\in\Lambda(M, \mathbb{I}_g),\;\forall \xi,\eta\in
T_\beta\Lambda(M, \mathbb{I}).
\end{equation}
Let $\alpha=\tilde{\mathbb{I}}_g(\beta)$. Note that $d\tilde{\mathbb{I}}_g(\beta)[\xi]\in T_\alpha
\Lambda(M, \mathbb{I}_g)$ is given by
$$
d\tilde{\mathbb{I}}_g(\beta)[\xi](t)=d{\mathbb{I}}_g(\beta(t))[\xi(t)]\;\forall t\in\R.
$$
Since $\mathbb{I}_g$ is an isometry on $(M, g)$ we have $d{\mathbb{I}}_g(x)\bigl(\nabla^g_vX)=
\nabla^g_{d{\mathbb{I}}_g(x)[v]}(\mathbb{I}_{g\ast}X)\bigr)$ for any $v\in T_xM$ and any
differentiable vector field $X$ near $x$. It follows that
\begin{eqnarray*}
&&\langle d{\mathbb{I}}_g(\beta(t))[\xi(t)], d{\mathbb{I}}_g(\beta(t))[\eta(t)]\rangle=
\langle \xi(t),\eta(t)\rangle\quad\forall t\in\R,\\
&&\langle (\nabla^g_{\dot\alpha}\mathbb{I}_{g\ast}\xi)(t)), (\nabla^g_{\dot\alpha}\mathbb{I}_{g\ast}\eta)(t))\rangle=
\langle (\nabla^g_{\dot\beta}\xi)(t),  (\nabla^g_{\dot\beta}\xi)(t)\rangle\quad\forall t\in\R,
\end{eqnarray*}
which imply (\ref{e:4.2}).

Let $G$ be the subgroup of $I(M,g)$ generated by $\mathbb{I}_g^k$, i.e.,
$G=\{\mathbb{I}_g^{lk}\,|\, \pm l\in\N\cup\{0\}\}$
and let $\tilde{G}=\{\tilde{\mathbb{I}}_g^{lk}\,|\, \pm l\in\N\cup\{0\}\}$.
The latter is a subgroup of the isometry group of $\Lambda(M, \mathbb{I}_g)$ and
${\rm Fix}(\tilde{\mathbb{I}}^k)$ is exactly the set of points fixed under the
action of $\tilde{G}$.

Observe that $\Lambda(S, \mathbb{I}_g)$ is a connected component of
${\rm Fix}(\tilde{\mathbb{I}}_g^k)$ containing $\gamma_0$.
(In fact, it is a totally geodesic submanifold  of $\Lambda(M, \mathbb{I}_g)$
by \cite[Cor.5.4.2]{PaT}).
For any $\gamma\in \Lambda(S, \mathbb{I}_g)$ we have
 $$T_\gamma\Lambda(S, \mathbb{I}_g)=\{\xi\in T_\gamma\Lambda(M, \mathbb{I}_g)\,|\,
d\tilde{\mathbb{I}}^k(\gamma)[\xi]=\xi\}.
$$
By (\ref{e:4.2}) $d\tilde{\mathbb{I}}_g^k(\gamma)$ restricts to a Hilbert space isometry from
$N{\cal O}_\gamma$ onto itself, denoted by $d\tilde{\mathbb{I}}_g^k(\gamma)'$, whose fixed point set
is $\tilde{N}{\cal O}_\gamma$. Moreover, since
\begin{eqnarray*}
{\rm EXP}(\gamma,d\tilde{\mathbb{I}}_g^k(\gamma)'[\xi])(t)&=&
\exp_{\gamma(t)}\Bigl(d{\mathbb{I}}_g^k(\gamma(t))[\xi(t)]\Bigr)\\
&=&\exp_{\mathbb{I}_g^k\gamma(t)}\Bigl(d{\mathbb{I}}_g^k(\gamma(t))[\xi(t)]\Bigr)
=\mathbb{I}_g^k\Bigl(\exp_{\gamma(t)}\xi(t)\Bigr)\;\forall
t\in\R
\end{eqnarray*}
for any $\xi\in N{\cal O}_\gamma$, we deduce that ${\cal F}_\gamma$ is
$d\tilde{\mathbb{I}}_g^k(\gamma)'$-invariant.
By the argument of \cite[Section~2]{Pa79}, for any $\xi\in N{\cal O}(\varepsilon)_\gamma$
and $l\in\N\cup\{0\}$ we have
$$(d\tilde{\mathbb{I}}_g^k(\gamma)')^l\bigl(\nabla{\cal F}_\gamma(\xi)\bigr)=
\nabla{\cal F}_\gamma((d\tilde{\mathbb{I}}_g^k(\gamma)')^l\xi).
$$
In particular, for any $\xi\in \tilde{N}{\cal O}_\gamma$ we have
$d\tilde{\mathbb{I}}_g^k(\gamma)'\bigl(\nabla{\cal F}_\gamma(\xi)\bigr)=
\nabla{\cal F}_\gamma(\xi)$ and hence
$\nabla{\cal F}_{\gamma}(\xi)\in\tilde{N}{\cal O}_\gamma$.
This implies (\ref{e:4.1}), because ${\cal L}\circ\widetilde{\rm EXP}_{\gamma}(\xi)=
{\cal F}_\gamma(\xi),\;\forall (\gamma, \xi)\in \tilde{N}{\cal O}(\varepsilon)$.

Finally, we check that condition (iv) is satisfied. Observe that
\begin{eqnarray*}
&&{\bf H}^0\bigl(d^2({\cal L}|_{{\cal X}(M, \mathbb{I}_g)})(\gamma_0)\bigr)=
{\bf H}^0\bigl(d^2\bigl({\cal L}\circ{\rm EXP}|_{T_{\gamma_0}{\cal X}(M,\mathbb{I}_g)(\varepsilon)}\bigr)(0)\bigr)\\
\bigl(\hbox{resp.}&&
{\bf H}^0\bigl(d^2({\cal L}|_{{\cal X}(S, \mathbb{I}_g)})(\gamma_0)\bigr)=
{\bf H}^0\bigl(d^2\bigl({\cal L}\circ{\rm EXP}|_{T_{\gamma_0}{\cal X}(S,\mathbb{I}_g)(\varepsilon)}\bigr)(0)\bigr)\;\bigr)
\end{eqnarray*}
consists of $\mathbb{I}_g$-invariant  Jacobi fields on $(M, F)$ (resp. $(S, F|_S)$) along $\gamma_0$.
For\linebreak $\xi\in {\bf H}^0\bigl(d^2({\cal L}|_{{\cal X}(S, \mathbb{I}_g)})(\gamma_0)\bigr)$ there is a geodesic variation $H$ of $\gamma_0$ on $(S,F|_S)$ whose variation field is equal to $\xi$. Since
 $\mathbb{I}_g\in I(M,F)$ implies that
$S$ is also a totally geodesic submanifold of $(M,F)$,   $H$ is also a geodesic variation  of $\gamma_0$ on $(M,F)$ and hence  $\xi$ is a Jacobi field  on $(M, F)$ along $\gamma_0$, that is,
$\xi\in{\bf H}^0\bigl(d^2({\cal L}|_{{\cal X}(M, \mathbb{I}_g)})(\gamma_0)\bigr)$. It follows that
${\bf H}^0\bigl(d^2({\cal L}|_{{\cal X}(S, \mathbb{I}_g)})(\gamma_0)\bigr)\subset
{\bf H}^0\bigl(d^2({\cal L}|_{{\cal X}(M, \mathbb{I}_g)})(\gamma_0)\bigr)$ and so
${\bf H}^0\bigl(d^2({\cal L}|_{{\cal X}(S, \mathbb{I}_g)})(\gamma_0)\bigr)=
{\bf H}^0\bigl(d^2({\cal L}|_{{\cal X}(M, \mathbb{I}_g)})(\gamma_0)\bigr)$ because\linebreak
 $m^0({\cal L}, {\cal O})=m^0({\cal L}|_{\Lambda(S, \mathbb{I}_g)}, {\cal O})$.
 \hfill$\Box$\vspace{2mm}

\noindent{\bf Proof of Theorem~\ref{th:1.12}}.\quad
Note that the iteration map $\varphi_m$ gives not only  Riemannian isometry (up to the factor $m^2$)
embedding from $\Lambda(M,\mathbb{I}_g)$ into $\Lambda(M,\mathbb{I}^m_g)$, but also a Banach manifold
embedding from ${\cal X}(M,\mathbb{I}_g)$ into ${\cal X}(M,\mathbb{I}^m_g)$.
So $\tilde\Lambda:=\varphi_m\bigl(\Lambda(M,\mathbb{I}_g)\bigr)$ (resp.
$\tilde{\cal X}=\varphi_m\bigl({\cal X}(M,\mathbb{I}_g)\bigr)$) is a Hilbert-Riemannian (resp. Banach)
submanifold of $\Lambda(M,\mathbb{I}^m_g)$ (resp. ${\cal X}(M,\mathbb{I}^m_g)$).
Moreover ${\cal L}\circ\varphi_m=m^2{\cal L}$
because $F$ is $\mathbb{I}_g$-invariant as we assumed below Theorem~\ref{th:1.5}. It follows that the diffeomorphism $\varphi_m:\Lambda(M,\mathbb{I}_g)
\to\tilde\Lambda$
induces isomorphisms
\begin{eqnarray}
&&(\varphi_m)_\ast: \mathscr{H}^0_\ast({\cal L}, \gamma;\K)\to
\hat{\mathscr{H}}^0_\ast({\cal L}|_{\tilde\Lambda}, \gamma^m;\K),
\quad \forall\gamma\in{\cal O},\label{e:4.3}\\
&&(\varphi_m)_\ast: C_\ast({\cal L}, {\cal O};\K)\to
C_\ast({\cal L}|_{\tilde\Lambda}, {\cal O}^m;\K).\label{e:4.4}
\end{eqnarray}

It remains to check that the conditions (ii)-(iv) in Theorem~\ref{th:1.9} are satisfied.
From the definition of ${\rm EXP}^m$ in (\ref{e:1.34}) it is easily seen that
(ii) holds. The other two ones can be proved in similar way to that of (6.13) and (6.15) in \cite{Lu3}.
Let us see them in detail.

Since $D\varphi_m(x):T_x\Lambda(M, \mathbb{I}_g)\to T_{x^m}\Lambda(M,\mathbb{I}^m_g)$ is given by
$D\varphi_m(x)[v]=v^m$, where $v^m(t)=v(mt)\;\forall t\in\R$, we have
$(D\varphi_m(x)[\xi], D\varphi_m(x)[\eta])_m=m^2\langle\xi,\eta\rangle_1$
for any $\xi,\eta\in T_x\Lambda(M, \mathbb{I}_g)$, and so
 $D\varphi_m(x)\bigl(N{\cal O}(\varepsilon/m)_x\bigr)=\tilde{N}{\cal O}^m(\varepsilon)_{x^m}\;\forall x\in{\cal O}$.
    Moreover ${\cal L}\circ\widetilde{\rm EXP}_{x^m}(v^m)= m^2{\cal
F}_{x}(v)\;\forall v\in N{\cal O}(\varepsilon)_x$. We can deduce that
\begin{eqnarray}\label{e:4.5}
\hat\nabla({\cal L}\circ\widetilde{\rm EXP}_{x^m})(v^m)=m\nabla{\cal
F}_{x}(v)\quad \forall (x,v)\in N{\cal O}(\varepsilon/m).
\end{eqnarray}

Consider the isometric action  on $\hat{N}{\cal
O}^m(\varepsilon)_{x^m}$ of group $\Z_m=\{e^{2\pi
ip/m}\,|\,p=0,\cdots,m-1\}$ given by
$(e^{2\pi ip/m}\cdot u)(t):=d\mathbb{I}^{-p}_g\bigl(x^m(t+ p/m)\bigr)[u(t+ p/m)]$ for all
$t\in\R$ and $p=0,\cdots,m-1$.
Indeed, for $u\in\hat{N}{\cal
O}^m(\varepsilon)_{x^m}$ we have $u(t)\in T_{x^m(t)}M=T_{x(mt)}M$. Since $\mathbb{I}_g(x(t))=x(t+1)\;\forall t$ it follows that
\begin{eqnarray*}
d\mathbb{I}^{-p}_g\bigl(x^m(t+ p/m)\bigr)[u(t+ p/m)]&=&d\mathbb{I}^{-p}_g\bigl(x(mt+p)\bigr)[u(t+ p/m)]\\
&=&d\mathbb{I}^{-p}_g\bigl(\mathbb{I}^p_g(x(mt))\bigr)[u(t+ p/m)]\in T_{x^m(t)}M.
\end{eqnarray*}
 Moreover a straightforward  computation shows
\begin{eqnarray*}
&&\tilde{\mathbb{I}}^m_g(e^{2\pi ip/m}\cdot u)(t):=d{\mathbb{I}}^m_g(x^m(t))[(e^{2\pi ip/m}\cdot u)(t)]\\
&&=d{\mathbb{I}}^m_g(x^m(t))\bigl[d\mathbb{I}^{-p}_g\bigl(x^m(t+ p/m)\bigr)[u(t+ p/m)]\bigr]\\
&&=d{\mathbb{I}}^m_g(x^m(t))\circ d\mathbb{I}^{-p}_g\bigl(x^m(t+ p/m)\bigr)[u(t+ p/m)]\\
&&=d{\mathbb{I}}^{-p}_g(x^m(t+1+ p/m))\circ d\mathbb{I}^{m}_g\bigl(x^m(t+ p/m)\bigr)[u(t+ p/m)]\\
&&=d{\mathbb{I}}^{-p}_g(x^m(t+1+ p/m))[u(t+ 1+ p/m)]\\
&&=(e^{2\pi ip/m}\cdot u)(t+1)\;\forall t\in\R.
\end{eqnarray*}
That is, $e^{2\pi ip/m}\cdot u\in \hat{N}{\cal
O}^m(\varepsilon)_{x^m}$. If $e^{2\pi ip/m}\cdot u=u$, $p=0,\cdots,m-1$, then
$$
d\mathbb{I}^{-p}_g\bigl(x^m(t+ p/m)\bigr)[u(t+ p/m)]=u(t)\quad\forall
t\in\R,\;p=0,\cdots,m-1,
$$
and in particular we have  $d\mathbb{I}_g\bigl(x^m(t)\bigr)[u(t)]=u(t+ 1/m)$, which
implies $v(t):=u(t/m)$ to give an element of $N{\cal O}(\varepsilon/m)_x$  with $v^m=u$.
Hence  $\tilde{N}{\cal O}^m(\varepsilon)_{x^m}$ is exactly the fixed
point set of the above $\Z_m$-action. It is not hard to see that the
functional $\hat{\cal F}_{x^m}$ is invariant under this
$\Z_m$-action.  By the argument on the page 23 of \cite{Pa79} we obtain
$$
e^{2\pi ip/m}\cdot \hat\nabla\hat{\cal
F}_{x^m}(v^m)=\hat\nabla\hat{\cal
F}_{x^m}(v^m)\quad\forall (x,v)\in N{\cal O}(\varepsilon/m),\;p=0,\cdots,m-1.
$$
This implies $\hat\nabla\hat{\cal
F}_{x^m}(v^m)\in\tilde{N}{\cal O}^m_{x^m}$, and so
\begin{eqnarray}\label{e:4.6}
\hat\nabla\hat{\cal
F}_{x^m}(v^m)=\hat\nabla({\cal L}\circ\widetilde{\rm EXP}_{x^m})(v^m),\quad \forall (x,v)\in N{\cal O}(\varepsilon/m).
\end{eqnarray}
because of (\ref{e:4.5}) and ${\cal L}\circ\widetilde{\rm EXP}_{x^m}=\hat{\cal
F}_{x^m}|_{\tilde{N}{\cal O}^m_{x^m}}$. The condition (iii) is satisfied.

Finally, we prove (iv) under the assumption $m^0({\cal L}, {\cal O}^m)=
m^0({\cal L}, {\cal O})$. Note that
 $$
 \hat{\cal F}^X_{x^m}(s e^{2\pi ip/m}\cdot v)=\hat{\cal F}^X_{x^m}(s v)
 $$
 for a given $v\in X\hat{N}{\cal O}^m_{x^m}$ and $s\in\R$ with $|s|$ small enough.
Taking the second derivative of both sides at $s=0$ we obtain
 $$
 d^2\hat{\cal F}^X_{x^m}(0)[e^{2\pi ip/m}\cdot v, e^{2\pi ip/m}\cdot v]
=d^2\hat{\cal F}^X_{x^m}(0)[v, v]
$$
  and hence
 $d^2\hat{\cal F}^X_{x^m}(0)[e^{2\pi ip/m}\cdot u, e^{2\pi ip/m}\cdot v]
 =d^2\hat{\cal F}^X_{x^m}(0)[u, v]\;\forall u,v\in X\hat{N}{\cal O}^m_{x^m}$, which
 leads, by the density of $X\hat{N}{\cal O}^m_{x^m}$ in $\hat{N}{\cal O}^m_{x^m}$,   to
 \begin{eqnarray*}
  (\hat{B}_{x^m}e^{2\pi ip/m}\cdot u, e^{2\pi ip/m}\cdot v)_m
 =(\hat{B}_{x^m}u, v)_m\quad\forall u,v\in \hat{N}{\cal O}^m_{x^m}.
 \end{eqnarray*}
  It follows that $\hat{B}_{x^m}e^{2\pi ip/m}\cdot u=e^{2\pi ip/m}\cdot(\hat{B}_{x^m} u)$, $p=0,1,\cdots,m$,
 and in particular,
  $$
 \hat{B}_{x^m}u=e^{2\pi ip/m}\cdot(\hat{B}_{x^m} u)\quad\forall u\in \tilde{N}{\cal O}^m_{x^m},\;   p=0,1,\cdots,m.
 $$
These show that  $\hat{B}_{x^m}u\in\tilde{N}{\cal O}^m_{x^m}$, that is,
\begin{eqnarray}\label{e:4.7}
\hat{B}_{x^m}(\tilde{N}{\cal O}^m_{x^m})\subset \tilde{N}{\cal O}^m_{x^m}.
\end{eqnarray}

 On the other hand,  (\ref{e:4.5}) and (\ref{e:4.6}) imply   $\hat A_{x^m}(D\varphi_m(x)[v])=D{\varphi}_{m}(x)[
A_{x}(v)]$ for any $(x,v)\in  N{\cal
O}(\varepsilon/m)\cap XN{\cal O}$,
 and hence
 \begin{eqnarray*}
 (\hat A_{x^m}(D\varphi_m(x)[v+su]), D\varphi_m(x)u)_m&=&(D{\varphi}_{m}(x)[
A_{x}(v+su)], D\varphi_m(x)u)_m\\
&=&m^2\langle
A_{x}(v+su), u\rangle_1
\end{eqnarray*}
 for a given $u\in XN{\cal O}_x$ and $s\in\R$ with $|s|$ small enough.
  Taking derivatives at $s=0$ and $v=0$ we obtain $d^2\hat{\cal F}^X_{x^m}(0)[v^m,u^m]=m^2d^2{\cal F}^X_x(0)[v,u]$,
and then as above arrive at
$$
(\hat{B}_{x^m}v^m, u^m)_m=m^2\langle B_xv,u\rangle_1\quad\forall u,v\in N{\cal O}_x.
$$
This and (\ref{e:4.7}) imply that $D\varphi_m(x)({\rm Ker}(B_x))\subset{\rm Ker}(\hat{B}_{x^m})\;\forall
x\in{\cal O}$, and thus
\begin{eqnarray}\label{e:4.8}
D\varphi_m(\gamma_0)\bigl({\bf H}^0\bigl(d^2\hat{\cal F}^X_{\gamma_0}(0)\bigr)\bigr)\subset
 {\bf H}^0\bigl(d^2\hat{\cal F}^X_{\gamma_0^m}(0)\bigr).
\end{eqnarray}
Since $\dim{\bf H}^0\bigl(d^2\hat{\cal F}^X_{\gamma_0^m}(0)\bigr)=m^0({\cal L}, {\cal O}^m)+1=
m^0({\cal L}, {\cal O})+1=\dim{\bf H}^0\bigl(d^2\hat{\cal F}^X_{\gamma_0}(0)\bigr)$ is finite and
$D\varphi_m(\gamma_0)$ is a linear injection, from (\ref{e:4.8}) we derive
$$
{\bf H}^0\bigl(d^2\hat{\cal F}^X_{\gamma_0^m}(0)\bigr)=D\varphi_m(\gamma_0)\bigl({\bf H}^0\bigl(d^2\hat{\cal F}^X_{\gamma_0}(0)\bigr)\bigr)\subset \tilde{N}{\cal O}^m_{\gamma_0^m}\subset T_{\gamma_0^m}\tilde\Lambda.
$$
The condition (iv) is satisfied.

Observe that $\varphi_m:\Lambda(M,\mathbb{I}_g)\to\Lambda(M,\mathbb{I}^m_g)$
is the composition of the diffeomorphism $\varphi_m:\Lambda(M,\mathbb{I}_g)
\to\tilde\Lambda$ and the inclusion $\tilde\Lambda\hookrightarrow \Lambda(M,\mathbb{I}^m_g)$.
The expected conclusions follow from (\ref{e:4.3})-(\ref{e:4.4}) and Theorem~\ref{th:1.9} immediately.
 \hfill$\Box$\vspace{2mm}

\section{Propositions and lemmas}\label{sec:5}
\setcounter{equation}{0}

\begin{proposition}\label{prop:5.1}
For $k\in\N$, $k>1$, let $({\cal H}, (\!(\cdot, \cdot)\!))$ be a $C^k$ Hilbert-Riemannian
manifold modeled on $H$, and let
 $({\cal X}, \|\cdot\|_{\cal X}$) be a $C^k$  Banach-Finsler manifold
modeled on  $X$. Suppose also that ${\cal X}\subset{\cal H}$ and the
inclusion ${\cal X}\hookrightarrow {\cal H}$ is $C^k$-smooth. Then
for a compact $C^k$ submanifold ${\cal O}$ of ${\cal H}$, which is
also contained in ${\cal X}$, we have\\
\noindent{\rm (i)} ${\cal O}$ is also  a $C^k$-smooth submanifold of
${\cal X}$, and ${\cal H}$ and ${\cal X}$ induce an equivalent
$C^k$-smooth manifold structure (including topology) on
$\mathcal{O}$;\\
\noindent{\rm (ii)} For the orthogonal normal bundle $N\mathcal{O}$
of ${\cal O}$ in ${\cal H}$ with respect to the metric $(\!(\cdot,
\cdot)\!)$, $XN\mathcal{O}:=(T_\mathcal{O}{\cal X})\cap N\mathcal{
O}$ is a $C^{k-1}$ subbundle of $T_\mathcal{O}{\cal X}$;\\
\noindent{\rm (iii)} If ${\cal F}\to{\cal O}$ is a $C^{k-1}$ subbundle
of $T_{\cal O}{\cal H}$ contained in $T_{\cal O}{\cal X}$ then it is
also a $C^{k-1}$ subbundle of $T_{\cal O}{\cal X}$; furthermore for the
orthogonal complementary bundle ${\cal F}^\bot$ of ${\cal F}$ in
$T_{\cal O}{\cal H}$ the bundle $X{\cal F}^\bot:={\cal F}^\bot\cap
T_{\cal O}{\cal X}$ is also a $C^1$ subbundle of $T_{\cal O}{\cal
X}$
\end{proposition}

\noindent{\bf Proof}. (i) Since  the inclusion ${\cal
X}\hookrightarrow {\cal H}$ is $C^k$-smooth, the
identity map  ${\cal X}\supset{\cal O}\ni x\mapsto x\in {\cal O}\subset
{\cal H}$ is an injective immersion with respect to the induced
submanifold structures from ${\cal H}$ and ${\cal X}$ respectively.
So it is a diffeomorphism by the compactness of ${\cal O}$ and
Proposition 3.3.4 in \cite[Page 149]{MaOu}. The first claim is
proved.

(ii)  Let $\Pi$ be the orthogonal bundle projection from
$T_\mathcal{ O}{\cal H}$ onto $T{\cal O}$ with respect to
$(\!(\cdot,\cdot)\!)$. It is a $C^{k-1}$ $\mathcal{O}$-bundle
morphism. By the  assumption
 the inclusion $T_\mathcal{ O}{\cal X}\hookrightarrow T_\mathcal{ O}{\cal H}$ is
a $C^{k-1}$-smooth $\mathcal{ O}$-bundle morphism. This and (i) assure
that the composition $\tilde\Pi$ of
$T_\mathcal{ O}{\cal X}\hookrightarrow T_\mathcal{O}{\cal H}$ and
$\Pi:T_\mathcal{O}{\cal H}\to T{\cal O}$
 is a $C^{k-1}$ $\mathcal{O}$-vector bundle morphism with kernel $XN\mathcal{ O}$. For
any $x\in\mathcal{ O}$  assume that the sequence $(v_n)\subset XN{\cal O}_x=T_x{\cal X}\cap
N\mathcal{ O}_x$ converges to $v\in T_x{\cal X}$ in $T_x{\cal X}$, i.e.,
 $\|v_n-v\|_{{\cal X},x}\to 0$. Then  $v\in N\mathcal{ O}_x$ because the inclusion
$T_x{\cal X}\hookrightarrow T_x{\cal H}$ is continuous and
\begin{eqnarray*}
|(\!(v, u)\!)_x|&=&|(\!(v_n, u)\!)-(\!(v, u)\!)_x|\le
\|v_n-v\|_x\cdot\|u\|_x\to 0\quad\forall u\in T_x\mathcal{ O}.
\end{eqnarray*}
Thus $T_x{\cal X}\cap N\mathcal{ O}_x$ is a closed
subspace of $T_x{\cal X}$. Observe that it also splits $T_x{\cal X}$
because of the topological direct sum decomposition $T_x{\cal
X}=T_x\mathcal{O}\oplus (T_x{\cal X}\cap N\mathcal{ O}_x)$.
That is, for each $x\in{\cal O}$ the continuous linear map $\tilde\Pi_x:T_x{\cal X}\to T_x{\cal O}$
is surjective and has a kernel that splits.  From
Proposition~6 of \cite[Page 52]{La} we derive
that the sequence $T_\mathcal{O}{\cal X} \stackrel{\tilde\Pi}{\to}T{\cal O}\to 0$ is exact.
Then we obtain that $XN\mathcal{ O}$ is a $C^{k-1}$
subbundle of $T_\mathcal{ O}{\cal X}$ by the argument above Proposition~5 in
\cite[Page 51]{La}.

(iii) The first claim may be obtained as in the proof of (ii) since
${\cal F}$ (as a subset of $T_{\cal O}{\cal X}$) is the kernel of
the composition of $C^1$ bundle morphisms $T_{\cal O}{\cal
X}\hookrightarrow T_{\cal O}{\cal H}$ and the orthogonal bundle
projection $T_{\cal O}{\cal H}\to{\cal F}^\bot$. As to the second one
note that $X{\cal F}^\bot$ is the kernel of the restriction to
$T_{\cal O}{\cal X}$ of the orthogonal bundle projection $T_{\cal
O}{\cal H}\to{\cal F}$. \hfill$\Box$\vspace{2mm}

\noindent{\bf Proof of Lemma~\ref{lem:2.2}}. Let ${\bf E}$ be as in
(\ref{e:2.7}). Since $d{\bf E}(t,\xi)[s,\eta]=(s, \eta(t),
\dot\eta(t))$ for $s\in\R$ and $\eta\in C^{k+1}_{V}(I, \R^n)=T_\xi
C^{k+1}_{V}(I, {\bf B}_{2\rho}(\R^n))$, it is easily checked that
${\bf E}$ is a submersion. In particular, it is transversal to the
submanifold ${\cal Z}\subset I\times{\bf
B}_{2\rho}(\R^n)\times\R^n$, and hence for each $(t,\xi)\in {\bf
E}^{-1}({\cal Z})$ the tangent space
$$
T_{(t,\xi)}{\bf E}^{-1}({\cal Z})=(d{\bf
E}(t,\xi))^{-1}(T_{\check{\xi}(t)}{\cal Z})
$$
has codimension $n$ in $\R\times T_\xi C^{k+1}_{V}(I, {\bf
B}_{2\rho}(\R^n))$. By Lemma~8.5.2 on the page 409 of \cite{MaOu}
the restriction ${\bf P}$ of the natural projection  ${\rm P}_2:
I\times C^{k+1}_{V}(I, {\bf B}_{2\rho}(\R^n))\to C^{k+1}_{V}(I, {\bf
B}_{2\rho}(\R^n))$ to ${\bf E}^{-1}({\cal Z})$ is a $C^k$ Fredholm
operator with ${\rm ind}({\bf P})=1-n$. Since $k\ge 2$ the Sard-Smale
theorem concludes that  all regular values of ${\bf P}$ forms a
residual subset
 $C^{k+1}_{V}(I, {\bf B}_{2\rho}(\R^n))^\circ_{\rm reg}$
in  $C^{k+1}_{V}(I, {\bf B}_{2\rho}(\R^n))$. Note that $\xi\in
C^{k+1}_{V}(I, {\bf B}_{2\rho}(\R^n))$ sits in $C^{k+1}_{V}(I, {\bf
B}_{2\rho}(\R^n))^\circ_{\rm reg}$ if and only if the map ${\bf
E}_\xi$ is transversal to the submanifold ${\cal Z}$, and in this
case ${\bf E}_\xi^{-1}({\cal Z})={\bf P}^{-1}(\xi)$ is a $C^k$
submanifold of dimension $1-n$. If $n>1$ we obtain that $\xi\in
C^{k+1}_{V}(I, {\bf B}_{2\rho}(\R^n))$ belongs to  $C^{k+1}_{V}(I,
{\bf B}_{2\rho}(\R^n))^\circ_{\rm reg}$ if and only if ${\bf
E}_\xi^{-1}({\cal Z})={\bf P}^{-1}(\xi)=\emptyset$, and so
$C^{k+1}_{V}(I, {\bf B}_{2\rho}(\R^n))^\circ_{\rm
reg}=C^{k+1}_{V}(I, {\bf B}_{2\rho}(\R^n))_{\rm reg}$. The first two
conclusions follow immediately.

Let  $\xi_i\in C^{k+1}_{V}(I, {\bf B}_{2\rho}(\R^n))_{\rm
reg}^\circ$, $i=1,2$. Denote by $\mathcal {P}_{k+1}(\xi_0,\xi_1)$
the set of all
 $C^{k+1}$ path $\mathbbm{p}:[0, 1]\to C^{k+1}_{V}(I,
{\bf B}_{2\rho}(\R^n))$ connecting $\xi_0$ to $\xi_1$. As above we
may prove that the $C^k$ map
$$
[0,1]\times I\times  \mathcal {P}_{k+1}(\xi_0,\xi_1)\to I\times{\bf
B}_{2\rho}(\R^n)\times\R^n,\; (s,t, \mathbbm{p})\mapsto {\bf E}(t,
\mathbbm{p}(s))
$$
is a submersion, and in particular transversal to the submanifold
${\cal Z}$. Hence we get a residual
 subset  $\mathcal
{P}_{k+1}(\xi_0,\xi_1)_{\rm reg}\subset \mathcal
{P}_{k+1}(\xi_0,\xi_1)$ such that for every  $\mathbbm{p}\in
\mathcal {P}_{k+1}(\xi_0,\xi_1)_{\rm reg}$  the set $\{(s,t)\in [0,
1]\times I\,|\, {\bf E}(t,\mathbbm{p}(s))\in{\cal Z}\}$ is a $C^k$
submanifold of $[0, 1]\times I$ of dimension $2-n$, and so empty for
$n>2$, and at most a finite  set for $n=2$. These lead to the third
conclusion. \hfill$\Box$\vspace{2mm}

We need the following special version of the omega lemma by
Piccione and Tausk  \cite[Theorem~4.3]{PiTa}.

\begin{lemma}\label{lem:5.2}
Let $N_1$ and $N_2$ be two manifolds of class $C^{m+l}$ ($m\ge 2$ and $l\ge 1$) and let
$T:N_1\to N_2$ be a $C^{m+l}$ map. Then for $I=[0,1]$ the maps
\begin{description}
\item[(i)] $H^1(I, N_1)\to H^1(I, N_1),\;\gamma\mapsto T\circ\gamma$,
is $C^{m-1}$ with $l=2$,
\item[(ii)] $C^1(I, N_1)\to C^1(I, N_1),\;\gamma\mapsto T\circ\gamma$,
is $C^{m-1}$ with $l=2$,
\item[(iii)] $C^0(I, N_1)\to C^0(I, N_1),\;\gamma\mapsto T\circ\gamma$,
is $C^{m-1}$ with $l=1$.
\end{description}
\end{lemma}

\noindent{\bf Proof of Lemma~\ref{lem:2.3}}. The first two claims and $C^1$ smoothness of $\tilde{A}$ can be proved as in the proof of \cite[Lemma 3.2]{Lu1}. The last one will be proved  as follows. Consider the smooth maps
\begin{eqnarray*}
&&\mathfrak{V}:{\bf B}_{2r}(X_V)\to C^0(I, I\times B_{2r}(\R^n)\times\R^n),\;\xi\mapsto\check{\xi},\\
&&\mathfrak{K}_j: C^0(I,\R^n)\to C^1(I,\R^n),\;\xi\mapsto \mathfrak{K}(\xi),\;j=1,2,3
\end{eqnarray*}
where $\mathfrak{K}_1(\xi)(t)=\int^t_0\xi(s)ds$, $\mathfrak{K}_2(\xi)(t)=\int^t_0e^s\xi(s)ds$, and
$\mathfrak{K}_3(\xi)(t)=\int^t_0e^{-2s}\xi(s)ds$. \linebreak By Lemma~\ref{lem:5.2}(iii) we have also
 the $C^{k-3}$ maps
\begin{eqnarray*}
&&\mathfrak{L}_q: C^0(I, I\times B_{2r}(\R^n)\times\R^n)\to C^0(I,\R^n),\;\eta\mapsto\partial_{q}\tilde{L}\circ\eta,\\
&&\mathfrak{L}_v: C^0(I, I\times B_{2r}(\R^n)\times\R^n)\to C^0(I,\R^n),\;\eta\mapsto\partial_{v}\tilde{L}\circ\eta
\end{eqnarray*}
induced by the $C^{k-1}$ maps $\partial_{q}\tilde{L},\partial_{v}\tilde{L}:I\times B_{2r}(\R^n)\times\R^n\to\R^n$,
respectively. Then   from (\ref{e:2.11}) and (\ref{e:2.12}) we deduce for any $\xi\in {\bf B}_{2r}(X_V)$ that
\begin{eqnarray*}
\tilde{A}(\xi)(t)&=&e^t\bigl(-\mathfrak{K}_3\circ\mathfrak{K}_2\circ\mathfrak{L}_q\circ\mathfrak{V}(\xi)+  \mathfrak{K}_3\circ\mathfrak{K}_2\circ\mathfrak{K}_1\circ\mathfrak{L}_v\circ\mathfrak{V}(\xi)\bigr)(t)\\
&&\qquad + \bigl(\mathfrak{K}_1\circ\mathfrak{L}_v\circ\mathfrak{V}(\xi)\bigr)(t)+   c_1e^t+
c_2e^{-t},\quad\forall t\in I.
\end{eqnarray*}
Hence $\tilde{A}$ is of class $C^{k-3}$.\hfill$\Box$\vspace{2mm}

\section{Proof of Theorem~\ref{th:1.15}}\label{sec:6}
\setcounter{equation}{0}

\subsection{Basic formulas}\label{sec:6.1}

In this section  we first derive the second variation of the energy under our boundary condition
 from recent work \cite{Jav13, Jav13+, Jav14} by Javaloyes and Soares.
 Let us begin with a briefly review about the Chern connection
and covariant derivative.
Let $(x^i,y^i)$ be the
canonical coordinates on  $TM$ associated to
local coordinates $(x^i)$ on $M$.
Let $\mathfrak{g}_{ij}(x,y)={g}^F(x,y)[\partial_{x^i}, \partial_{x^j}]$
and matrices $(\mathfrak{g}^{ij}(x,y))=(\mathfrak{g}_{ij}(x,y))^{-1}$.
Following \cite[Chapter 2]{BaChSh} the Cartan tensor $\mathcal{C}^F$ of $F$ has components given by
\begin{equation}\label{e:6.1}
\mathcal{C}_{ijm}(x,y)=\mathcal{C}^F_y\left(\partial_{x^i}, \partial_{x^j}, \partial_{x^m}\right)=
\frac{1}{4}\frac{\partial^3}{\partial y^i\partial y^j\partial y^m}F^2\left(x,y^l\partial_{x^l}|_x\right).
\end{equation}
Set $\mathcal{C}^i_{jm}:= \mathfrak{g}^{il}\mathcal{C}_{ljm}$,
$\gamma^i_{jm}(x,y):=\frac{1}{2}\mathfrak{g}^{is}\left(\partial_{x^m} \mathfrak{g}_{sj}-\partial_{x^s}\mathfrak{g}_{jm}+
\partial_{x^j}\mathfrak{g}_{ms}\right)$
and 
$$
N^i_j(x,y):=\gamma^i_{jm}(x,y)y^m-\mathcal{C}^i_{jm}(x,y)\gamma^m_{rs}(x,y)y^ry^s.
$$
 The \textsf{Chern connection} is a linear connection $\nabla$ on
 the pulled-back bundle $\pi^\ast(TM)$ over $TM\setminus 0_{TM}$
with ($C^{k-3}$) \textsf{Christoffel symbols}
\begin{equation}\label{e:6.2}
\Gamma^i_{jm}(x,y)=\gamma^i_{jm}- \mathfrak{g}^{il}\left(\mathcal{C}_{ljs}N^s_m-\mathcal{C}_{jms}N^s_i+ \mathcal{C}_{mls}N^s_j\right),
\end{equation}
that is,
$\nabla_{\partial_{x^i}}\partial_{x^j}(x,y)=\Gamma^i_{jm}(x,y)\partial_{x^m}$.
Let $R_{j\,\,ml}^{\,\,\,i}(x,y)$ be defined by
\cite[(3.3.2)]{BaChSh}, and 
$$
P_{j\,\,ml}^{\,\,\,i}(x,y)=-\partial_{\partial y^l}\Gamma_{\,\,jm}^i(x,y).
$$
For $y\in T_xM\setminus\{0\}$ the trilinear map $R_y$ (resp. $P_y$)
from $T_xM\times T_xM\times T_xM$ to $T_xM$ given by
$R_y(v,u)w=v^m u^l w^j R_{j\,\,ml}^{\,\,\,i}(x,y) \partial_{x^i}|_x$
(resp. $P_y(v,u,w)=v^j u^m w^l P_{j\,\,ml}^{\,\,\,i}(x,y)\partial_{x^i}|_x$)
define  the \textsf{Chern-Riemannian curvature tensor} (or
\textsf{$hh$-curvature tensor} \cite[Exercise 3.9.6]{BaChSh})  $R_Y$
and the \textsf{Chern-Nonriemannian curvature tensor} (or the Chern curvature
\cite[Page 112]{ShenW}) $P_Y$ with every nowhere vanishing vector field $Y$
on an open subset of $M$. According to
 \cite[\S 2]{Ma} (see \cite[\S5.2]{ShenW} and \cite{Ra042} for details ), for each nowhere vanishing $C^{k-2}$ vector field $V$ on an open subset $\Omega\subset M$ the \textsf{Chern connection} of $(M,F)$ at $V$ in $\Omega\subset M$
is the unique affine connection
$\nabla^V$ on $T\Omega$ which is torsion-free and
almost $g_V$-compatible, where $g_V(x)=\mathfrak{g}^F(x,V(x))\;\forall x\in\Omega$.  Let $\tilde\Gamma_{ij}^l$ be
the formal Christoffel symbols of $\nabla^V$ in the coordinates $(x^i)$, i.e.
$\nabla^V_{\partial_{x^i}}\partial_{x^j}(x)=\tilde\Gamma_{ij}^l(x) \partial_{x^l}\;\forall i,j$.
Then
$\tilde\Gamma^i_{jl}(x)=\Gamma^i_{jl}(x,V(x))$. The \textsf{curvature tensor} $R^V$ of $\nabla^V$ is defined by
$R^V(X,Y)Z=\nabla^V_X\nabla^V_YZ-\nabla^V_Y\nabla^V_XZ-
\nabla^V_{[X,Y]}Z$ for vector fields $X,Y,Z$ in $\Omega$. It and $R_V$, $P_V$ are related by \cite[Theorem~2.1]{Jav13+}.


For a $W^{1,2}_{\rm loc}$ curve $c$ from $I$ or $\R$ to $M$,
and $l\in\{0,1\}$ let
$W^{l,2}_{\rm loc}(c^\ast TM)$
denote the space of all   $W^{l,2}_{\rm loc}$ vector fields along $c$.
Then  $\dot{c}\in L^2_{\rm loc}(c^\ast TM):=W^{0,2}_{\rm loc}(c^\ast TM)$. Let
$(x^i, y^i)$ be the canonical coordinates around $\dot{c}(t)\in TM$.
Write
$\dot c(t)=\dot{c}^i(t) \partial_{x^i}|_{c(t)}$ and
$\zeta(t)=\zeta^i(t)\partial_{x^i}|_{c(t)}$
for $\zeta\in W^{1,2}_{\rm loc}(c^\ast TM)$. Given a nowhere vanishing $\xi\in C^0(c^\ast TM)$
the \textsf{covariant derivative} of $\zeta$ along  $c$ (with $\xi$ as reference vector)  is defined by
\begin{equation}\label{e:6.3}
D^\xi_{\dot{c}}\zeta(t):= \bigl(\dot{\zeta}^m(t)
+  \zeta^i(t)\dot{c}^j(t)\Gamma_{ij}^m(c(t), \xi(t))\bigr)\partial_{x^m}|_{c(t)}.
\end{equation}
$D^\xi_{\dot{c}}\zeta$ belongs to  $L^{2}_{\rm loc}$, and sits in $C^{\min\{k-3,r\}}(c^\ast TM)$ provided
 $c$ is of class $C^{r+1}$, $\zeta\in C^{r+1}(c^\ast TM)$ and $\xi\in C^r(c^\ast TM)$ for some $r\in\N\cup\{0,\infty\}$.  $D^\xi_{\dot{c}}\zeta(t)$ depends only on $\xi(t)$, $\dot{c}(t)$ and behavior of $\zeta$ near $t$; and $D^\xi_{\dot{c}}\zeta(t)=\nabla^{\tilde\xi}_{\dot{c}}\tilde\zeta(c(t))$ if $\dot{c}(t)\ne 0$
 and $\tilde\xi$ and $\tilde\zeta$ are any extensions of $\xi$ and $\zeta$ near $c(t)$.
When the above $\xi$ belongs to $W^{1,2}_{\rm loc}(c^\ast TM)$, $D^\xi_{\dot{c}}$  is \textsf{ almost $g_\xi$-compatible}, that is, for any $\eta,\zeta\in W^{1, 2}_{\rm loc}(c^\ast TM)$ we have
\begin{equation}\label{e:6.4}
\frac{d}{dt}g_{\xi}(\zeta,\eta)=g_{\xi}\bigl(D^\xi_{\dot{c}}\zeta,\eta\bigr)
+g_{\xi}\bigl(\zeta, D^\xi_{\dot{c}}\eta\bigr)
+2\mathcal{C}^F_{\xi}\bigl(D^\xi_{\dot{c}}\xi,\zeta,\eta\bigr)\quad{\rm a.e.}
\end{equation}
(cf. \cite[(4)]{Jav13}). If $c$ is  $C^2$ and regular,  then $c$ is a constant speed geodesic if and only if $D^{\dot c}_{\dot{c}}\dot{c}(t)\equiv 0$. For $\gamma\in\Lambda(M,\mathbb{I})$ and a nowhere vanishing $\xi\in T_\gamma\Lambda(M,\mathbb{I})$ we have also
\begin{eqnarray}\label{e:6.5}
D^{\xi}_{\dot\gamma}(\mathbb{I}_\ast\eta)(t)=\mathbb{I}_\ast\bigl(D^{\xi}_{\dot\gamma}\eta(t)\bigr)
\quad\hbox{a.e.}\;t\in\R,\;\forall \eta\in W^{1, 2}_{\rm loc}(\gamma^\ast TM),
\end{eqnarray}
 which implies  $D^{\xi}_{\dot\gamma}\zeta(t+1) =\mathbb{I}_\ast\bigl(D^\xi_{\dot\gamma}\zeta(t)\bigr)$ for a.e.
$t\in\R$ and for every $\zeta\in T_\gamma\Lambda(M,\mathbb{I})$.


Recall that $F$ is of class $C^k$, $k\ge 5$. $\mathcal{L}|_{\mathcal{X}}$ is $C^{k-2}$ by Lemma~\ref{lem:2.3}.
Each nontrivial critical point $\gamma$ of ${\cal L}$ in $\Lambda(M, \mathbb{I})$
is a nonconstant $C^{k}$ smooth $\mathbb{I}$-invariant geodesic of constant speed.
 For $u,w\in T_{\gamma({t})}M$ let $V$ be an extension of $\dot\gamma(t)$, and let $U$ and $W$ be extensions of  $u$ and $w$ onto an open neighborhood of $\gamma({t})$, respectively; and  define
 $R^{\dot\gamma(t)}(\dot\gamma({t}),u)w:=(R^V(V,U)W)({\gamma({t})})$.
The right side is independent of the chosen extensions.

\begin{proposition}\label{prop:6.1}
For a nontrivial critical point $\gamma$ of ${\cal L}$ and any $\xi,\eta\in T_\gamma{\cal X}(M,\mathbb{I})$,
\begin{eqnarray}\label{e:6.6}
d^2(\mathcal{L}|_{\mathcal{X}})(\gamma)[\xi,\eta]
&=&\int_0^1 \bigl(g_{\dot\gamma}\bigl(D^{\dot\gamma}_{\dot\gamma}\xi(t), D^{\dot\gamma}_{\dot\gamma}\eta(t)\bigr)
- g_{\dot\gamma}\bigl(R^{\dot\gamma}(\dot\gamma,\xi)\eta,\dot\gamma\bigr)\bigr)dt\nonumber\\
&=&\int_0^1 \bigl(g_{\dot\gamma}\bigl(D^{\dot\gamma}_{\dot\gamma}\xi(t),D^{\dot\gamma}_{\dot\gamma}\eta(t)\bigr)
- g_{\dot\gamma}\bigl(R^{\dot\gamma}(\xi, \dot\gamma)\dot\gamma, \eta\bigr)\bigr)dt.
\end{eqnarray}
\end{proposition}

This formula can also be seen on the pages 35-36 of \cite{Ma}.

\noindent{\bf Proof of Proposition~\ref{prop:6.1}}.
Firstly, \textsf{we assume that $\xi\in T_\gamma\Lambda(M,\mathbb{I})$ is $C^3$ and has no zeroes.}
Taking  a $C^3$ variation
$\Gamma:\R\times(-\epsilon, \epsilon)\to M$, $\Gamma=\Gamma(t,s)$
of the curve $\gamma = \Gamma_0=\Gamma(\cdot,0)$
 such that $U(\cdot,s):=\partial_s\Gamma(\cdot,s)\in
T_{\Gamma_s}\Lambda(M,\mathbb{I})\cap C^2(\Gamma_s^\ast TM)$ and  $U(t,0)=\xi(t)$
for all $t\in\R$ and all $s\in (-\epsilon,\epsilon)$, where $\Gamma_s(t):=\Gamma(t,s)$.
We can assume that each $\Gamma_s$ is  regular by shrinking $\epsilon>0$ and
   that each $\Gamma^t:(-\epsilon,\epsilon)\to M$ given by
$\Gamma^t(s)=\Gamma(t,s)$ is regular. Then  $T:=\partial_t\Gamma$
 and $U=\partial_s\Gamma$ are $C^2$ vector fields along $\Gamma$ without any zeros.
For every fixed $t$ or $s$ we have
covariant derivatives  along $\Gamma^t$ or $\Gamma_s$:
\begin{equation}\label{e:6.7}
D^T_{U(t,\cdot)}U(t,s)\quad\hbox{or}\quad D^T_{T(\cdot,s)}U(t,s)
\end{equation}
From the proof of \cite[Proposition~3.2]{Jav13} we have (in our notations)
\begin{eqnarray*}
\frac{d^2}{ds^2}\mathcal{L}(\Gamma_s)
=\int_0^1 \left(g_{\dot\Gamma_s}\Bigl(D^{\dot\Gamma_s}_{\dot\Gamma^t}D^{\dot\Gamma_s}_{\dot\Gamma_s}\dot\Gamma^t,
\dot\Gamma_s\Bigr)+
g_{\dot\Gamma_s}\Bigl(D^{\dot\Gamma_s}_{\dot\Gamma_s}\dot\Gamma^t, D^{\dot\Gamma_s}_{\dot\Gamma^t}\dot\Gamma_s\Bigl)\right)dt
\end{eqnarray*}
But according to \cite[Theorem~1.1]{Jav13+} we have
$$
D^{\dot\Gamma_s}_{\dot\Gamma_s}D^{\dot\Gamma_s}_{\dot\Gamma^t}\dot\Gamma^t-
D^{\dot\Gamma_s}_{\dot\Gamma^t}D^{\dot\Gamma_s}_{\dot\Gamma_s}\dot\Gamma^t=R^\Gamma(\dot\Gamma^t)=
R_{\dot\gamma}(\dot\gamma, \xi)\dot\Gamma^t-P_{\dot\gamma}(\xi, \dot\Gamma^t, D^{\dot\gamma}_{\dot\gamma}\dot\gamma)
-P_{\dot\gamma}(\dot\gamma, \dot\Gamma^t, D^{\dot\gamma}_{\dot\gamma}\xi)
$$
and hence using $D^{\dot\Gamma_s}_{\dot\Gamma^t}\dot\Gamma_s=D^{\dot\Gamma_s}_{\dot\Gamma_s}\dot\Gamma^t$
by \cite[Proposition~3.2]{Jav13} we deduce
\begin{eqnarray}\label{e:6.8}
&&\frac{d^2}{ds^2}\mathcal{L}(\Gamma_s)|_{s=0}
=\int_0^1 \left(
g_{\dot\Gamma_s}\Bigl(D^{\dot\Gamma_s}_{\dot\Gamma_s}\dot\Gamma^t, D^{\dot\Gamma_s}_{\dot\Gamma_s}\dot\Gamma^t\Bigl)+ g_{\dot\Gamma_s}\Bigl(D^{\dot\Gamma_s}_{\dot\Gamma^t}D^{\dot\Gamma_s}_{\dot\Gamma_s}\dot\Gamma^t,
\dot\Gamma_s\Bigr)\right)\Bigm|_{s=0}dt\nonumber\\
&&=\int_0^1g_{\dot\gamma}\bigl(D^{\dot\gamma}_{\dot\gamma}\xi, D^{\dot\gamma}_{\dot\gamma}\xi\bigl)dt+
\int_0^1 g_{\dot\gamma}\bigl(D^{\dot\gamma}_{\dot\gamma}D^{\dot\gamma}_{\xi}\xi,
\dot\gamma\bigr)dt-\int_0^1 g_{\dot\gamma}\bigl(R_{\dot\gamma}(\dot\gamma, \xi)\xi, \dot\gamma\bigr)dt
\nonumber\\
&&+\int_0^1 g_{\dot\gamma}\bigl(P_{\dot\gamma}(\xi, \xi, D^{\dot\gamma}_{\dot\gamma}\dot\gamma),
\dot\gamma\bigr)dt+
\int_0^1 g_{\dot\gamma}\bigl(P_{\dot\gamma}(\dot\gamma, \xi, D^{\dot\gamma}_{\dot\gamma}\xi),\dot\gamma\bigr)dt.
\end{eqnarray}
Moreover $g_{\dot\gamma}\bigl(D^{\dot\gamma}_{\dot\gamma}D^{\dot\gamma}_{\xi}\xi,
\dot\gamma\bigr)=\frac{d}{dt}g_{\dot\gamma}\bigl(D^{\dot\gamma}_{\xi}\xi,
\dot\gamma\bigr)$, $D^{\dot\gamma}_{\dot\gamma}\xi=D^{\dot\gamma}_{\xi}\dot\gamma$ and
$$
R^{\dot\gamma}(\dot\gamma, \xi)\xi=
R_{\dot\gamma}(\dot\gamma, \xi)\xi+P_{\dot\gamma}(\xi, \xi, D^{\dot\gamma}_{\dot\gamma}\dot\gamma)
-P_{\dot\gamma}(\dot\gamma, \xi, D^{\dot\gamma}_\xi{\dot\gamma})
$$
by \cite[Theorem~2.1]{Jav13+}. These and (\ref{e:6.8}) lead to
\begin{eqnarray}\label{e:6.9}
\frac{d^2}{ds^2}\mathcal{L}(\Gamma_s)\left|_{s=0}\right.
&=&\int_0^1 \bigr(g_{\dot\gamma}\bigl(D^{\dot\gamma}_{\dot\gamma}\xi(t), D^{\dot\gamma}_{\dot\gamma}\xi(t)\bigr)-
g_{\dot\gamma}(R^{\dot\gamma}(\dot\gamma,\xi)\xi,\dot\gamma)\bigl)dt\nonumber\\
&&+ g_{\dot\gamma(t)}\bigl(D^T_{U(t,\cdot)}U(t,s)|_{s=0},\dot\gamma(t)\bigr)\Bigm|^{t=1}_{t=0}.
\end{eqnarray}

\begin{claim}\label{cl:6.2}
$g_{\dot\gamma(1)}\bigl(D^{T(1,\cdot)}_{U(1,\cdot)}U(1,s)|_{s=0},\dot\gamma(1)\bigr)
=g_{\dot\gamma(0)}\bigl(D^{T(0,\cdot)}_{U(0,\cdot)}U(0,s)|_{s=0},\dot\gamma(0)\bigr)$.
\end{claim}
In fact, let us  choose a system of coordinates $(x^i)$ on an open neighborhood $\Omega$ of $\gamma(0)$.
Then $(\tilde x^i)=(x^i\circ\mathbb{I}^{-1})$ forms
a  system of coordinates  on that  of $\gamma(1)$, $\tilde\Omega:=\mathbb{I}(\Omega)$, and
$\mathbb{I}_\ast(\partial_{x^i}|_p)=\partial_{\tilde x^i}|_{\tilde{p}}$ for $p\in\Omega$ and $\tilde{p}=\mathbb{I}(p)$, $i=1,\cdots,n$.
For $v=y^i\partial_{x^i}|_p$ we have $\tilde{v}:=\mathbb{I}_\ast(v)=
\partial_{\tilde x^i}|_{\tilde p}$, i.e.,
$(\tilde{y}^i(\mathbb{I}_\ast v))=(y^i(v))$.
Hence $\tilde x^i(\tilde{p}, \tilde{v})=x^i(p,v)$ and $\tilde{y}^i(\tilde{p}, \tilde{v})= y^i(p,v)$
for $i=1,\cdots,n$.
 Since
$$
\tilde{\mathfrak{g}}_{ij}(\tilde x,\tilde y)={g}^F(\tilde x,\tilde y)\left[\partial_{\tilde x^i}, \partial_{\tilde x^j}\right]
={g}^F(\mathbb{I}x, \mathbb{I}_\ast y)\left[\mathbb{I}_\ast\left(\partial_{x^i}\right),
\mathbb{I}_\ast\left(\partial_{x^j}\right)\right]=\mathfrak{g}_{ij}(x,y),
$$
we get
$\partial_{\tilde x^m}\tilde{\mathfrak{g}}_{ij}(\tilde x,\tilde y)
=\partial_{x^m}\mathfrak{g}_{ij}(x, y)$ and thus $\gamma^i_{jm}(x,y)=\tilde\gamma^i_{jm}(\tilde x,\tilde y)$.
Using (\ref{e:6.1}) it easily proved that
$\tilde{\mathcal{C}}_{ijm}(\tilde x,\tilde y)=\mathcal{C}_{ijm}(x,y)$
and hence
$\tilde{\mathcal{C}}^i_{jm}(\tilde x,\tilde y)=\tilde{\mathfrak{g}}^{il}(\tilde x,\tilde y)\mathcal{C}_{ljm}(\tilde x,\tilde y)=\mathfrak{g}^{il}(x,y)\mathcal{C}_{ljm}(x,y)=\mathcal{C}^i_{jm}(x,y)$ and
\begin{eqnarray*}
\tilde N^i_j(\tilde x,\tilde y)&=&\tilde\gamma^i_{jm}(\tilde x,\tilde y)\tilde y^m-\tilde{\mathcal{C}}^i_{jm}(\tilde x,
\tilde y)\tilde\gamma^m_{rs}(\tilde x,\tilde y)y^r\tilde y^s\\
&=&\gamma^i_{jm}(x,y)y^m-\mathcal{C}^i_{jm}(x,y)\gamma^m_{rs}(x,y)y^ry^s=N^i_j(x,y).
\end{eqnarray*}
Let $\tilde\Gamma^i_{jl}(\tilde x,\tilde y)$ be the  Christoffel symbols in the local coordinates
$(\tilde{x}^i, \tilde{y}^i)$. We obtain
\begin{equation}\label{e:6.10}
\tilde\Gamma^i_{jl}(\tilde x,\tilde y)=\Gamma^i_{jl}(x,y).
\end{equation}
Let $|t|+|s|$ be small enough. We may write
$T(t,s)=T^i(t,s)\partial_{x^i}|_{\Gamma(t,s)}$,
$U(t,s)=U^i(t,s)\partial_{x^i}|_{\Gamma(t,s)}$,
$\gamma^i(t)=x^i\circ\gamma(t)$, $i=1,\cdots,n$, $T(t,0)=\dot{\gamma}(t)=\gamma^i(t)\partial_{x^i}|_{\gamma(t)}$.
Since $D^{\dot\gamma(0)}_{\partial_{x^i}}\partial_{x^j}(\gamma(0))=\Gamma_{ij}^k(\dot\gamma(0)) \partial_{x^k}|_{\gamma(0)}$
for $i,j=1,\ldots,n$, by  (\ref{e:6.7}) we have
$$
D^{T(0,\cdot)}_{U(0,\cdot)}U(0,s)|_{s=0}=\frac{{\rm \partial}U^i}{{\rm \partial}s}(0,0)\partial_{x^i}|_{\gamma(0)}+  U^i(0,0)\frac{{\rm \partial}U^j}{{\rm \partial}s}(0,0)\Gamma_{\,\,ij}^k(\dot{\gamma}(0))\partial_{x^k}|_{\gamma(0)},
$$
\begin{eqnarray}\label{e:6.11}
&&g_{\dot\gamma(0)}\bigl(D^{T(0,\cdot)}_{U(0,\cdot)}U(0,s)|_{s=0},\dot\gamma(0)\bigr)
=\frac{{\rm \partial}U^i}{{\rm \partial}s}(0,0)
T^l(0,0)g_{\dot\gamma(0)}\left(\partial_{x^i}|_{\gamma(0)},\partial_{x^l}|_{\gamma(0)}\right)\nonumber\\
&&+U^i(0,0)\frac{{\rm \partial}U^j}{{\rm \partial}s}(0,0)\Gamma_{\,\,ij}^m(\dot{\gamma}(0))T^l(0,0)g_{\dot\gamma(0)}\left(
\partial_{x^m}|_{\gamma(0)},\partial_{x^l}|_{\gamma(0)}\right).
\end{eqnarray}
From $\Gamma(t+1,s)=\mathbb{I}(\Gamma(t,s))\;\forall (t,s)$ it follows that
$T(t+1,s)=\mathbb{I}_\ast(T(t,s))$ and $U(t+1,s)=\mathbb{I}_\ast(U(t,s))$ for all
$(t,s)$.
If $|t-1|+|s|$ is small enough we get
\begin{eqnarray*}
&&T(t+1,s)=T^i(t,s)\partial_{\tilde x^i}|_{\Gamma(t+1,s)},\quad
U(t+1,s)=U^i(t,s)\partial_{\tilde x^i}|_{\Gamma(t+1,s)},\\
&&\tilde\gamma^i(t+1):=\tilde x^i\circ\gamma(t+1)=\tilde x^i(\mathbb{I}(\gamma(t)))=x^i(\gamma(t))=\gamma^i(t),\;i=1,\cdots,n,\\ &&T(t+1,0)=\dot{\gamma}(t+1)=\gamma^i(t)\partial_{\tilde x^i}|_{\gamma(t+1)}.
\end{eqnarray*}
These and $D^{\dot\gamma(1)}_{\partial_{\tilde x^i}}\partial_{\tilde x^j}(\gamma(1))=\tilde\Gamma_{ij}^m(\dot\gamma(1)) \partial_{\tilde x^m}|_{\gamma(1)}$
for $i,j=1,\ldots,n$, lead to
$$
D^{T(1,\cdot)}_{U(1,\cdot)}U(1,s)|_{s=0}=\frac{{\rm \partial}U^i}{{\rm \partial}s}(0,0)\partial_{\tilde x^i}|_{\gamma(1)}+  U^i(0,0)\frac{{\rm \partial}U^j}{{\rm \partial}s}(0,0)\tilde\Gamma_{\,\,ij}^m(\dot{\gamma}(1))\partial_{\tilde x^m}|_{\gamma(1)}.
$$
 By (\ref{e:6.10})
 $\tilde\Gamma_{ij}^m(\dot\gamma(1))=\tilde\Gamma_{ij}^m(\mathbb{I}_\ast\dot\gamma(0))=\Gamma_{ij}^m(\dot\gamma(0))$
for any $i,j,m=1,\cdots,n$. Hence
\begin{eqnarray*}
&&g_{\dot\gamma(1)}\bigl(D^{T(1,\cdot)}_{U(1,\cdot)}U(1,s)|_{s=0},\dot\gamma(1)\bigr)
=\frac{{\rm \partial}U^i}{{\rm \partial}s}(0,0)
T^l(0,0)g_{\dot\gamma(1)}\left(\partial_{\tilde x^i}|_{\gamma(1)},\partial_{\tilde x^l}|_{\gamma(1)}\right)\nonumber\\
&&+U^i(0,0)\frac{{\rm \partial}U^j}{{\rm \partial}s}(0,0)\Gamma_{\,\,ij}^m(\dot{\gamma}(0))T^l(0,0)g_{\dot\gamma(1)}\left(
\partial_{\tilde x^m}|_{\gamma(1)},\partial_{\tilde x^l}|_{\gamma(1)}\right).
\end{eqnarray*}
Since
$g_{\dot\gamma(0)}\left(
\partial_{x^m}|_{\gamma(0)},\partial_{x^l}|_{\gamma(0)}\right)
=g_{\dot\gamma(1)}\left(\partial_{\tilde x^m}|_{\gamma(1)},\partial_{\tilde x^l}|_{\gamma(1)}\right)$,
Claim~\ref{cl:6.2} follows from this  and (\ref{e:6.11}) immediately.

Now (\ref{e:6.9}) and Claim~\ref{cl:6.2} yield
\begin{eqnarray}\label{e:6.12}
d^2(\mathcal{L}|_{\mathcal{X}})(\gamma)[\xi,\xi]
=\int_0^1 \bigl(g_{\dot\gamma}\bigl(D^{\dot\gamma}_{\dot\gamma}\xi, D^{\dot\gamma}_{\dot\gamma}\xi\bigr)-
g_{\dot\gamma}(R^{\dot\gamma}(\dot\gamma,\xi)\xi, \dot\gamma)
\bigr)dt.
\end{eqnarray}
This also holds  for any $\xi\in T_\gamma\Lambda(M,\mathbb{I})\cap C^3(\gamma^\ast TM)$ because of the density of the set 
$$
\{\xi\in T_\gamma\Lambda(M,\mathbb{I})\cap C^3(\gamma^\ast TM)\,|\, \xi(t)\ne 0\;\forall t\}
$$
in $T_\gamma\Lambda(M,\mathbb{I})\cap C^3(\gamma^\ast TM)$ with respect to  $C^3$-topology.
It follows for any $\xi,\eta\in T_\gamma\Lambda(M,\mathbb{I})\cap C^3(\gamma^\ast TM)$ that
\begin{eqnarray}\label{e:6.13}
d^2(\mathcal{L}|_{\mathcal{X}})(\gamma)[\xi,\eta]&=&\frac{1}{2}\bigl(
d^2(\mathcal{L}|_{\mathcal{X}})(\gamma)[\xi+\eta,\xi+\eta]-d^2(\mathcal{L}|_{\mathcal{X}})(\gamma)[\xi,\xi]-
d^2(\mathcal{L}|_{\mathcal{X}})(\gamma)[\eta,\eta]\bigr)\nonumber\\
&=&\int_0^1 g_{\dot\gamma}\left(D^{\dot\gamma}_{\dot\gamma}\xi(t), D^{\dot\gamma}_{\dot\gamma}\eta(t)\right)dt
-\frac{1}{2}\int^1_0g_{\dot\gamma}(R^{\dot\gamma}(\dot\gamma,\xi)\eta,\dot\gamma)dt\nonumber\\
&&-\frac{1}{2}\int^1_0g_{\dot\gamma}(R^{\dot\gamma}(\dot\gamma,\eta)\xi,\dot\gamma)dt.
\end{eqnarray}
Observe that \cite[(11)]{Jav13} and \cite[Lemma~3.10]{Jav14} imply
\begin{eqnarray*}
&&g_{\dot\gamma}(R^{\dot\gamma}(\dot\gamma, \eta)\xi,\dot\gamma)-g_{\dot\gamma}(R^{\dot\gamma}(\xi,\dot\gamma)\dot\gamma,\eta)
=B^{\dot\gamma}(\xi,\eta,\dot\gamma,\dot\gamma)+B^{\dot\gamma}(\dot\gamma,\xi,\eta,
\dot\gamma)+B^{\dot\gamma}(\dot\gamma,\dot\gamma,
\xi,\eta)\\
&&\hspace{20mm}+ B^{\dot\gamma}(\eta,\dot\gamma,\xi,\dot\gamma)+B^{\dot\gamma}(\dot\gamma,\xi,\dot\gamma,
\eta)+B^{\dot\gamma}(\dot\gamma,\eta,\xi,\dot\gamma)=0.
\end{eqnarray*}
Hence $g_{\dot\gamma}(R^{\dot\gamma}(\dot\gamma, \eta)\xi,\dot\gamma)=g_{\dot\gamma}(R^{\dot\gamma}(\xi, \dot\gamma)\dot\gamma, \eta)=-g_{\dot\gamma}(R^{\dot\gamma}(\dot\gamma, \xi)\dot\gamma, \eta)
=g_{\dot\gamma}(R^{\dot\gamma}(\dot\gamma,\xi)\eta, \dot\gamma)$ by \cite[(14)]{Jav13+}.
This and (\ref{e:6.13}) show that the expected equalities hold
for any $\xi,\eta\in T_\gamma\Lambda(M,\mathbb{I})\cap C^3(\gamma^\ast TM)$.
Finally, by the density of $T_\gamma\Lambda(M,\mathbb{I})\cap C^3(\gamma^\ast TM)$
in $T_\gamma{\cal X}(M,\mathbb{I})$ we complete the proof.
\hfill$\Box$\vspace{2mm}

Since $\gamma$ is a nonconstant geodesic (and so regular), for each fixed $t$ we choose
 a $C^3$ vector field  $V$ near $\gamma({t})$ such that
  $V(\gamma(s))=\dot\gamma(s)$ for $|s-{t}|\ll 1$ and that it is also
  a geodesic vector field (i.e. $\nabla^V_VV=0$). Then
   $R^V$ is the Riemannian curvature tensor of the metric $g_V$ near  $\gamma({t})$.
 It follows that
 \begin{equation}\label{e:6.14}
|g_{\dot\gamma}(R^{\dot\gamma}(\xi(t),\dot\gamma({t}))\dot\gamma(t),\eta(t))|\le
3\max|K(\gamma(t))||v|^2_{g_{\dot\gamma({t})}}
|\xi(t)|_{g_{\dot\gamma({t})}}|\eta(t)|_{g_{\dot\gamma({t})}},
 \end{equation}
 where $|v|_{g_{\dot\gamma({t})}}=(g_{\dot\gamma({t})}(v,v))^{1/2}$ and $\max|K(\gamma(t))|$ is the maximum of the sectional curvatures of  $g_V$ at $\gamma({t})$.
This implies that the right side of (\ref{e:6.6}) can be extended
into a continuous symmetric bilinear form on $T_\gamma\Lambda(\mathbb{I}, M)$, still
denoted by $d^2(\mathcal{L}|_{\mathcal{X}})(\gamma)$ unless otherwise stated.

Introduce another inner-product $\langle\!\langle \cdot,\cdot\rangle\!\rangle_{1,\gamma}$
on $T_\gamma\Lambda(M,\mathbb{I})$ given by
$$
\langle\!\langle\xi,\eta\rangle\!\rangle_{1,\gamma}:=\int_0^1 \bigl(g_{\dot\gamma}(\xi, \eta)+ g_{\dot\gamma}\bigl(D^{\dot\gamma}_{\dot\gamma}\xi(t), D^{\dot\gamma}_{\dot\gamma}\eta(t)\bigr)
\bigr)dt.
$$
It and the inner-product in (\ref{e:1.1}) induce equivalent norms.
Similarly, we have also an equivalent inner-product $\langle\!\langle \cdot,\cdot\rangle\!\rangle_{0,\gamma}$ to  $\langle \cdot,\cdot\rangle_{0,\gamma}$.
Since $\gamma$ is a geodesic,  (\ref{e:6.4}) yields
\begin{eqnarray*}
\frac{d}{dt}g_{\dot\gamma}\bigl(D^{\dot\gamma}_{\dot\gamma}\xi,\eta\bigr)&=&
g_{\dot\gamma}\bigl((D^{\dot\gamma}_{\dot\gamma})^2\xi, \eta\bigr)+g_{\dot\gamma}\bigl(D^{\dot\gamma}_{\dot\gamma}\xi, D^{\dot\gamma}_{\dot\gamma}\eta\bigr)
+2\mathcal{C}^F_{\dot\gamma}\bigl(D^{\dot\gamma}_{\dot\gamma}\dot\gamma, D^{\dot\gamma}_{\dot\gamma}\xi,\eta\bigr)\nonumber\\
&=&g_{\dot\gamma}\bigl((D^{\dot\gamma}_{\dot\gamma})^2\xi, \eta\bigr)+g_{\dot\gamma}\bigl(D^{\dot\gamma}_{\dot\gamma}\xi, D^{\dot\gamma}_{\dot\gamma}\eta\bigr)
\end{eqnarray*}
for any $\xi,\eta\in T_\gamma\Lambda(\mathbb{I}, M)\cap C^{k-2}(\gamma^\ast TM)$.
From this and (\ref{e:6.5})-(\ref{e:6.6}) we derive
\begin{eqnarray}\label{e:6.15}
d^2(\mathcal{L}|_{\mathcal{X}})(\gamma)[\xi,\eta]
&=&\int_0^1 \bigl(g_{\dot\gamma}\bigl(-R^{\dot\gamma}(\xi, \dot\gamma)\dot\gamma, \eta\bigr)- g_{\dot\gamma}\bigl((D^{\dot\gamma}_{\dot\gamma})^2\xi(t),\eta(t)\bigr)
\bigr)dt\nonumber\\
&&\qquad+g_{\dot\gamma(1)}\bigl(D^{\dot\gamma}_{\dot\gamma}\xi(1),\eta(1)\bigr)- g_{\dot\gamma(0)}\bigl(D^{\dot\gamma}_{\dot\gamma}\xi(0),\eta(0)\bigr)\nonumber\\
&=&\int_0^1 \bigl(g_{\dot\gamma}\bigl(-R^{\dot\gamma}(\xi, \dot\gamma)\dot\gamma, \eta\bigr)- g_{\dot\gamma}\bigl((D^{\dot\gamma}_{\dot\gamma})^2\xi(t),\eta(t)\bigr)
\bigr)dt\nonumber\\
&=&\int_0^1 g_{\dot\gamma}(\mathfrak{L}_\gamma\xi, \eta)dt=\langle\!\langle \mathfrak{L}_\gamma\xi,\eta\rangle\!\rangle_{0,\gamma}.
\end{eqnarray}
 Here $\mathfrak{L}_\gamma:T_\gamma\Lambda(\mathbb{I}, M)\cap W^{2,2}((\gamma|_I)^\ast TM)\to
 L^{2}((\gamma|_I)^\ast TM)$ with $I=[0,1]$  is  defined by
$\mathfrak{L}_\gamma\xi=-(D^{\dot\gamma}_{\dot\gamma})^2\xi - R^{\dot\gamma}(\xi, \dot\gamma)\dot\gamma$.
 which maps $T_\gamma\Lambda(\mathbb{I}, M)\cap C^{k-2}(\gamma^\ast TM)$ into
 $T_\gamma\Lambda(\mathbb{I}, M)\cap C^{k-4}(\gamma^\ast TM)$ by (\ref{e:6.5}). Observe that (\ref{e:6.6}) and (\ref{e:6.14}) imply
$$
d^2(\mathcal{L}|_{\mathcal{X}})(\gamma)[\xi,\xi]\ge\langle\!\langle\xi,\xi\rangle\!\rangle_{1,\gamma}-
{\rm const}\langle\!\langle\xi, \xi\rangle\!\rangle_{0,\gamma}\quad\forall \xi\in T_\gamma\Lambda(M, \mathbb{I}).
$$
Hence  $\mathfrak{L}_\gamma$ is
an essentially self-adjoint  elliptic operator. By Theorem~22.G in \cite[\S 22.16b]{Ze}
each spectrum point of $\mathfrak{L}_\gamma$ is a real eigenvalue of finite multiplicity and all
eigenvalues of it (counted according to their multiplicity) are
$$
-{\rm const}\le\mu_1=\min\{d^2(\mathcal{L}|_{\mathcal{X}})(\gamma)[\xi,\xi]\,|\,\xi\in T_\gamma\Lambda(M, \mathbb{I}),\; \langle\!\langle\xi, \xi\rangle\!\rangle_{0,\gamma}=1\}\le\mu_2\le\cdots,
$$
 and $\mu_m\to\infty\;\hbox{as}\;m\to\infty$.
All  eigenvectors are in $T_\gamma\Lambda(\mathbb{I}, M)\cap C^{k-2}(\gamma^\ast TM)$,
any two eigenvectors  corresponding with two distinct eigenvalues are orthogonal
 with respect to the inner-product $\langle\!\langle\cdot, \cdot\rangle\!\rangle_{0,\gamma}$,
and there exists a complete orthogonal system of eigenvectors in the completion of
 $T_\gamma\Lambda(\mathbb{I}, M)\cap C^{k-2}(\gamma^\ast TM)$ with respect to
$\langle\!\langle\cdot, \cdot\rangle\!\rangle_{0,\gamma}$.

Call $\xi\in C^{k-2}(\gamma^\ast TM)$  a $C^{k-2}$ \textsf{Jacobi field}
along $\gamma$ if it satisfies the \textsf{Jacobi equation}:
$(D^{\dot\gamma}_{\dot\gamma})^2\xi(t)+ R^{\dot\gamma}(\xi, \dot\gamma)\dot\gamma=0$.
Such a field has the same properties as in  Riemannian case, see
 \cite[Lemmas~3.16,~3.17 and Proposition~3.18]{Jav14}. Let $\mathscr{J}_\gamma$ be the set of $C^{k-2}$ Jacobi fields along $\gamma$. It is a vector space isomorphic to $T_{\gamma(0)}M\times T_{\gamma(0)}M$. By (\ref{e:6.15})
$$
{\bf H}^0(d^2(\mathcal{L}|_{\mathcal{X}})(\gamma))=\{\xi\in \mathscr{J}_\gamma\,|\,
\mathbb{I}_\ast\xi(t)=\xi(t+1),\;\forall t\in\R\}={\rm Ker}(\mathfrak{L}_\gamma).
$$
Since $\mathfrak{L}_\gamma\dot\gamma=0$, each
eigenvector $\xi$ corresponding to a nonzero eigenvalue sits in
$$
T_\gamma\Lambda(M, \mathbb{I})\cap \mathcal{V}_\gamma=\{\xi\in C^{k-2}(\gamma^\ast TM)\,|\, \langle\!\langle\xi, \dot\gamma\rangle\!\rangle_{0,\gamma}=0,\;\mathbb{I}_\ast\xi={\bf T}_1\xi\}
$$
where ${\bf T}_1\xi$ is defined as  in (\ref{e:1.15}) and
$\mathcal{V}_\gamma=\{\xi\in C^{k-2}(\gamma^\ast TM)\,|\,
\langle\!\langle\xi, \dot\gamma\rangle\!\rangle_{0,\gamma}=0\}$.
For $\xi\in\mathscr{J}_\gamma$, since  $g_{\dot\gamma(t)}(\xi(t),\dot\gamma(t))=at+b,\;\forall t\in\R$,
by the proof of \cite[Lemma~3.17]{Jav14}, where $a$ and $b$ are real constants, we obtain that
$\langle\!\langle\xi, \dot\gamma\rangle\!\rangle_{0,\gamma}=0$ if and only if $g_{\dot\gamma(t)}(\xi(t),\dot\gamma(t))=0\;\forall t$. It follows that
$\dim\mathscr{J}_\gamma\cap\mathcal{V}_\gamma=\dim\mathscr{J}_\gamma-2=2\dim M-2$.
Clearly, every $\xi\in\mathscr{J}_\gamma$ has a decomposition $({\langle\!\langle\xi, \dot\gamma\rangle\!\rangle_{0,\gamma}/\langle\!\langle\dot\gamma, \dot\gamma\rangle\!\rangle_{0,\gamma}})\dot\gamma+\xi^\bot$,
where $\xi^\bot\in \mathscr{J}_\gamma\cap\mathcal{V}_\gamma$. Define
\begin{equation}\label{e:6.16}
\mathbb{L}_\gamma: \mathcal{V}_\gamma\to \mathcal{V}_\gamma\cap C^{k-4}(\gamma^\ast TM),\quad\xi\mapsto
-(D^{\dot\gamma}_{\dot\gamma})^2\xi - R^{\dot\gamma}(\xi, \dot\gamma)\dot\gamma,
\end{equation}
which coincides with $\mathfrak{L}_\gamma$ on $T_\gamma\Lambda(M, \mathbb{I})\cap \mathcal{V}_\gamma$.
 Define  the \textsf{index} $\lambda(\gamma, \mathbb{I})$
(resp. the \textsf{nullity} $\nu(\gamma, \mathbb{I})$)  of the geodesic $\gamma$
  as the index of the quadratic form
$d^2(\mathcal{L}|_{\mathcal{X}})(\gamma)$ on $T_\gamma\Lambda(\mathbb{I}, M)$
(resp. the nullity of $d^2(\mathcal{L}|_{\mathcal{X}})(\gamma)$ on $T_\gamma\Lambda(\mathbb{I}, M)$
minus one). Then
\begin{eqnarray}
\lambda(\gamma, \mathbb{I})&=&\sum_{\mu<0}\dim\bigl\{\xi\in T_\gamma\Lambda(M, \mathbb{I})\cap C^{k-2}(\gamma^\ast TM)\,|\,
 \mathfrak{L}_\gamma\xi=\mu\xi\bigr\}\nonumber\\
&=&\sum_{\mu<0}\dim\bigl\{\xi\in \mathcal{V}_\gamma\,|\,
\mathbb{I}_\ast \xi={\bf T}_1\xi,\; \mathbb{L}_\gamma\xi=\mu\xi\bigr\},\label{e:6.17}\\
\nu(\gamma, \mathbb{I})&=&\dim\bigl\{\xi\in T_\gamma\Lambda(M, \mathbb{I})\cap C^{k-2}(\gamma^\ast TM)\,|\, \mathfrak{L}_\gamma\xi=0\bigr\}-1\nonumber\\
&=&\dim\bigl\{\xi\in \mathcal{V}_\gamma\,|\, \mathbb{I}_\ast \xi={\bf T}_1\xi,\;
\mathbb{L}_\gamma\xi=0\bigr\}.\label{e:6.18}
\end{eqnarray}
Moreover, if the $\R$-orbit of $\gamma$ in $\Lambda(\mathbb{I}, M)$ is an embedded $S^1$,
 by (\ref{e:1.20}) we have
\begin{equation}\label{e:6.19}
\lambda(\gamma, \mathbb{I})=m^-({\cal L}, \R\cdot\gamma)\quad\hbox{and}\quad \nu(\gamma, \mathbb{I})=m^0({\cal L}, \R\cdot\gamma).
\end{equation}

\textsf{From now on} we assume that
  $\gamma$ is a closed
$\mathbb{I}$-invariant $F$-geodesic of constant speed and with least period $\alpha>0$. Then
for all $m\in\N$, $t\in\R$ we have
$\gamma^{m\alpha+1}(t+1)=\mathbb{I}(\gamma^{m\alpha+1}(t))$ and $({\bf T}_{m\alpha+1}\gamma)(t)=\mathbb{I}(\gamma(t))=\mathbb{I}\circ \gamma(t)=\gamma(t+1)=({\bf T}_1\gamma)(t)$.
It follows that  ${\bf T}_\tau\gamma^{m\alpha+1}$, $m\in\N$ and $\tau\in\R$, are critical points of ${\cal L}$
in $\Lambda(M,\mathbb{I})$. Moreover both $\mathbb{I}_\ast$ and ${\bf T}_\tau$ induce \textsf{bijective linear maps}, still denoted by $\mathbb{I}_\ast$ and
${\bf T}_\tau$,
 ${\mathbb{I}}_\ast:{\cal V}_\gamma\to {\cal V}_{{\bf T}_1(\gamma)}$ and ${\bf T}_\tau: {\cal V}_\gamma\to {\cal V}_{{\bf T}_\tau(\gamma)}$.
 For each  $\mu\le 0$ and each integer $m\in\N$ let
\begin{equation}\label{e:6.20}
\begin{array}{rl}
 &{\cal V}_{\gamma^{m\alpha+1}}(\mu,1, {\mathbb{I}}):=\bigl\{\xi\in {\cal V}_{\gamma^{m\alpha+1}}\,|\,
 \mathbb{L}_{\gamma^{m\alpha+1}}\xi=\mu \xi,\;{\mathbb{I}}_\ast\xi={\bf T}_1\xi\bigr\},
 \\[12pt]
 &{\cal V}_{\gamma}(\mu, m\alpha+1, {\mathbb{I}}):=\bigl\{\xi\in {\cal V}_{\gamma}\,|\, \mathbb{L}_{\gamma}\xi=\mu\xi,\;
 {\mathbb{I}}_\ast\xi={\bf T}_{m\alpha+1}\xi\bigr\}.
 \end{array}
\end{equation}
  Here $\mathbb{L}_{\gamma^{m\alpha+1}}$ is defined as in (\ref{e:6.16}).  By  (\ref{e:6.17}) and (\ref{e:6.18}) we obtain
$$
\lambda(\gamma^{m\alpha+1}, \mathbb{I})=\sum_{\mu<0}\dim {\cal V}_{\gamma^{m\alpha+1}}(\mu, 1, {\mathbb{I}})\quad\hbox{and}\quad
\nu(\gamma^{m\alpha+1}, \mathbb{I})=\dim {\cal V}_{\gamma^{m\alpha+1}}(0,1,{\mathbb{I}}).
$$
They imply that the spaces in (\ref{e:6.20}) are finite dimensional.
 For $\xi\in{\cal V}_\gamma$ define $\xi^{m\alpha+1}(t)=\xi((m\alpha+1)t),\;\forall t\in\R$, and
 $ \iota_m(\xi)=\xi^{m\alpha+1}/(m\alpha+1)^2$.
  It is easily checked that $\iota_m\bigl({\cal V}_{\gamma}(\mu, m\alpha+1, {\mathbb{I}})\bigr)\subset {\cal V}_{\gamma^{m\alpha+1}}(\mu,1,{\mathbb{I}})$ and that   $\iota_m:{\cal V}_\gamma\to {\cal V}_{\gamma^{m\alpha+1}}$
   is a linear injection.   For  $\eta\in {\cal V}_{\gamma^{m\alpha+1}}(\mu,1,{\mathbb{I}})$ set
 $\xi(t):=(m\alpha+1)^2\eta(t/(m\alpha+1)^2),\;\forall t\in\R$.
  it is easy to prove that $\xi\in {\cal V}_{\gamma}(\mu, m\alpha+1,{\mathbb{I}})$ and  $\iota_m(\xi)=\eta$.
  These show that  $\iota_m$ restricts to an  isomorphism from ${\cal V}_{\gamma}(\mu, m\alpha+1,{\mathbb{I}})$ onto
   ${\cal V}_{\gamma^{m\alpha+1}}(\mu,1,{\mathbb{I}})$.
   Thus
 $\dim{\cal V}_{\gamma}(\mu, m\alpha+1,{\mathbb{I}})=\dim{\cal V}_{\gamma^{m\alpha+1}}(\mu,1,{\mathbb{I}})$
  and hence
\begin{equation}\label{e:6.21}
\begin{array}{rl}
&\lambda(\gamma^{m\alpha+1}, \mathbb{I})=\sum_{\mu<0}\dim {\cal V}_\gamma(\mu, m\alpha+1,{\mathbb{I}})\\
&\hspace{25mm}=\sum_{\mu<0}\dim\bigl\{\xi\in{\cal V}_\gamma\,|\, \mathbb{L}_\gamma\xi=\mu\xi,\;\mathbb{I}_\ast \xi={\bf T}_{m\alpha+1}\xi\bigr\}
\\[12pt]
&\nu(\gamma^{m\alpha+1}, \mathbb{I})=\dim {\cal V}_\gamma(0, m\alpha+1,{\mathbb{I}})\\
&\hspace{25mm}=\dim\bigl\{\xi\in {\cal V}_\gamma\,|\, \mathbb{L}_\gamma\xi=0,\;\mathbb{I}_\ast \xi={\bf T}_{m\alpha+1}\xi\bigr\}.
\end{array}
\end{equation}
Set $J_\gamma(\mu):=\bigl\{\xi\in {\cal V}_\gamma\,|\, \mathbb{L}_\gamma\xi=\mu \xi\bigr\}$
for $\mu\in\R$. Its dimension is finite because
 $J_\gamma(\mu)\cap T_\gamma\Lambda(\mathbb{I}, M)=
\bigl\{\xi\in {\cal V}_\gamma\,|\, \mathbb{L}_\gamma\xi=\mu \xi\;\&\;\mathbb{I}_\ast\xi={\bf T}_1\xi\bigr\}$
is contained in the finite dimension  vector space
$\{\xi\in T_\gamma\Lambda(\mathbb{I}, M)\cap C^{k-2}(\gamma^\ast TM)\,|\, \mathfrak{L}_\gamma\xi=\mu \xi\}$
which vanishes if $\mu\notin\{\mu_k\}_{k=1}^\infty$
by the argument below (\ref{e:6.15}).
Moreover, $J_\gamma(0)={\rm Ker}(\mathbb{L}_\gamma)=\mathscr{J}_\gamma\cap\mathcal{V}_\gamma$ and so
$\dim J_\gamma(0)=2\dim M-2$. It is clear  that
${\cal V}_{\gamma}(\mu, m\alpha+1,{\mathbb{I}})=\bigl\{\xi\in J_\gamma(\mu)\,|\,\mathbb{I}_\ast\xi={\bf T}_{m\alpha+1}\xi\bigr\}$.
 Define a \textsf{positive definite bilinear form} $\omega_m$ on ${\cal V}_\gamma$
 (resp. ${\cal V}_{{\bf T}_1(\gamma)}$)   by
 $$
 \omega^\gamma_m(\xi,\eta)=\int^{m\alpha+1}_0g_{\dot\gamma}(\xi(t),\eta(t))dt\quad{\rm (resp.}\;
 \omega^{{\bf T}_1\gamma}_m(\xi,\eta)=\int^{m\alpha+1}_0g_{{\bf T}_{1}\dot\gamma}(\xi(t),\eta(t))dt
 \;{\rm )}
 $$
 Then  $\mathbb{I}_\ast:{\cal V}_\gamma\to {\cal V}_{{\bf T}_1(\gamma)}$ preserves the inner-product, and the restriction of ${\bf T}_1$ to  $J_\gamma(\mu)$ also preserves inner-products. It follows that
 ${\bf T}_1^{-1}\circ \mathbb{I}_\ast$ is an orthogonal transformation on the inner product space
 $\bigl({\cal V}_{\gamma},\omega_m\bigr)$ and hence on
 $\bigl({\cal V}_{\gamma}(\mu, m\alpha+1,{\mathbb{I}}),\omega_m\bigr)$ because the latter
 is an invariant subspace of ${\bf T}_1^{-1}\circ \mathbb{I}_\ast$.
Consider the complexifications of ${\cal V}_{\gamma}$ and $J_\gamma(\mu)$,
${\cal V}_{\gamma}\otimes\C$ and $J_\gamma(\mu)\otimes\C$,
and also use $\mathbb{L}_\gamma$, $\mathbb{I}_\ast$ and ${\bf T}_\tau$ to denote  the $\C$-linear extensions
of  $\mathbb{L}_\gamma$, $\mathbb{I}_\ast$ and ${\bf T}_\tau$, respectively.
 For $z\in\C$ consider the complex linear space
$$
S[\mu,m\alpha+1, z\mathbb{I}]:=\bigl\{\xi\in J_\gamma(\mu)\otimes\C\,|\,z\mathbb{I}_\ast \xi={\bf T}_{m\alpha+1}\xi\bigr\}
$$
 with Hermitian extension inner product $\tilde\omega^\gamma_m$ of $\omega^\gamma_m$  to
$J_\gamma(\mu)\otimes\C$ given by
$$
\tilde\omega^\gamma_m(\alpha+i\beta, a+ib):=\omega_m(\alpha,a)+\omega_m(\beta,b)+i\omega_m(\beta,a)-i\omega_m(\alpha,b)
$$
for any $\alpha+i\beta, a+ib\in J(\mu)\otimes\C$. Note that
$S[\mu,m\alpha+1, \mathbb{I}]={\cal V}_\gamma(\mu,m\alpha+1, \mathbb{I})\otimes\C$.
It is easy to check that
 ${\bf T}_1^{-1}\circ \mathbb{I}_\ast$  is a unitary transformation on the Hermitian  inner product space
$\bigl(S[\mu,m\alpha+1, \mathbb{I}], \tilde\omega^\gamma_m\bigr)$.
As in the proof of  \cite[Lemma~1.2]{Tan82} we have a direct sum decomposition of
 finitely many of its subspaces:
 \begin{equation}\label{e:6.22}
 S[\mu,m\alpha+1, \mathbb{I}]=\bigoplus_{|z|=1}\bigoplus_{\rho^m=z}\left\{
\xi\in S[\mu,1,z^{-1}\mathbb{I}]\;|\;{\bf T}_\alpha \xi=\rho \xi\right\},
 \end{equation}
 where $z$ takes over all eigenvalues of ${\bf T}_1^{-1}\circ \mathbb{I}_\ast$ and
 $\oplus_{\rho^m=z}\left\{\xi\in S[\mu,1,z^{-1}\mathbb{I}]\;|\;{\bf T}_\alpha \xi=\rho \xi\right\}$
 is the corresponding eigenspace with $z$; moreover  $|z|$ is the modulus of $z$.
For each $z\in{\bf U}:=\{\rho\in\C\,|\,|\rho|=1\}$ define functions
\begin{eqnarray*}
&&\Lambda^z(\rho)=\sum_{\mu<0}\dim_{\C}\left\{\xi\in  S[\mu,1,z^{-1}\mathbb{I}]\,|\, {\bf T}_\alpha \xi=\rho \xi\right\}\quad\forall\rho\in{\bf U},\\
&&N^z(\rho)=\dim_{\C}\left\{\xi\in  S[0,1, z^{-1}\mathbb{I}]\,|\, {\bf T}_\alpha \xi=\rho \xi\right\}\quad\forall\rho\in{\bf U}.
\end{eqnarray*}
Let $({\cal V}_{\gamma}\otimes\C)_z=\{\xi\in{\cal V}_{\gamma}\otimes\C\,|\, z\mathbb{I}_\ast\xi={\bf T}_1\xi\}$
and denote by $\mathbb{L}^z_\gamma$ the restriction of  $\mathbb{L}_\gamma$ to $({\cal V}_{\gamma}\otimes\C)_z$.
Then $S[0,1, z^{-1}\mathbb{I}]={\rm Ker}(\mathbb{L}_\gamma^z)\subset \{\xi\in {\cal V}_{\gamma}\otimes\C\,|\,
  \mathbb{L}_\gamma\xi=0\}=J_\gamma(0)\otimes\C$ and so $\dim_{\C}S[0,1, z^{-1}\mathbb{I}]\le 2\dim M-2$.
As in \cite{GroTa78, Tan77I, Tan82} we can obtain

\begin{lemma}\label{lem:6.3}
\begin{description}
\item[(i)] $\Lambda^z$ and $N^z$ are identically zero except for a finite number of $z$'s,
\item[(ii)] For each $z$, $N^z(\rho)=0$ except for at most $2(\dim M-1)$ points $\rho$, which
will be called $\mathbb{L}^z_\gamma$ Poincare points.
\item[(iii)] For each $z$, $\Lambda^z$ is locally constant except possibly at the $\mathbb{L}^z_\gamma$ Poincar\'e
points.
\item[(iv)] For each $z$, the inequality $\lim_{\rho\to\rho_0}\Lambda^z(\rho)\ge\Lambda^z(\rho_0)$
 holds for any $\rho_0$.
\item[(v)] (\ref{e:6.21}) and (\ref{e:6.22}) imply
\begin{eqnarray*}
\lambda(\gamma^{m\alpha+1}, \mathbb{I})=\sum_{|z|=1}\sum_{\rho^m=z}\Lambda^z(\rho)\quad\hbox{and}\quad
\nu(\gamma^{m\alpha+1}, \mathbb{I})=\sum_{|z|=1}\sum_{\rho^m=z}N^z(\rho).
\end{eqnarray*}
\end{description}
\end{lemma}

As in \cite{GroTa78} Lemma~\ref{lem:3.1}(iii)-(v) lead to:

\begin{lemma}\label{lem:6.4}
 Either $\lambda(\gamma^{m\alpha+1}, \mathbb{I})=0$ for all nonnegative integers $m$ or there
exist positive numbers $\varepsilon$ and $a$ such that for all integers $m_1\ge m_2\ge 0$,
$$
\lambda(\gamma^{m_1\alpha+1}, \mathbb{I})-\lambda(\gamma^{m_2\alpha+1}, \mathbb{I})\ge (m_1- m_2)\varepsilon-a.
$$
\end{lemma}

\subsection{Proof of Theorem~\ref{th:1.15}}\label{sec:6.2}

Let $\Lambda(M, \mathbb{I})$ be equipped with the Riemannian-Hilbert structure
given by (\ref{e:1.1}) with the $\mathbb{I}$-invariant Riemannian metric $g$ on $M$.
  Since $\mathbb{I}$ has only finitely many $\mathbb{I}$-invariant $F$-geodesics by  assumption,
   all of them must be closed  and thus there exist finitely many $\mathbb{I}$-invariant closed
   $F$-geodesics, $\gamma_1,\cdots,\gamma_r$, such that any $\mathbb{I}$-invariant $F$-geodesic lies in one
   of critical orbits $\R\cdot\gamma_i^{m\alpha_i+1}$ of ${\cal L}$ in $\Lambda(M, \mathbb{I})$, $i=1,\cdots,r$, $m=0,1,2,\cdots$, where $\alpha_i$ denotes the least period of $\gamma_i$.
 Of course  these critical orbits are all isolated.

\begin{proposition}\label{prop:6.5}
 There exists a constant $B$ such that
$\dim\mathscr{H}^0_k({\cal L}, \gamma_i^{m\alpha_i+1};\K)\le B$  for all integers $k, m\ge 0$
and $i=1,\cdots,r$.
\end{proposition}

  We postpone the proof of it until next subsection.
  This and (\ref{e:1.30}) lead to
   $\dim C_k({\cal L}; \R\cdot\gamma_i^{m\alpha_i+1};\K)\le 4B$ for any $k,m\ge 0$ and $i=1,\cdots,r$.
   Moreover, by (\ref{e:1.29}) we have for all $q\ge 2\dim M-1$, all $m\ge 0$ and
    $i=1,\cdots,r$,
    $$C_{q+\nu(\gamma_i^{m\alpha_i+1}, \mathbb{I})}({\cal L}; \R\cdot\gamma_i^{m\alpha_i+1};\K)=\mathscr{H}_q^0({\cal L},\gamma_i^{m\alpha+1};\K)=0.
    $$
   As in \cite{GroTa78} it follows from this and Lemma~\ref{lem:6.4} that
    $$
\mathfrak{C}:=\sup\Bigl\{\sum_m \dim C_k({\cal L}, \R\cdot{\gamma}_i^{m\alpha_i+1};\K)\,\Bigm|\, k\ge 2\dim M,\;i=1,\cdots,r\Bigr\}<\infty,
$$
 which, as in \cite{GM2, GroTa78}, leads to
 $\dim H_k\left(\Lambda(M,\mathbb{I})^{a_2}, \Lambda(M, \mathbb{I})^{a_1}\right)\le 4B\mathfrak{C}\;\forall k\ge 2\dim M$ for all regular values $0<a_1<a_2$. Since ${\rm Crit}({\cal L})\cap{\cal L}^{-1}(0, a]=\emptyset$ for $0<a<\min_i{\cal L}(\gamma_i)$,  ${\rm Fix}(\mathbb{I})=\{{\cal L}\le 0\}$ is a strong deformation retract of $\mathcal{L}^{-1}([0,a])=\Lambda(M,\mathbb{I})^{a}$ by  \cite[Page 21, Lemma~3.2]{Ch93} or \cite[Page 181, Lemma~8.3]{MW}.
These and $\dim{\rm Fix}(\mathbb{I})\le\dim M$ yield
$$
\sup\{\dim H_k(\Lambda(M, \mathbb{I}))\,|\, k\ge 2\dim M\}\le 4B\mathfrak{C}
$$
as in the proof of \cite[Theorem~4.1]{GroTa78}. The desired conclusion follows from
the fact that the inclusion $\Lambda(M,\mathbb{I})\subset C^0(I, M)^{\mathbb{I}}$
 is a homotopy equivalence \cite[Theorem~1.3]{Gro73I}.
\hfill$\Box$\vspace{2mm}

\subsection{Proof of Proposition~\ref{prop:6.5}}\label{sec:6.3}

Let $(\gamma,\alpha)$ be one of $(\gamma_i, \alpha_i)$, $i=1,\cdots,\gamma_r$.
We proceed in two cases as in \cite{Tan82}.

\subsubsection{$\alpha$ is rational}\label{sec:6.3.1}
There exist  unique
relatively prime positive integers $m_0$ and $s_0$ with $\alpha=s_0/m_0$. ($s_0$ is actually
 the smallest positive integer with the property $\gamma(\R)\subset{\rm Fix}(\mathbb{I}^{s_0})$.)
 Since $\dim_{\C}J_\gamma(0)\otimes\C=2(\dim M-1)$ and $\mathbb{I}_\ast\circ \mathbb{L}_\gamma=\mathbb{L}_{\mathbb{I}(\gamma)}\circ
 \mathbb{I}_\ast$ the unitary transformation
 $\mathbb{I}_\ast^{s_0}: J_\gamma(0)\otimes\C\to J_\gamma(0)\otimes\C$
 has at most $2(\dim M - 1)$ eigenvalues,
 $ \exp(2\pi i\theta_1),\cdots, \exp(2\pi i\theta_k)$, where $0\le\theta_1<\theta_2<\cdots<\theta_k<1$.
If $\theta_i>0$ is rational, choose positive integers $p_i$ and $q_i$ which are relatively prime and satisfies
$p_i/q_i=\theta_i$.  Define $s=s_0$ if all $\theta_i$'s are irrational, and
$s$ to be the product of $s_0$ and the least
common multiple of the $q_i's$ with $\theta_i$ rational otherwise. Then
 $\mathbb{I}^{s}(\gamma(t))=\gamma(t)\;\forall t$ because
 $\mathbb{I}^{s_0}(\gamma(t))=\gamma(t),\;\forall t$.
 Let $\mathbb{J}$ be the restriction of $\mathbb{I}$ to ${\rm Fix}(\mathbb{I}^s)$.
 Observe that ${\rm Fix}(\mathbb{I}^s)$ is a collection of closed totally geodesic
submanifolds of $(M,g)$ and that  $\gamma$ may be viewed as a critical point of ${\cal L}$ on the Hilbert manifold $\Lambda({\rm Fix}(\mathbb{I}^s), \mathbb{J})$. Actually,  any critical point of ${\cal L}|_{\Lambda({\rm Fix}(\mathbb{I}^s), \mathbb{J})}$ is one of ${\cal L}$.

\begin{lemma}\label{lem:6.6}
$\exists\, k_0\in\N$ such that
$\nu(\gamma^{m\alpha+1}, \mathbb{I})=\nu(\gamma^{m\alpha+1}, \mathbb{J})
\;\forall m\ge k_0$.
\end{lemma}

\noindent{\bf Proof}. By Lemma~\ref{lem:6.3}(i)  there exist only finitely many $z_j\in{\bf U}$, $j=1,\cdots,l$,
such that $N^{z_j}\ne 0$ for $j=1,\cdots,l$. Hence Lemma~\ref{lem:6.3}(v) leads to
\begin{eqnarray}\label{e:6.23}
\nu(\gamma^{m\alpha+1},\mathbb{I})=\sum^l_{j=1}\sum_{\rho^m=z_j}N^{z_j}(\rho).
\end{eqnarray}
Write $z_j=\exp(i{\rm arg}(z_j))$ with ${\rm arg}(z_j)\in [0,2\pi)$. It is not hard to prove that
 there are at most finitely many integers $m$ with
$N^{z_j}(\rho)\ne 0$ and $\rho^m=z_j$ provided $\beta_j:={\rm arg}(z_j)/(2\pi)$ is irrational.
This implies that there  exists
 $k_0\in\N$ such that for every $m\ge k_0$ and every $j$ with irrational $\beta_j$ we have
$N^{z_j}(\rho)=0$ for every $\rho$ with $\rho^m=z_j$.
Thus (\ref{e:6.23}) gives rise to
\begin{eqnarray}\label{e:6.24}
\nu(\gamma^{m\alpha+1},\mathbb{I})=\sum_{\hbox{rational}\;\beta_j}\sum_{\rho^m=z_j}N^{z_j}(\rho),\quad\forall
m\ge k_0.
\end{eqnarray}
For the finite set $\mathcal{Q}:=\{\beta_j\,|\,\beta_j={\rm arg}(z_j)/(2\pi)\;\hbox{is rational}\}$ we
can find a positive integer $r$ such that $r\beta_j\in\Z\;\forall\beta_j\in\mathcal{Q}$ and so
$(z_j)^r=\exp(ir{\rm arg}(z_j))=\exp\left(2\pi i\beta_jr\right)=1$ for every $\beta_j\in\mathcal{Q}$.
This and (\ref{e:6.24}) yield
\begin{eqnarray}\label{e:6.25}
\nu(\gamma^{m\alpha+1},\mathbb{I})=\sum_{z^r=1}\sum_{\rho^m=z}N^{z}(\rho),\quad\forall
m\ge k_0.
\end{eqnarray}
Suppose that  $N^z(\rho)$ in the right side is positive for some $z$ and $\rho$. Then
$$
\{\xi\in S[0,1,z^{-1}\mathbb{I}]\,|\, {\bf T}_\alpha \xi=\rho \xi\}=
S[0,1,z^{-1}\mathbb{I}]\cap S[0,\alpha,\rho {\rm id}_M]\ne\emptyset
$$
and so  $\rho^{m_0}z^{s_0}$ is an eigenvalue
 of $\mathbb{I}_\ast^{s_0}$ on $\{\xi\in {\cal V}_\gamma\otimes\C\,|\, \mathbb{L}_\gamma\xi=0\}$.
These lead to
\begin{eqnarray*}
\bigoplus_{z^r=1}\bigoplus_{\rho^m=z}
S[0,1,z^{-1}\mathbb{I}]\cap S[0,\alpha,\rho {\rm id}_M]
&\subset&\left\{\xi\in
S[0, m\alpha+1,\mathbb{I}]\,|\,\mathbb{I}_\ast^{s_0}\xi=\xi\right\}\\
&\subset&\left\{\xi\in
S[0, m\alpha+1,\mathbb{I}]\,|\,\mathbb{I}_\ast^{s}\xi=\xi\right\}
\end{eqnarray*}
because $s=s_0\tau$ with $\tau\in\N$ and hence
$\mathbb{I}_\ast^{s_0}\xi=\xi$ gives $\mathbb{I}_\ast^{s}\xi=\xi$.
By (\ref{e:6.22})
\begin{eqnarray*}
 && \left\{\xi\in
S[0, m\alpha+1,\mathbb{I}]\,|\,\mathbb{I}_\ast^{s}\xi=\xi\right\}= S[0,m\alpha+1, \mathbb{I}]\cap{\rm Ker}(\mathbb{I}^s_\ast-id)\\
&=&\bigoplus_{|z|=1}\bigoplus_{\rho^m=z}\left\{
 X\in S[0,1,z^{-1}\mathbb{I}]\;|\;{\bf T}_\alpha \xi=\rho \xi\right\}\cap{\rm Ker}(\mathbb{I}^s_\ast-id).
 \end{eqnarray*}
From these and (\ref{e:6.25}) it follows that
\begin{eqnarray*}
\nu(\gamma^{m\alpha+1}, \mathbb{I})&=&\bigoplus_{z^r=1}\bigoplus_{\rho^m=z}
\dim_{\C}S[0,1,z^{-1}\mathbb{I}]\cap S[0,\alpha,\rho {\rm id}_M]\\
&\le&
\bigoplus_{|z|=1}\bigoplus_{\rho^m=z}\dim_{\C}\left\{
 \xi\in S[0,1,z^{-1}\mathbb{I}]\;|\;{\bf T}_\alpha \xi=\rho \xi\right\}\cap{\rm Ker}(\mathbb{I}^s_\ast-id)\\
 &=&\nu(\gamma^{m\alpha+1}, \mathbb{J})\quad\hbox{for all}\;m\ge k_0.
 \end{eqnarray*}
The definition of nullities yields the converse inequality  $\nu(\gamma^{m\alpha+1},\mathbb{J})\le \nu(\gamma^{m\alpha+1}, \mathbb{I})$.
 \hfill$\Box$\vspace{2mm}

By Lemma~\ref{lem:6.6} and  Corollary~\ref{cor:1.11}
 we obtain
\begin{eqnarray}\label{e:6.26}
\mathscr{H}^0_\ast({\cal L}, \gamma^{m\alpha+1};\K)=\mathscr{H}^0_\ast\left({\cal L}\bigm|_{\Lambda({\rm Fix}(\mathbb{I}^s), \mathbb{J})}, \gamma^{m\alpha+1};\K\right),\quad\forall m\ge k_0.
\end{eqnarray}
Put $\widehat{M}=\Lambda({\rm Fix}(\mathbb{I}^s), \mathbb{J})$ and $\widehat{\cal L}={\cal L}\bigm|_{\Lambda({\rm Fix}(\mathbb{I}^s), \mathbb{J})}$. Then $\mathbb{J}\in I(\widehat{M}, F|_{\hat{M}})\cap I(\widehat{M}, g|_{\widehat{M}})$ is of order $s$.
As in the proof of \cite[Lemma~2.9]{GroTa78} we can show

\begin{lemma}\label{lem:6.7}
There exist distinct positive integers ${k}_1,\cdots, k_q$,
where $q$  has an upper bound only depending
 on $\dim \widehat{M}$ and $s/s_0$, and for each integer $j\in [1, q]$ a strictly increasing infinite sequence
$\{{m}^i_j| j=1,2,\cdots\}$ of positive integers
 such that
\begin{description}
\item[(i)]  the sets $N_j:=\{{m}^i_jk_j| i=1,2,\cdots\}$ form a partition of
$\{ms_0+m_0\,|\, m\in\N\cup\{0\}\}=
\{m_0(m\alpha+ 1)\,|\, m\in\N\cup\{0\}\}$;
\item[(ii)] for the maximal integer  ${s}^i_j$ relatively prime to ${m}^i_j$
  and dividing $s/s_0$,
 $$
 \nu\left(\bar{\gamma}^{{m}^i_j{k}_j}, \mathbb{J}\right)=\nu\left(\bar{\gamma}^{{m}^i_j{k}_j}, \mathbb{J}\bigm|{\rm Fix}\left(\mathbb{J}^{s_0s^i_j}\right)  \right)=\nu\left(\bar{\gamma}^{{k}_j}, \mathbb{J}^r\bigm|{\rm Fix}\left(\mathbb{J}^{s_0{s}^i_j}\right)  \right),
 $$
  where $\bar{\gamma}=\gamma^{1/m_0}$ and $r$ is an integer
with the property $r{m}^i_j\equiv 1\;{\rm mod}\;s_0{s}^i_j$.
 \end{description}
\end{lemma}

By the first equality in Lemma~\ref{lem:6.7}(ii) and Corollary~\ref{cor:1.11} we get
$$
\mathscr{H}^0_\ast(\widehat{\cal L}, \bar{\gamma}^{{m}^i_j{k}_j})\cong\mathscr{H}^0_\ast\bigl(\widehat{\cal L}|\Lambda(S,  \mathbb{J}), \bar{\gamma}^{{m}^i_j{k}_j}\bigr)\quad\hbox{with}\quad
S={\rm Fix}\bigl(\mathbb{J}^{s_0{s}^i_j}\bigr)
$$
Moreover the second equality in Lemma~\ref{lem:6.7}(ii) and Theorem~\ref{th:1.12}.
lead to
$$
\mathscr{H}^0_\ast\bigl(\widehat{\cal L}|\Lambda(S, \mathbb{J}^r|_S), \bar{\gamma}^{{k}_j}\bigr)\cong\mathscr{H}^0_\ast\bigl(\widehat{\cal L}|\Lambda\bigl(S,  \mathbb{J}\bigr), \bar{\gamma}^{{m}^i_j{k}_j}\bigr)\quad\hbox{with}\quad S={\rm Fix}\bigl(\mathbb{J}^{s_0{s}^i_j}\bigr)
$$
since
$\left(\mathbb{J}^r|_S\right)^{{m}^i_j}=\mathbb{J}^{r{m}^i_j}|_S
=\mathbb{J}|_S$
by the condition $r{m}^i_j\equiv 1\;{\rm mod}\;s_0{s}^i_j$. Hence
\begin{eqnarray}\label{e:6.27}
\mathscr{H}^0_\ast(\widehat{\cal L}, \bar{\gamma}^{{m}^i_j{k}_j})\cong
\mathscr{H}^0_\ast\bigl(\widehat{\cal L}|\Lambda(S, \mathbb{J}^r|_S), \bar{\gamma}^{{k}_j}\bigr)\quad
\hbox{with}\quad S={\rm Fix}\bigl(\mathbb{J}^{s_0{s}^i_j}\bigr).
\end{eqnarray}
Observe that $1\le s_0{s}^i_j\le s$, $\mathbb{J}^s=id_{\widehat{M}}$ and
$\mathscr{H}^0_l\bigl(\widehat{\cal L}|\Lambda(S, \mathbb{J}^r|_S), \bar{\gamma}^{{k}_j}\bigr)=0$
for $l\ge 2\dim S-1$ by the shifting theorem. We deduce that
$$
\left\{\dim\mathscr{H}^0_\ast\bigl(\widehat{\cal L}|\Lambda(S, \mathbb{J}^r|_S), \bar{\gamma}^{{k}_j}\bigr)\bigm|
r{m}^i_j\equiv 1\;{\rm mod}\;s_0{s}^i_j,\;i\in\N, j=1,\cdots,q\right\}
$$
is a finite set, and so by (\ref{e:6.27}) the set
$\bigl\{\dim\mathscr{H}^0_\ast(\widehat{\cal L}, \bar{\gamma}^{{m}^i_j{k}_j})\bigm|
i\in\N, j=1,\cdots,q\bigr\}$
is finite. Since $\bar{\gamma}^{m_0}=\gamma$ and $\{{m}^i_j{k}_j\,|\, i\in\N, j=1,\cdots,q\}=
\{m_0(m\alpha+ 1)\,|\, m\in\N\cup\{0\}\}$ by Lemma~\ref{lem:3.5}(i) we obtain
$\sharp\bigl\{\dim\mathscr{H}^0_\ast(\widehat{\cal L}, {\gamma}^{m\alpha+1})\bigm|
m\in\N\cup\{0\}\bigr\}<\infty$.
 This and (\ref{e:6.26}) imply
that $\dim\mathscr{H}^0_k({\cal L}, \gamma^{m\alpha+1};\K)\le B$ for some
constant $B>0$ and  all integers $k, m\ge 0$.

\subsubsection{$\alpha$ is irrational}\label{sec:6.3.2}

As in the proof of (\ref{e:6.25}) we can derive from (i),(ii) and (v) in Lemma~\ref{lem:6.3} that
there exist positive integers  $\mathfrak{m}_0$ and $\mathfrak{s}$ such that
$$
\nu(\gamma^{m\alpha+1}, \mathbb{I})=\sum_{z^{\mathfrak{s}}=1}\sum_{\rho^m=z}N^z(\rho)
\quad\forall m\ge \mathfrak{m}_0.
$$
Let $J(\rho,z):=S[0,1,z^{-1}\mathbb{I}]\cap S[0,\alpha,\rho{\rm id}_M]$. By (\ref{e:6.22}) this means
\begin{eqnarray*}
 S[0,m\alpha+1, \mathbb{I}]&=&\bigoplus_{|z|=1}\bigoplus_{\rho^m=z}\left\{
 \xi\in S[0,1,z^{-1}\mathbb{I}]\;|\;{\bf T}_\alpha \xi=\rho \xi\right\}\nonumber\\
 &=&\bigoplus_{z^{\mathfrak{s}}=1}\bigoplus_{\rho^m=z} J(\rho,z)\quad\forall m\ge \mathfrak{m}_0.
\end{eqnarray*}
Carefully checking the proof of \cite[Lemma~3.1]{Tan82} it is not difficult to see
that the corresponding result also holds in our case. Namely we have

  \begin{lemma}\label{lem:6.8}
  There exist distinct positive integers ${k}_1,\cdots, k_q$,
where $q$  has an upper bound only depending
 on $\dim M$ and $\mathfrak{s}$, and for each integer $j\in [1, q]$ a strictly increasing infinite sequence
$\{{m}^i_j| j=1,2,\cdots\}$ of positive integers
 such that
\begin{description}
\item[(i)]  the sets $N_j:=\{{m}^i_jk_j| i=1,2,\cdots\}$ form a partition of
$\N$;
\item[(ii)] for the maximal integer $s^i_j$ relatively prime $m^i_j$ dividing $\mathfrak{s}$,
\begin{eqnarray*}
\bigoplus_{z^{\mathfrak{s}}=1}\bigoplus_{\rho^{m^i_jk_j}=z}J(\rho,z)=
\bigoplus_{z^{{\mathfrak{s}}^i_j}=1}\bigoplus_{\rho^{k_j}=z}J(\rho,z).
\end{eqnarray*}
\end{description}
\end{lemma}

 For each
$z=\exp(2\pi iu/v)$  with $u,v\in\N$ satisfying $(u, v) = 1$ and $u\le v$,  let
$$
Q^z:=\left\{a\in\N\,\Bigm|\,\exists b\in\N\;{\rm s.t.}\;(b,av)=1,
 N^z\left(\exp\Bigl(\frac{2\pi ib}{av}\Bigr)\right)>0\right\}.
 $$
Then the set
$Q=\bigcup_{z^{\mathfrak{s}}=1}Q^z\cup\{1\}$ is
 finite, and that $a\in Q$ and $a|{m}^i_j{k}_j$ imply $a|{k}_j$.
 (See the proof of \cite[Lemma~3.1]{Tan82}).

\begin{lemma}\label{lem:6.9}
$m^0\bigl(^{m\alpha+1}\!{\cal L}, \R\cdot\gamma\bigr)
=m^0\bigl(\,^{m\alpha+1}\!{\cal L}^{\mathfrak{X}(m)\alpha}, \R\cdot\gamma\bigr)$
for every integer  $m=m^i_jk_j\ge m_0$ and $\mathfrak{X}(m):=s^i_jk_j$.
\end{lemma}

\noindent{\bf Proof}.
Let $H_\gamma(^{m\alpha+1}\!{\cal L})$ and $H_\gamma(^{m\alpha+1}\!{\cal L}^{\mathfrak{X}(m)\alpha})$ be
the continuous bilinear extensions of
\begin{eqnarray*}
d^2\bigl(\,^{m\alpha+1}\!{\cal L}|_{{\cal X}^{m\alpha+1}(M, \mathbb{I})}\bigr)(\gamma)\quad\hbox{and}\quad
d^2\bigl(\,^{m\alpha+1}\!{\cal L}^{\mathfrak{X}(m)\alpha}|_{{\cal X}^{m\alpha+1}(M,\mathbb{I})\cap{\cal X}^{\mathfrak{X}(m)\alpha}(M)}\bigr)(\gamma)
\end{eqnarray*}
onto $T_\gamma\Lambda^{m\alpha+1}(M,\mathbb{I})$
and $T_\gamma(\Lambda^{m\alpha+1}(M,\mathbb{I})\cap\Lambda^{\mathfrak{X}(m)\alpha}(M))$, respectively.
 Since $^{m\alpha+1}\!{\cal L}^{\mathfrak{X}(m)\alpha}$
 is the restriction of $^{m\alpha+1}\!{\cal L}$ to $\Lambda^{\mathfrak{X}(m)\alpha}(M)\cap\Lambda^{m\alpha+1}(M,\mathbb{I})$,
 $H_\gamma(^{m\alpha+1}\!{\cal L}^{\mathfrak{X}(m)\alpha})$ is
  that of
  $H_\gamma(^{m\alpha+1}\!{\cal L})$ to
$T_\gamma(\Lambda^{m\alpha+1}(M,\mathbb{I})\cap\Lambda^{\mathfrak{X}(m)\alpha}(M))$. Let $\mathscr{A}_\gamma$ be the bounded linear operator defined by
$$
H_\gamma(^{m\alpha+1}\!{\cal L})[\xi, \eta]=\langle\!\langle\mathscr{A}_\gamma\xi,\eta\rangle\!\rangle_{m\alpha+1}\quad
\hbox{for all}\;\xi,\eta\in T_\gamma\Lambda^{m\alpha+1}(M,\mathbb{I}).
$$
Because
 $ \mathscr{A}_{{\bf T}_u(\gamma)}\circ {\bf T}_u={\bf T}_u\circ \mathscr{A}_\gamma\;\forall u\in\R$,
   $\mathscr{A}_\gamma$ maps ${\bf T}_\gamma(\Lambda^{m\alpha+1}(M,\mathbb{I})\cap\Lambda^{\mathfrak{X}(m)\alpha}(M))$
 into itself, i.e., $T_\gamma(\Lambda^{m\alpha+1}(M, \mathbb{I})\cap\Lambda^{\mathfrak{X}(m)\alpha}(M))$ is an invariant subspace of $\mathscr{A}_\gamma$.
  Hence
${\rm Ker}\bigl(H_\gamma(^{m\alpha+1}\!{\cal L}^{\mathfrak{X}(m)\alpha})\bigr)$ defined by
\begin{eqnarray*}
\left\{\xi\in T_\gamma\left(\Lambda^{m\alpha+1}(M, \mathbb{I})\cap\Lambda^{\mathfrak{X}(m)\alpha}(M)\right)\,\Bigm|\,
H_\gamma\left(^{m\alpha+1}\!{\cal L}\right)[\xi,\eta]=0\;\forall \eta\in T_\gamma\Lambda^{m\alpha+1}(M,\mathbb{I})\right\}
\end{eqnarray*}
is equal to
$T_\gamma\Lambda^{\mathfrak{X}(m)\alpha}(M)\cap{\rm Ker}\left(H_\gamma(^{m\alpha+1}\!{\cal L})\right)$.
It is clear that
$$
{\cal V}_\gamma\cap{\rm Ker}\left(H_\gamma(^{m\alpha+1}\!{\cal L})\right)
=\left\{\xi\in
 {\cal V}_\gamma\,|\, \mathbb{L}_\gamma\xi=0,\;\mathbb{I}_\ast \xi={\bf T}_{m\alpha+1}\xi\right\}
={\cal V}_\gamma(0, m\alpha+1,\mathbb{I})
$$
and has the complexified  space $S[0,m\alpha+1, \mathbb{I}]$
by the definition  above (\ref{e:6.22}).
Moreover
\begin{eqnarray*}
{\cal V}_\gamma(0,\mathfrak{X}(m)\alpha, {\rm id}_M):&=&\bigl\{\xi\in J_\gamma(0)\,|\, \xi={\bf T}_{\mathfrak{X}(m)\alpha}\xi\bigr\}\\
&=&
\bigl\{\xi\in {\cal V}_\gamma\,|\, \xi={\bf T}_{\mathfrak{X}(m)\alpha}X,\;\mathbb{L}_\gamma\xi=0\bigr\}\\
&=&\bigl\{\xi\in {\cal V}_\gamma\,|\,\mathbb{L}_\gamma\xi=0\bigr\}\cap T_\gamma\Lambda^{\mathfrak{X}(m)\alpha}(M)
\end{eqnarray*}
and this and
${\cal V}_\gamma(0, m\alpha+1,{\mathbb{I}})=\bigl\{\xi\in {\cal V}_\gamma\,|\, \mathbb{L}_\gamma\xi=0,\;\mathbb{I}_\ast \xi={\bf T}_{m\alpha+1}\xi\bigr\}$ lead to
\begin{eqnarray*}
&&{\cal V}_\gamma\cap T_\gamma\Lambda^{\mathfrak{X}(m)\alpha}(M)\cap{\rm Ker}\left(H_\gamma(^{m\alpha+1}\!{\cal L})\right)\\
&=&{\cal V}_\gamma\cap {\rm Ker}\left(H_\gamma(^{m\alpha+1}\!{\cal L})\right)\cap
T_\gamma\Lambda^{\mathfrak{X}(m)\alpha}(M)\\
&=&{\cal V}_\gamma(0, m\alpha+1, \mathbb{I})
\cap\bigl\{\xi\in {\cal V}_\gamma\,|\,\mathbb{L}_\gamma\xi=0\bigr\}\cap T_\gamma\Lambda^{\mathfrak{X}(m)\alpha}(M)\\
&=&{\cal V}_\gamma(0, m\alpha+1,\mathbb{I})\cap {\cal V}_\gamma(0,\mathfrak{X}(m)\alpha, {\rm id}_M).
\end{eqnarray*}
Observe that ${\cal V}_\gamma\cap T_\gamma\Lambda^{\mathfrak{X}(m)\alpha}(M)\cap{\rm Ker}\left(H_\gamma(^{m\alpha+1}\!{\cal L})\right)$
is  the  orthogonal complementary space of $\dot\gamma$ in
${\rm Ker}\bigl(H_\gamma(^{m\alpha+1}\!{\cal L}^{\mathfrak{X}(m)\alpha})\bigr)=
T_\gamma\Lambda^{\mathfrak{X}(m)\alpha}(M)\cap{\rm Ker}\left(H_\gamma(^{m\alpha+1}\!{\cal L})\right)$.
We obtain
$$
\left({\cal V}_\gamma\cap T_\gamma\Lambda^{\mathfrak{X}(m)\alpha}(M)\cap{\rm Ker}\left(H_\gamma(^{m\alpha+1}\!{\cal L})\right)\right)\otimes\C=
S[0, m\alpha+1,\mathbb{I}]\cap S[0,\mathfrak{X}(m)\alpha,{\rm id}_M].
$$
Here $S[0,\mathfrak{X}(m)\alpha,{\rm id}_M]={\cal V}_\gamma(0,\mathfrak{X}(m)\alpha, {\rm id}_M)\otimes\C=\bigl\{\xi\in J_\gamma(0)\otimes\C\,|\, \xi={\bf T}_{\mathfrak{X}(m)\alpha}\xi\bigr\}
$
by the definition of $S[\mu,m\alpha+1, z\mathbb{I}]$ above (\ref{e:6.22}).

 On the other hand, for  $\xi\in S[0, m\alpha+1,\mathbb{I}]\setminus\{0\}$
  (\ref{e:3.8}) and (\ref{e:3.9}) imply that 
  $$
  \xi\in\bigoplus_{z^{s^i_j}=1}\bigoplus_{\rho^{k_j}=z}J(\rho,z).
  $$
It follows that $T_{\xi(m)\alpha}\xi=\rho^{s^i_jk_j}\xi=\xi$,
i.e., $\xi\in S[0,\mathfrak{X}(m)\alpha,{\rm id}_M]$.
Hence $S[0, m\alpha+1,\mathbb{I}]\subset S[0,\mathfrak{X}(m)\alpha,{\rm id}_M]$
and thus
$S[0, m\alpha+1,\mathbb{I}]\cap S[0,\mathfrak{X}(m)\alpha,{\rm id}_M]=S[0, m\alpha+1,\mathbb{I}]$.
Note that $S[0, m\alpha+1,\mathbb{I}]$ is the complexified orthogonal space of $\dot \gamma$ in ${\rm Ker}\bigl(H_\gamma(^{m\alpha+1}\!{\cal L})\bigr)$. We arrive at
${\rm Ker}\bigl(H_\gamma(^{m\alpha+1}\!{\cal L})\bigr)={\rm Ker}\bigl(H_\gamma(^{m\alpha+1}\!{\cal L}^{\mathfrak{X}(m)\alpha})\bigr)$
and hence the desired result.
 \hfill$\Box$\vspace{2mm}

Lemma~\ref{lem:6.9} and (\ref{e:1.45}) lead to
 $m^0\bigl({\cal L}, \R\cdot\gamma^{m\alpha+1}\bigr)
=m^0\bigl(\,^{m\alpha+1}\!{\cal L}^{\mathfrak{X}(m)\alpha}, \R\cdot\gamma\bigr)$
for $m\ge \mathfrak{m}_0$, and so by Theorem~\ref{th:1.13} we  arrive at

\begin{proposition}\label{prop:6.10}
For  $m=m^i_jk_j\ge \mathfrak{m}_0$ and $\mathfrak{X}(m)=s^i_jk_j$ it holds that
$$
\mathscr{H}^0_\ast(\,^{m\alpha+1}\!{\cal L}, \gamma;\K)\cong \mathscr{H}^0_\ast(\,^{m\alpha+1}\!{\cal L}^{\mathfrak{X}(m)\alpha}, \gamma;\K).
$$
\end{proposition}

Replacing   $g(\dot x(t), \dot x(t))$ by
 $F(\dot x(t), \dot x(t))$ in the proof of \cite[Lemma~3.3]{Tan82} we get

\begin{lemma}\label{lem:6.11}
 For $m=m^i_jk_j\ge\mathfrak{m}_0$ and $\mathfrak{X}(m)=s^i_jk_j$ as defined by Lemma~\ref{lem:6.8},
 $$
 \mathscr{H}^0_\ast\left(^{m\alpha+1}\!{\cal L}^{\mathfrak{X}(m)\alpha},\gamma;\K\right)
 \cong\mathscr{H}^0_\ast\left(^{\mathfrak{X}(m)\alpha}\!{\cal L}^{m\alpha+1},\gamma;\K\right).
 $$
\end{lemma}

 For each fixed integer $m=m^i_jk_j\ge \mathfrak{m}_0$ and
 $\mathfrak{X}({m}):=s^i_jk_j\le \mathfrak{s}\cdot\max_jk_j$ we have
 \begin{equation}\label{e:6.28}
\mathscr{H}^0_\ast({\cal L}, \gamma^{{m}\alpha+1};\K)\cong \mathscr{H}^0_\ast(\,^{{m}\alpha+1}\!{\cal L},\gamma;\K)
\cong\mathscr{H}^0_\ast\left(^{\mathfrak{X}({m})\alpha}\!{\cal L}^{{m}\alpha+1},\gamma;\K\right)
 \end{equation}
by (\ref{e:1.46}), Proposition~\ref{prop:6.10} and Lemma~\ref{lem:6.11}.
Let us write $m=l_m\mathfrak{X}(m)+ \eta(m)$ with
  $l_m\in\Z$ and $0\le\eta(m)\le\mathfrak{X}(m)\le \mathfrak{s}\max_jk_j$. Since $x\in\Lambda^{m\alpha+1}(M,\mathbb{I})\cap\Lambda^{\mathfrak{X}(m)\alpha}(M)$
  satisfies   
  $$
  \mathbb{I}(x(t))=x(t+m\alpha+1)=x(t+l_m\mathfrak{X}(m)\alpha+ \eta(m)\alpha+1)=x(t+ \eta(m)\alpha+1)
  \;\hbox{for any}\; t\in\R,
  $$ 
  $\Lambda^{m\alpha+1}(M,\mathbb{I})\cap\Lambda^{\mathfrak{X}(m)\alpha}(M)
   =\Lambda^{\eta(m)\alpha+1}(M,\mathbb{I})\cap\Lambda^{\mathfrak{X}(m)\alpha}(M)$
    as Hilbert manifolds and
    \begin{eqnarray*}
    ^{\mathfrak{X}(m)\alpha}\!{\cal L}^{m\alpha+1}=\, ^{\mathfrak{X}(m)\alpha}\!{\cal L}|_{\Lambda^{m\alpha+1}(M,\mathbb{I})\cap\Lambda^{\mathfrak{X}(m)\alpha}(M)}=\,
    ^{\mathfrak{X}(m)\alpha}\!{\cal L}|_{\Lambda^{\eta(m)\alpha+1}(M,\mathbb{I})\cap\Lambda^{\mathfrak{X}(m)\alpha}(M)}.
    \end{eqnarray*}
    Then $\eta(m)\le\mathfrak{X}(m)\le \mathfrak{s}\max_jk_j$ implies that the set
  \begin{eqnarray*}
 && \Bigl\{\mathscr{H}^0\left(^{\mathfrak{X}(m)\alpha}\!{\cal L}^{m\alpha+1},\gamma;\K\right)\,\Bigm|\,
 m\ge \mathfrak{m}_0\Bigr\}\\
 &=&\Bigl\{\mathscr{H}^0\bigl(\,^{\mathfrak{X}(m)\alpha}\!{\cal L}|_{\Lambda^{\eta(m)\alpha+1}(M,\mathbb{I})\cap\Lambda^{\mathfrak{X}(m)\alpha}(M)},
 \gamma;\K\bigr)\,\Bigm|\,m\ge \mathfrak{m}_0\Bigr\}
 \end{eqnarray*}
  is finite, and hence
$\{\dim\mathscr{H}^0_k({\cal L}, \gamma^{m\alpha+1};\K)\,|\,k\ge 0,\,m\ge 0\}$
is bounded by  (\ref{e:6.28}).


\begin{thebibliography}{L3}



\bibitem{AbF} A. Abbondandolo, A. Figalli,  High action orbits for
 Tonelli Lagrangians and superlinear Hamiltonians on compact configuration spaces.
 {\it J. Differential Equations}, {\bf 234}(2007), no. 2, 626--653.


\bibitem{ArKe} B. Aradi, D. CS. Kertesz,
Isometries, submetries and distance coordinates on Finsler manifolds.
{\it Acta Math. Hungar.}, DOI: 10.1007/s10474-013-0381-1.


\bibitem{BKl} V. Bangert, W. Klingenberg, Homology generated by iterated closed geodesics.
{\it Topology}, {\bf 22}(1983), 379-388.

\bibitem{BaChSh} D. Bao, S. S. Chern, Z. Shen, {\it An Introduction to
Riemann-Finsler Geometry}. Springer, Berlin (2000).




\bibitem{Bre} H. Brezis,  {\it Functional Analysis,  Sobolev Spaces and
Partial Differential Equations}. Springer 2011.

\bibitem{Bri} F. Brickell, On the differentiability of affine and projective transformations.
{\it Proc. Amer. Math. Soc.}, {\bf 16}(1965), 567-574.



\bibitem{Ca} E. Caponio, The index of a geodesic in a randers space
and some remarks about the lack of regularity of the energy
functional of a Finsler metric. {\it Acta Mathematica Academiae
Paedagogicae Ny\'iregyh\'aziensis}, {\bf 26} (2010), 265-274.


\bibitem{CaJaMa1} E. Caponio, M. A. Javaloyes and A.
Masiello, Morse theory of causal geodesics in a stationary spacetime
via Morse theory of geodesics of a Finsler metric. {\it Ann. Inst. H. Poincar¨¦ Anal. Non Lin\'eaire},
 {\bf 27}(2010), no. 3, pp. 857-876.

\bibitem{CaJaMa2} E. Caponio, M. A. Javaloyes and A. Masiello,  Addendum to ``Morse theory of
 causal geodesics in a
stationary spacetime via Morse theory of geodesics of a Finsler
metric". {\it Ann. Inst. H. Poincar\'e Anal. Non Lin\'eaire}, {\bf 30}(2013)
961-968.

\bibitem{CaJaMa3} E. Caponio, M. A. Javaloyes and A.
Masiello, On the energy functional on Finsler manifolds and
applications to stationary spacetimes. {\it Math.Ann.}, {\bf
351}(2011), no.2, 365-392.


\bibitem{Ch93} K. C. Chang, {\it Infinite Dimensional Morse Theory and Multiple Solution Problem}.
Birkh\"{a}user, 1993.


\bibitem{DenHo} S. Deng, Z. Hou, The group of isometries of a Finsler space,
{\it Pacific J. Math.}, {\bf 207}(2002), 149-155.

\bibitem{FHT01} Y. F\'elix,
S. Halperin,  and J.-C. Thomas,
{\it Rational Homotopy Theory}.  Graduate Texts in
Mathematics, 205. Springer-Verlag, New York, 2001.

\bibitem{FOT08} Y. F\'elix, J. Oprea, D. Tanr\'e,
 {\it Algebraic models in geometry}. Oxford Graduate Texts in Mathematics,
  17. Oxford University Press, Oxford, 2008.

\bibitem{Fr} T. Frankel, On theorems of Hurwitz and Bochner. {\it J. Math. Mech.}, {\bf 15}(1966), 373-377.

\bibitem{GM1} D. Gromoll, W. Meyer, On differentiable functions with isolated
critical points, {\it Topology}, {\bf 8}(1969), 361--369.

\bibitem{GM2} D. Gromoll, W. Meyer,  Periodic geodesics on compact
          Riemannian manifolds. {\it J. Diff. Geom.}, {\bf 3}(1969), 493-510.


\bibitem{Gro73I} K. Grove, Condition (C)  for the energy integral on certain path spaces and
applications to the theory of geodesics. {\it J. Diff. Geom.}, {\bf 8}(1973), 207-223.


\bibitem{Gro74} K. Grove,  Isometry-invariant geodesics, {\it Topology}, {\bf 13}(1974), 281-292.

\bibitem{Gro85} K. Grove,  The isometry-invariant geodesics problem: closed and
open. Geometry and topology (College Park, Md., 1983/84), 125-140,
Lecture Notes in Math., 1167, Springer, Berlin, 1985.

\bibitem{GroHaV78} K. Grove, S. Halperin; M.  Vigu\'e-Poirrier,  The rational homotopy theory
  of certain path spaces with applications to geodesics. {\it Acta Math.} {\bf 140} (1978), no. 3-4, 277-303.

\bibitem{GroTa78} K. Grove, M, Tanaka,  On the number of invariant closed
geodesics. {\it Acta Math.} {\bf 140 }(1978), no. 1-2, 33-48.

\bibitem{GroHa82} K. Grove, S. Halperin,  Contributions of rational homotopy
theory to global problems in geometry. {\it Inst. Hautes \'eudes
Sci. Publ. Math.} No. {\bf 56} (1982), 171-177 (1983).

\bibitem{GroHa91} K. Grove, S. Halperin, Elliptic isometrics, condition (C) and proper
maps. {\it Arch. Math.}, {\bf 56}(1991), 288-299.


\bibitem{Hin88} N.~Hingston,  Isometry-invariant geodesics on spheres, {\it Duke Math. J.},
  {\bf 57}(1988), no.~3, 761--768.



\bibitem{Hi} M. Hirsch, {\it Differential Topology}.
Springer-Verlag, 1976.

\bibitem{HrMa} U. Hryniewicz, L. Macarini, Local contact homology and applications.
arXiv:1202.3122.

\bibitem{JoMcC} J.D.S. Jones, J. McCleary,  String homology, and closed geodesics on manifolds which are elliptic spaces. arXiv:1409.8643.



\bibitem{Jav13}  M.~{\'A}. Javaloyes,  Chern connection of a pseudo-{F}insler metric
  as a family of affine connections. {\it Publ. Math. Debrecen}, {\bf 84}(2014), pp. 29-43.


 \bibitem{Jav13+}  M.~{\'A}. Javaloyes,
 Corrigendum to ``Chern connection of a pseudo-Finsler metric as a family of affine connections".
  {\it Publ. Math. Debrecen}, {\bf 84}(2014), pp. 481-487.


\bibitem{Jav14}  M.~{\'A}. Javaloyes, B.~L.~Soares,
Geodesics and Jacobi fields of pseudo-Finsler manifolds.
{\it Publ. Math. Debrecen}, {\bf 87}(2015), pp. 57-78.


\bibitem{Ker} J. Kern, Lagrange Geometry. {\it Arch. Math.}, {\bf 25}(1974), 438-443.


\bibitem{KoKrVa} L. Kozma, A. Krist\'aly, C. Varga, Critical point
theorems on Finsler manifolds. {\it Beitr\"age zur Algebra und}
({\it Geometrie Contributions to Algebra and Geometry}), {\bf
45}(2004), no.1, 47-59.

\bibitem{La} S. Lang, {\it Differential Manifolds}.
Springer-Verlag, 1985.


\bibitem{Lu1} G. Lu, Corrigendum: The Conley conjecture for Hamiltonian
systems on the cotangent bundle and its analogue for Lagrangian
systems. {\it J. Funct. Anal.} {\bf 261}(2011), 542-589.


\bibitem{Lu2} G. Lu, The splitting lemmas for nonsmooth functionals
on Hilbert spaces I. {\it Discrete and continuous dynamical
systems}, {\bf 33}(2013), no.7, 2939-2990.
 arXiv:1211.2127.


\bibitem{Lu3} G. Lu, Methods of infinite dimensional Morse theory
for geodesics on Finsler manifolds. {\it Nonlinear Analysis}, {\bf 113}(2015),
230-282, DOI:10.1016/j.na.2014.09.016. arXiv:1212.2078v7.


\bibitem{Lu4} G. Lu, Nonsmooth Morse theory and applications to higher elliptic equations,
 harmonic maps and minimal surfaces.  A book in preparation.


\bibitem{MaOu} J. Margalef-Roig, E. Outerelo Dominguez, {\it Differential Topology}.
                  North-Holland Mathematics Studies {\bf 173},
                  North-Holland-Amsterdam.

\bibitem{Ma} H. H. Matthias, {\it Zwei Verallgemeinerungen eines Satzes von Gromoll und Meyer}.
Bonn Mathematical Publications, Universit\"at Bonn Mathematisches
Institut, 1980.


\bibitem{MW} J. Mawhin, M. Willen, {\it Critical point theory and
Hamiltonian systems}. Applied Mathematical Science, Springer-Verlag,
New York, 1989.

\bibitem{Maz1} M. Mazzucchelli, On the multiplicity of isometry-invariant  geodesics on product manifolds.
{\it Algebraic \& Geometry Topology}, {\bf 14}(2014), 135-156.

\bibitem{Maz2} M. Mazzucchelli, Isometry-invariant  geodesics and the fundamental group.
{\it Math. Ann.}, {\bf  362}(2015),  265¨C280.


\bibitem{McCZ} J. McCleary and W. Ziller, On the free loop space of homogeneous spaces. {\it Amer.
J. Math.}, {\bf 109}(1987), 765-782. Corrections to: ``On the free loop space of homogeneous spaces"
{\it Amer. J. Math.}, {\bf 113}(1991), 375-377.

\bibitem{McSa} D. McDuff, D. Salamon, {\it $J$-holomorphic curves and
symplectic topology}. American Mathematical Society, Providence,
RI, 2004.

\bibitem{Mcl} M. Mclean, Local Floer homology and infinitely many
simple Reeb orbits. {\it Geom. Topl.}, {\bf 12}(2012), no.4, 1901-1923.

\bibitem{Me} F. Mercuri, The critical points theory for the closed
geodesics problem. {\it Math. Z.}, {\bf 156}(1977), 231-245.


\bibitem{MirA} R. Miron, M. Anastasiei, {\it The Geometry of Lagrange Spaces:
Theory and Applications}. Kluwer Academic Publishers, 1994.



\bibitem{Pa} R. S. Palais, Homotopy theory of infinite dimensional
manifolds. {\it Topology}, {\bf 5}(1966), 1-66.

\bibitem{Pa79} R. S. Palais, The Principle of Symmetric Criticality. {\it Commun. Math. Phys.},
{\bf 69}(1979), 19-30.

\bibitem{PaT} R. S. Palais, Chuu-lian Terng, {\it Critical Point Theory and Submanifold Geometry}.
Lecture Notes in Math., 1353, Springer, Berlin, 1988.

\bibitem{PP05} S. Papadima; L. Paunescu,  Isometry-invariant geodesics
and nonpositive derivations of the cohomology. {\it J. Differential
Geom.} {\bf 71}(2005), no. 1, 159-176.




\bibitem{PiTa} P. Piccione, D. Tausk, On the Banach differential
structure for sets of maps on non-compact domains. {\it Nonlinear Analysis},
{\bf 46}(2001), 245-265.




\bibitem{Rad89} H. B. Rademacher, Metrics with only finitely many isometry invariant
  geodesics, {\it Math. Ann.}, {\bf 284}(1989), no.~3, 391--407.



\bibitem{Ra042} H. B. Rademacher,  {\it Nonreversible Finsler metrics of positive flag curvature}, in A sampler of Riemann-Finsler geometry, vol. 50 of Math. Sci. Res. Inst. Publ., Cambridge Univ. Press,
Cambridge, 2004, pp. 261-302.



\bibitem{ShenW} Z. Shen, {\it Lectures on Finsler Geometry}.
      World Scientific Publishing Co., New Jersey (2001).



\bibitem{Tan77I} M. Tanaka, On invariant closed geodesics under isometries.
{\it Kodai Math. Sem. Rep.} {\bf 28} (1976/77), no. 2-3, 262-277.



\bibitem{Tan82} M. Tanaka, On the existence of infinitely many isometry-invariant geodesics.
 {\it J. Diff. Geom.}  {\bf 17}(1982), no.2, 171-184.


\bibitem{VPS76} M. Vigu\'e-Poirrier, D. Sullivan, The homology theory of closed geodesic problem.
{\it J. Diff. Geom.}, {\bf 11}(1976), 633-644.


\bibitem{Ze} E. Zeidler, {\it Nonlinear functional analysis and its applications. II/A. Nonlinear monotone operators.} Translated from the German by the author and Leo F. Boron. Springer-Verlag, New York, 1990.




\end{thebibliography}
\end{document}